%%%%%%%%%%%%%%%%%%%%%%%%%%%%%%%%%%%%%%%%%%%%%%%%%%%%%%%%%%
%%%%%%%%%%%%%%%%%%%%%%%%%%%%%%%%%%%%%%%%%%%%%%%%%%%%%%%%%%
%%
%%     This is the AMS-LaTeX file:
%%
%%     Colli-Signori-Sprekels
%%     Second-order analysis of an optimal control problem in a phase field 
%%     tumor growth model with singular potentials and chemotaxis
%%
%%%%%%%%%%%%%%%%%%%%%%%%%%%%%%%%%%%%%%%%%%%%%%%%%%%%%%%%%%
%%%%%%%%%%%%%%%%%%%%%%%%%%%%%%%%%%%%%%%%%%%%%%%%%%%%%%%%%%

\def\LaTeX{%
  \let\Begin\begin
  \let\End\end
  \let\salta\relax
  \let\finqui\relax
  \let\futuro\relax}

\def\UK{\def\our{our}\let\sz s}
\def\USA{\def\our{or}\let\sz z}

\UK
%\USA

%%%%%%%%%%%%%%%%%%%%%%%%%%%%%%%%%

% scegliere fra \TeX e \LaTeX  e fra  \UK oppure \USA

%\TeX
\LaTeX

%\UK
\USA

%%%%%%%%%%%%%%%%%%%%%%%%%%%%%%%%%
%% page layout
%%%%%%%%%%%%%%%%%%%%%%%%%%%%%%%%%

\salta

\documentclass[twoside,12pt]{article}
\setlength{\textheight}{24cm}
\setlength{\textwidth}{16cm}
\setlength{\oddsidemargin}{2mm}
\setlength{\evensidemargin}{2mm}
\setlength{\topmargin}{-15mm}
\parskip2mm

%%%%%%%%%%%%%%%%%%%%%%%%%%%%%%%%%
%% packages
%%%%%%%%%%%%%%%%%%%%%%%%%%%%%%%%%

%\usepackage{color}
\usepackage[usenames,dvipsnames]{color}
\usepackage{amsmath}
\usepackage{amsthm}
\usepackage{amssymb,bbm}
\usepackage[mathcal]{euscript}

\usepackage{cite}
\usepackage{hyperref}
\usepackage{enumitem}

\usepackage[ulem=normalem,draft]{changes}
%
%		COLORS FOR CORRECTIONS
%
% do the same, please (i.e., don't use the standard {\color{red} text} or similar): 
% just choose the color you prefer in \def\yourname

% EXAMPLE OF USE:  \fredi{I want this to become blue}
%
%IF YOU LATER WANT TO LET THE COLOR DISAPPEAR, ACTIVATE \def\fredi #1{{#1}} BELOW
 
\definecolor{viola}{rgb}{0.3,0,0.7}
\definecolor{ciclamino}{rgb}{0.5,0,0.5}
\definecolor{blu}{rgb}{0,0,0.7}
\definecolor{rosso}{rgb}{0.85,0,0}

\def\elvis #1{{\color{green}#1}}
\def\juerg #1{{\color{green}#1}}
\def\anold #1{{\color{magenta}#1}}

\def\pier #1{{\color{red}#1}} 
 
\def\juergen #1{{\color{blue}#1}}
\def\newju #1{{\color{red}#1}}
\def\pcol #1{{\color{blue}#1}}
\def\an #1{{\color{rosso}#1}}

\def\elvis #1{#1}
\def\juergen #1{#1}
\def\juerg #1{#1}
\def\anfirst #1{{#1}}
\def\andrea #1{{#1}}
\def\pier #1{{#1}}
\def\newju #1{{#1}}
\def\anold #1{{#1}}
\def\pcol #1{{#1}}
\def\an #1{{#1}}

%%%%%%%%%%%%%%%%%%%%%%%%%%%%%%%%%
%% you may adjust the baseline
%%%%%%%%%%%%%%%%%%%%%%%%%%%%%%%%%

%\renewcommand{\baselinestretch}{0.975}

%%%%%%%%%%%%%%%%%%%%%%%%%%%%%%%%%
%% bibliographystyle
%%%%%%%%%%%%%%%%%%%%%%%%%%%%%%%%%

\bibliographystyle{plain}

%%%%%%%%%%%%%%%%%%%%%%%%%%%%%%%%%
%% environments
%%%%%%%%%%%%%%%%%%%%%%%%%%%%%%%%%

%
\newtheorem{theorem}{Theorem}[section]

\newtheorem{definition}[theorem]{Definition}

\finqui

\def\Bthm{\Begin{theorem}}
\def\Ethm{\End{theorem}}
\def\Blem{\Begin{lemma}}
\def\Elem{\End{lemma}}

\def\Brem{\Begin{remark}\rm}
\def\Erem{\End{remark}}

\def\Bdim{\Begin{proof}}
\def\Edim{\End{proof}}
\def\Bcenter{\Begin{center}}
\def\Ecenter{\End{center}}
\let\non\nonumber

%%%%%%%%%%%%%%%%%%%%%%%%%%%%%%%%%
%% macros
%%%%%%%%%%%%%%%%%%%%%%%%%%%%%%%%%

% macro salvate

% sottosezioni non numerate

\def\step #1 \par{\medskip\noindent{\bf #1.}\quad}
\def\jstep #1: \par {\vspace{2mm}\noindent\underline{\sc #1 :}\par\nobreak\vspace{1mm}\noindent}

\def\Lip{Lip\-schitz}
\def\Holder{H\"older}
\def\Frechet{Fr\'echet}

\def\lhs{left-hand side}
\def\rhs{right-hand side}
\def\sfw{straightforward}

% versioni inglesi (UK) o americane (USA)

% bold, cal e mathop

\def\multibold #1{\def\arg{#1}%
  \ifx\arg\pto \let\next\relax
  \else
  \def\next{\expandafter
    \def\csname #1#1#1\endcsname{{\bf #1}}%
    \multibold}%
  \fi \next}

\def\pto{.}

\def\multical #1{\def\arg{#1}%
  \ifx\arg\pto \let\next\relax
  \else
  \def\next{\expandafter
    \def\csname cal#1\endcsname{{\cal #1}}%
    \multical}%
  \fi \next}

% operatori

\def\multimathop #1 {\def\arg{#1}%
  \ifx\arg\pto \let\next\relax
  \else
  \def\next{\expandafter
    \def\csname #1\endcsname{\mathop{\rm #1}\nolimits}%
    \multimathop}%
  \fi \next}

\multibold
qwertyuiopasdfghjklzxcvbnmQWERTYUIOPASDFGHJKLZXCVBNM.

\multical
QWERTYUIOPASDFGHJKLZXCVBNM.

\multimathop
diag dist div dom mean meas sign supp .

% accorpamenti di formule citate:
% uso  \accorpa {prima}{seconda}
%      \Accorpa\cs prima seconda (con il comodo blank anche dopo)
% NB: \Accorpa definisce \cs come l'accorpamento delle due citazioni
% e scrive sul file.log

\def\Accorpa #1#2 #3 {\gdef #1{\eqref{#2}--\eqref{#3}}%
  \wlog{}\wlog{\string #1 -> #2 - #3}\wlog{}}

% macro comode

\def\infess{\mathop{\rm inf\,ess}}
\def\supess{\mathop{\rm sup\,ess}}

\def\<#1>{\mathopen\langle #1\mathclose\rangle}
\def\norma #1{\mathopen \| #1\mathclose \|}

\def\I2 #1{\int_{Q_t}|{#1}|^2}
\def\IT2 #1{\int_{Q_t^T}|{#1}|^2}
\def\IO2 #1{\norma{{#1(t)}}^2}
\def\ov #1{{\overline{#1}}}
\def\next{\\ & \quad}

\def\iot {\int_0^t}

\def\intQt{\int_{Q_t}}
\def\intQ{\int_Q}
\def\iO{\int_\Omega}

\def\Qtt{\int_{Q_t^T}}

\def\dt{\partial_t}
\def\dn{\partial_{\bf n}}
\def\S{{\cal S}}

\def\X{{\cal X}}
\def\Y{{\cal Y}}
\def\Uh{{\cal U}}

\def\checkmmode #1{\relax\ifmmode\hbox{#1}\else{#1}\fi}

% insiemi numerici

\def\bu{{\bf u}}
\def\bh{{\bf h}}
\def\bk{{\bf k}}
\def\xih{{\xi^{\bf h}}}
\def\xik{{\xi^{\bf k}}}

\def\erre{{\mathbb{R}}}
\def\rz{{\mathbb{R}}}

\def\enne{{\mathbb{N}}}

\def\CC{{\mathbb{C}}}

\def\J{{\cal J}}
\def\Z{{\cal Z}}
\def\Jred{{\J}_{\rm red}}
\def\Pn{{\mathbb{P}_n}}

% spazi di funzioni a valori vettoriali su [0,T], [0,t], [0,s], [0,+\infty), [\delta,T]

% Come ricordare: in generale i simboli L H W  C da soli per gli spazi su (0,T)
% gli stessi raddoppiati per (0,+\infty)
% aggiunta di t o s al simbolo per (0,t) e (0,s)
% aggiunta di d al simbolo semplice o doppio per intervalli (\delta,T) e (\delta,+\infty)
% il simbolo C e i suoi derivati mettono le quadre anziche' le tonde

% Esempi   \L2V   \L\infty\Vp   \W{1,1}H   \C0H   \LL2V   \CC0\Vp   \Ld2V  \CCdH

\def\genspazio #1#2#3#4#5{#1^{#2}(#5,#4;#3)}
\def\spazio #1#2#3{\genspazio {#1}{#2}{#3}T0}

\def\L {\spazio L}
\def\H {\spazio H}

\def\C #1#2{C^{#1}([0,T];#2)}

%\def\LL {\spazioinf L}
%\def\HH {\spazioinf H}
%\def\WW {\spazioinf W}
%\def\CC #1#2{C^{#1}([0,+\infty);#2)}

% spazi di funzioni su \Omega, \Gamma, Q e \Sigma

\def\Lx #1{L^{#1}(\Omega)}
\def\Hx #1{H^{#1}(\Omega)}

\def\Ldue{\Lx 2}

\def\Huno{\Hx 1}
\def\Hdue{\Hx 2}

\def\Liq{{L^\infty(Q)}}

% spazi di funzioni su Q e S

% lettere greche

%\let\badtheta\theta
%\let\theta\vartheta
\let\eps\varepsilon
\let\vp\varphi
\let\lam\lambda

\def\a{\alpha}	%%%alpha
\def\b{\beta}	%%%beta
\def\d{\delta}  %%%delta
   %%%eta
\def\th{\theta} %%%theta
\def\r{\rho}    %%%rho
\def\s{\sigma}  %%%sigma
\def\m{\mu}	    %%%mu
\def\ph{\varphi}	%%%phi
\def\z{\zeta}      %%%%%zeta
      %%%%%chi
      %%%%%psi
\def\om{\omega}
\def\cd{c_{\d}}
\def\h{\mathbbm{h}}

\let\TeXchi\chi                         % new \chi, exactly on the baseline
\newbox\chibox
\setbox0 \hbox{\mathsurround0pt $\TeXchi$}
\setbox\chibox \hbox{\raise\dp0 \box 0 }
\def\chi{\copy\chibox}

% quadratino di fine dimostrazione

% abbreviazioni specifiche del lavoro

\def\ubar{\overline{\bf u}}
\def\uebar{\overline{u}_1}
\def\uzbar{\overline{u}_2}

\let\hat\widehat

\def\uad{{\cal U}_{\rm ad}}
\def\UR{{\cal U}_{R}}
\def\bmu{\overline\mu}
\def\bvp{\overline\varphi}
\def\bsigma{\overline\sigma}

\def\bph{{\ov \ph}}
\def\bm{{\ov \m}}   
\def\bs{{\ov \s}}
\def\diff{\bph^\bk-\bph}

\def\CC{{\mathscr{C}}}

\usepackage{amsmath}
\DeclareFontFamily{U}{mathc}{}
\DeclareFontShape{U}{mathc}{m}{it}%
{<->s*[1.03] mathc10}{}

\DeclareMathAlphabet{\mathscr}{U}{mathc}{m}{it}
%%%%%%%%%%%%%%%%%%%%%%%%%%%%%%
\Begin{document}
%%%%%%%%%%%%%%%%%%%%%%%%%%%%%%%%%

%%%%%%%%%%%%%%%%%%%%%%%%%%%%%%%%%
%% front page
%%%%%%%%%%%%%%%%%%%%%%%%%%%%%%%%%

%
\title{Second-order analysis of an optimal control problem in a phase field 
tumor growth model with singular potentials and chemotaxis}
\author{}
\date{}
\maketitle
\Bcenter
\vskip-1.5cm
{\large\sc Pierluigi Colli$^{(1)}$}\\
{\normalsize e-mail:{\tt pierluigi.colli@unipv.it}}\\[0.25cm]
{\large\sc Andrea Signori$^{(2)}$}\\
{\normalsize e-mail:{\tt andrea.signori02@universitadipavia.it}}\\[0.25cm]
{\large\sc J\"urgen Sprekels$^{(3)}$}\\
{\normalsize e-mail: {\tt juergen.sprekels@wias-berlin.de}}\\[.5cm]
$^{(1)}$
{\small Dipartimento di Matematica ``F. Casorati''}\\
{\small Universit\`a di Pavia}\\
{\small via Ferrata 5, I-27100 Pavia, Italy}\\[.3cm] 
$^{(2)}$
{\small Dipartimento di Matematica e Applicazioni}\\
{\small Universit\`{a} di Milano-Bicocca}\\
{\small via Cozzi, I-20125 Milan, Italy}\\[.3cm]
$^{(3)}$
{\small Department of Mathematics}\\
{\small Humboldt-Universit\"at zu Berlin}\\
{\small Unter den Linden 6, D-10099 Berlin, Germany}\\[2mm]
{\small and}\\[2mm]
{\small Weierstrass Institute for Applied Analysis and Stochastics}\\
{\small Mohrenstrasse 39, D-10117 Berlin, Germany}\\[10mm]
\Ecenter
\Begin{abstract}
\noindent 
This paper concerns a distributed optimal control problem for a tumor growth model of Cahn--Hilliard type including chemotaxis with possibly singular \an{potentials, where the} control and state variables are nonlinearly coupled.
First, we discuss the weak well-posedness of the system under very general assumptions for the potentials, which may be singular and \pcol{nonsmooth}. 
Then, we establish the strong well-posedness \pcol{of the system in a reduced setting,
which however admits the logarithmic potential: this analysis will 
lay the foundation for the study of} the corresponding optimal control problem.
Concerning the optimization problem, we address the existence of \anfirst{minimizers} and 
establish both first-order necessary and second-order sufficient conditions for optimality. 
\pcol{The mathematically challenging second-order analysis is completely performed here, after 
showing that the solution mapping is twice continuously differentiable between suitable 
Banach spaces via the implicit function theorem. Then, we completely identify the 
second-order Fr\'echet derivative of the control-to-state operator and carry out a thorough and 
detailed investigation about the related properties.}

\End{abstract}
\salta
\pagestyle{myheadings}
\newcommand\testopari{\sc Colli--Signori--Sprekels}
\newcommand\testodispari{\sc \an{Second-order analysis in a tumor growth control problem}}
\markboth{\testopari}{\testodispari}
\finqui
%
%%%%%%%%%%%%%%%%%%%%%%%%%%%%%%%%%
%% very beginning
%%%%%%%%%%%%%%%%%%%%%%%%%%%%%%%%%

\section{Introduction}
\label{INTRO}
\setcounter{equation}{0}

Lots of disclosures have been obtained in the \juerg{past} decades concerning tumor growth modeling: \anfirst{see, e.g., the pioneering works \cite{CLLW,CL,WLFC}.} 
The main advantage of \juerg{a} mathematical approach {is to be capable of predicting and analyzing} tumor growth behavior without inflicting any harm to the patients, thus helping medical practitioners to plan the clinical medications.

\juerg{The phase} field approach to tumor modeling consists in describing the tumor fraction by means of an order parameter $\ph$ representing the concentration of the tumor,
which \elvis{usually is normalized} to range between $-1$ and $1$.
Namely, the \anfirst{level} sets $\{\ph=1\}$ and $\{\ph=-1\}$ \pcol{may} describe the regions of pure phases\an{:} the tumorous phase and the healthy phase, respectively.
Moreover, the diffuse interface approach postulates the existence of a \juerg{thin} transition layer $\{-1< \ph< 1\}$ \anfirst{in which the phase variable passes rapidly, but continuously, from one phase to the other. We} assume the growth and proliferation of the tumor to be driven by the absorption and consumption of some nutrient, \anfirst{so that the equation for the phase variable, which has a Cahn--Hilliard type structure, is coupled with a reaction-diffusion equation for the}
variable $\s$ captur\anfirst{ing} the evolution of an unknown species nutrient (e.g., oxygen, glucose) in which the tissue 
\juerg{under} consideration is embedded.
%Thus, the model we are going to consider consists of a Cahn--Hilliard type equation with two relaxation terms an additional source term accounting for phase segregation
%coupled with a reaction-diffusion equation for the nutrient variable.

Let  $\a>0,~\b>0$, \juerg{and let} $\Omega\subset\erre^3$ denote some open and bounded domain having a smooth boundary $\Gamma=\partial\Omega$ and the unit outward normal $\,{\bf n}$.
We denote by $\dn$ the outward normal derivative to $\Gamma$.  Moreover, we fix some final time $T>0$ and
introduce for every $t\in (0,T]$ the sets $Q_t:=\Omega\times (0,t)$, $Q_t^T:=\Omega\times (t,T)$,
 and $\Sigma_t:=\Gamma\times (0,t)$,
 where we put, for the sake of brevity, $Q:=Q_T$ and $\Sigma:=\Sigma_T$.
We then consider the following optimal control problem: 

\vspace{3mm}\noindent
(${\cal CP}$) \quad Minimize the {\it cost functional}
\begin{align} 
	%\non	
	\J((\mu,\vp,\sigma),{\bf u})
   :=&\,  \frac{{\anfirst{b_1}}}2 \intQ |\vp-\widehat \vp_Q|^2
  + \frac{{\anfirst{b_2}}}2 \iO |\vp(T)-\widehat\vp_\Omega|^2
  + \frac{{\anfirst{b_0}}}2 \intQ |{\bf u}|^2 
 % \,+\,\kappa\,g({\bf u})
 % \\
  \label{cost} 
  % =& \J_1((\mu,\vp,\sigma),{\bf u}) + \kappa g(\bu)
\end{align} 
subject to  the {\it state system }
\begin{align}
\label{ss1}
&\alpha\dt\mu+\dt\ph-\Delta\mu=P(\ph)(\sigma+\chi(1-\ph)-\mu) - \h(\ph)u_1 \quad&&\mbox{in }\,Q\,,\\[1mm]
\label{ss2}
&\beta\dt\vp-\Delta\vp+F'(\vp)=\mu+\chi\,\sigma \quad&&\mbox{in }\,Q\,,\\[1mm]
\label{ss3}
&\dt\sigma-\Delta\sigma=-\chi\Delta\vp-P(\ph)(\sigma+\chi(1-\ph)-\mu)+u_2\quad&&\mbox{in }\,Q\,,\\[1mm]
\label{ss4}
&\dn \mu=\dn\vp=\dn\sigma=0 \quad&&\mbox{on }\,\Sigma\,,\\[1mm] 
\label{ss5}
&\mu(0)=\mu_0,\quad \vp(0)=\vp_0,\quad \sigma(0)=\sigma_0\,,\quad &&\mbox{in }\,\Omega\,,
\end{align}
\Accorpa\Statesys {ss1} {ss5}
and to the control constraint
\begin{equation}
\label{uad}
{\bf u}=(u_1,u_2)\in\uad\,.
\end{equation}
Here, the constants ${\anfirst{b_1}},{\anfirst{b_2}}$ are nonnegative, while $\,{\anfirst{b_0}}\,$  
is positive. Moreover,
 $\widehat\vp_Q$ and $\widehat \vp_\Omega$ are given target functions,
and the set of {\it admissible controls} $\uad$ is a nonempty, closed and convex subset of the control space
\begin{equation}
\label{defU}
\Uh:= L^\infty(Q)^2.
\end{equation}

The state system \Statesys\ constitutes a simplified and relaxed version of the four-species thermodynamically consistent model for tumor growth
originally proposed by Hawkins-Daruud et al.\ in \cite{HZO} that additionally includes chemotaxis effects.
Let us briefly review the role of the occurring symbols. The primary variables $\ph, \m$ 
and $\s$ denote the phase field, the associated chemical potential, and the nutrient concentration, respectively.
Furthermore, we stress that the additional term $\a\dt\m$ 
corresponds to a parabolic regularization for equation \eqref{ss1},
whereas the term $\b\dt\ph$ is the viscosity contribution \juerg{to the} Cahn--Hilliard equation.
The key idea behind these regularizations \juerg{originates} from the fact that 
their presence allows us to take into account more general potentials
that may be singular and possibly nonregular.
The nonlinearity $P$ denotes a proliferation function, whereas the positive constant $\chi$
represents the chemotactic sensitivity. Lastly, as a common feature of phase field models, $F$ is a nonlinearity which is assumed 
to possess a double-well shape. Typical examples are given by the regular, logarithmic,
and double obstacle potentials, which are defined, in this order, by
\begin{align}
\label{regpot}
&F_{\rm reg}(r)=\frac 14 \left(1-r^2\right)^2 \quad\mbox{for $\,r\in\rz$,}\\[1mm]
\label{logpot}
&F_{\rm log}(r)=\left\{\begin{array}{ll}
(1+r)\,\ln(1+r)+(1-r)\,\ln(1-r)-k_1 r^2\quad&\mbox{for $\,r\in (-1,1)$}\\
2\,\ln(2)-k_1\quad&\mbox{for $\,r\in\{-1,1\}$\,,}\\
+\infty\quad&\mbox{for $\,r\not\in[-1,1]$}
\end{array}\right.\\[1mm]
&F_{\rm obs}(r)=\left\{\begin{array}{ll}
k_2(1-r^2)\quad&\mbox{for $\,r\in [-1,1]$}\\
+\infty\quad&\mbox{for $\,r\not\in[-1,1]$}
\end{array}\right.,
\label{2obspot}
\end{align}
where $k>1$ and $k_2>0$ so that $F_{\rm log}$ and $F_{\rm obs}$ are nonconvex. 
Observe that $F_{\rm log}$ is very relevant in the 
applications, where $F'_{\rm log}(r)$ becomes unbounded as $r\to\pm 1$,
and that in the case of \eqref{2obspot} the second equation \eqref{ss2} has to be interpreted
 as a differential inclusion,
where $F'(\ph)$ \elvis{is} understood in  the sense of subdifferentials.

In this paper, we take two distributed controls that act in the phase equation and in the nutrient equation,
respectively. The control variable $u_1$, \newju{which is nonlinearly coupled to the state variable $\varphi$} \pcol{in the phase equation \eqref{ss1}, models} the application of a cytotoxic drug into 
the system; it is multiplied by a truncation function $\h(\anfirst{\cdot})$ in order to have the action
only in the spatial region where the tumor cells are \an{located. For instance,
it can be assumed that $\h(-1)=0, \h(1)=1 ,\h(\ph)$ is in between if $-1<\ph<1$;
% and $\h=0$ elsewhere; 
see \cite{GLSS, GARL_1, HKNZ, KL} for some insights on possible choices of $\h$.}
On the other hand, the control $u_2$ can model  
either an external medication or some nutrient supply.

As far as well-posedness \anfirst{is concerned}, the above model has already been investigated in the case $\chi=0$ in \cite{CGH,CGRS1,CGRS2,CRW},
and in \cite{FGR} with $\a=\b=\chi=0$.
There the authors also pointed out how $\a$ and $\b$ can be set to zero, by providing
the proper framework in which a limit system can be identified and uniquely solved.
We also note that in \cite{CGS24} a version has been studied in which the Laplacian in the equations 
\eqref{ss1}--\eqref{ss3} has been replaced by fractional powers of a more general class of selfadjoint operators
having compact resolvents. 

For some nonlocal variations of the above model we refer to \cite{FLRS, FLS, SS}.
Moreover, in order to better emulate in-vivo tumor growth,
it is possible to include in similar models the effects generated by the fluid flow development
by postulating a Darcy\anfirst{'s} law or a Stokes--Brinkman\anfirst{'s} law.
In this direction\anfirst{,} we refer to \cite{WLFC,GLSS,DFRGM ,GARL_1,GARL_4,GAR, EGAR, FLRS,GARL_2, GARL_3}, and \anfirst{we also mention} \cite{GLS}\anfirst{,} where elastic effects are included.
For further models, discussing the case of multispecies,
we address the reader to \cite{DFRGM,FLRS}.

The investigation of the associated optimal control problem also \pcol{presents} a 
wide number of results of which we mention\pcol{\cite{CGRS3,EK, EK_ADV, S, S_a, S_b, S_DQ, SigTime, FLS, ST, CGS24, CSS1, GARLR, KL, SW}.}
Notice that, despite the \juerg{number of contributions, only 
\cite{EK_ADV}} established second-order optimality conditions under suitable restrictions on the considered model.
In particular, the authors \pcol{of \cite{EK_ADV}} avoid considering the chemotaxis effects and allow only regular potentials to be considered.

\pcol{In this paper, first we discuss the weak well-posedness of the system \eqref{ss1}--\eqref{ss5} in a very general framework for the potentials, which includes 
all \elvis{of} the cases in \eqref{regpot}--\eqref{2obspot}. Then, we turn our attention to the strong well-posedness 
of~\eqref{ss1}--\eqref{ss5} in the cases of the regular $F_{\rm reg}$ and logarithmic $F_{\rm log}$ potentials. This is done in Section~\ref{STATE}, while the corresponding optimal control problem is investigated in the following sections. Section~\ref{MINIMIZER} 
is concerned with the existence of \anfirst{minimizers}, then the intensive and 
crucial~Section~\ref{SEC:DIFF} establishes the differentiability properties of the control-to-state 
operator and contains a number of results on the concerned linearized problems and 
the basic stability estimates for the solutions.  
The last two Sections~\ref{SEC:FOC} and~\ref{SEC:SOC} treat in some detail the 
first-order necessary and second-order sufficient conditions for optimality, respectively.
\anfirst{L}et us point out that the second-order analysis is 
challenging from the mathematical viewpoint and \juerg{demands} to prove 
that the solution mapping is twice continuously differentiable between suitable Banach spaces. By taking advantage of the regularizing effect due to the aforementioned relaxation terms, we can deal with a complete study of the second-order analysis, 
still covering the case of singular potentials and chemotaxis. Moreover, we are able to  identify the second-order Fr\'echet derivative of the control-to-state operator and 
investigate the related properties in a sharp and profound way.}  
 
Throughout the paper, we make repeated use of H\"older's inequality, of the elementary Young's inequality
\begin{equation}
\label{Young}
a b\,\le \delta |a|^2+\frac 1{4\delta}|b|^2\quad\forall\,a,b\in\erre, \quad\forall\,\delta>0,
\end{equation}
as well as the continuity of the embeddings $H^1(\Omega)\subset L^p(\Omega)$ for $1\le p\le 6$ and 
$\Hdue\subset C^0(\overline\Omega)$. Notice that the latter embedding is also compact, while this holds true
for the former embeddings only if $p<6$. 

Lastly, let us introduce a convention that will be tacitly employed in the rest of the paper:
the symbol small-case $c$ is used to indicate every constant
that depends only on the structural data of the problem (such as
$T$, $\Omega$, $\a$ or $\b$, the shape of the 
nonlinearities, and the norms of the involved functions),
so that its meaning may change from line to line.
\pier{When a parameter $\delta$ enters the computation, then
the symbol $c_\delta$ denotes constants
that depend on~$\delta$ in addition.
On the contrary, precise constants we could refer 
to are treated in a different way.}

%%%%%%%%%%%%%%%%%%%%%%%%%%%%%%%%%%%%%%%%%%%%%%%%%%%%%%%%%%%%%%%%%%%%%%

\section{General Setting and Properties of the State System}
\label{STATE}
\setcounter{equation}{0}

In this section, we introduce the general setting of our control 
problem and state some results on the state system \eqref{ss1}--\eqref{ss5}. 
To begin with, for a Banach space $\,X\,$ we denote by $\|\cdot\|_X$
the norm in the space $X$ or in a power \juerg{thereof}, and by $\,X^*\,$ its dual space. 
The only \anfirst{exeption} from this rule applies to the norms of the
$\,L^p\,$ spaces and of their powers, which we often denote by $\|\cdot\|_p$, for
$\,1\le p\le +\infty$. As usual, for Banach spaces $\,X\,$ and $\,Y\,$ we introduce the linear space
$\,X\cap Y\,$ which becomes a Banach space when equipped with its natural norm $\,\|u\|_{X \cap Y}:=
\|u\|_X\,+\,\|u\|_Y$, for $\,u\in X\cap Y$.
Moreover, we recall the definition~\eqref{defU} of~${\cal U}$ and introduce the spaces
\begin{align}
  & H := \Ldue \,, \quad  
  V := \Huno\,,   \quad
  W_{0} := \{v\in\Hdue: \ \dn v=0 \,\mbox{ on $\,\Gamma$}\}.
  \label{defHVW}
\end{align}
\anfirst{Furthermore, by} $\,(\,\cdot\,,\,\cdot\,)$, $\,\Vert\,\cdot\,\Vert$\anfirst{, and $\<\cdot,\cdot>$}, we denote the standard inner product 
\pier{and related norm} in $\,H$, \anfirst{as well as the dual product between $V$ and its dual $V^*$}. For given final time $T>0$, we introduce the spaces
\begin{align}\label{defZ}
&\newju{Z := \H1 H \cap \L\infty V \cap \L2 {W_0}, \quad \Z := Z \times Z \times Z},\\[0.5mm]
\label{defV}
&\newju{\mathcal{V}:=\bigl(L^\infty(0,T;H)\cap L^2(0,T;V)\bigr)\,\times\,Z\,\times\,
\bigl(L^\infty(0,T;H)\cap L^2(0,T;V)\bigr),}
\end{align}
\newju{which are Banach spaces when endowed with their natural norms.}

\vspace{2mm}
\pier{Some assumptions on the data are stated here.} 
\begin{enumerate}[label={\bf (W\arabic{*})}, ref={\bf (W\arabic{*})}]
\item \label{const:weak}
	$\alpha,\beta $ and $\chi$ are positive constants.
\item \label{F:weak}
	$F=F_1+F_2$, where $\,F_1:\erre\to [0,+\infty]\,$ is convex and lower semicontinuous with
	$\,F_1(0)=0$\pier{, and where $F_2 \in C^1(\erre)$ has a \Lip\ continuous derivative} $F'_2$.
\item \label{P:weak}
	$P \pier{{}\in C^0(\erre)}$ is nonnegative, bounded, and \Lip\ continuous.
\item \label{h:weak}
	$\h \pier{{}\in C^0(\erre)}$ is nonnegative, bounded, and \Lip\ continuous.
%\item \label{u:weak}
%	$ \pier{(u_1, u_2) \in L^\infty (Q) \times L^2(Q)}$.
\end{enumerate}
\pier{For the sake of} simplicity, we indicate with a common notation\pier{%
\begin{equation} 
\text{$L$ as a \Lip\ constant for $F_2', \, P, $ \an{\,and\,} $ \h$.}   \label{pier1} 
\end{equation}
Let us note that all of the choices \eqref{regpot}--\eqref{2obspot} are admitted for the potentials. 
In fact, the assumption \ref{F:weak} implies that the subdifferential $\partial F_1$ of $F_1$ is a 
maximal monotone graph in $\erre \times\erre$ with effective domain $D(\partial F_1 ) \subset D(F_1 ) $, and, since $F_1$
\juerg{attains the minimum value~$0$ at $0$, it} turns out that $0\in D(\partial F_1 )$ and $0\in\partial F_1(0)$. 
Now, in the general setting depicted by~\ref{const:weak}--\anfirst{\ref{h:weak}}}, we are able 
to provide a first well-posedness result for \pier{the} system~\Statesys. 
\pier{First}, let us present the notion of weak solution to \Statesys.
\begin{definition}
\label{DEF:WEAK}
A quadruplet $(\m,\ph,\xi,\s)$ \juerg{is called} a weak solution to the initial boundary value problem \Statesys\ if
\begin{align}
	\ph & \in \H1 H \cap \L\infty V \cap \L2 {W_0}, \label{pier2-1}
	\\
	\m,\s & \in \H1 {V^*} \cap \L\infty H \cap \L2 V, \label{pier2-2}
	\\ 
	\xi & \in \L2 H, \label{pier2-3}
\end{align}
and if $(\m,\ph,\xi,\s)$ satisfies the corresponding weak formulation \juerg{given by}
\begin{align}
	 \label{var:1}& \<\dt(\alpha\mu + \ph), v > 
	+ \iO \nabla \mu \cdot \nabla v
	= \iO P(\ph)(\sigma+\chi(1-\ph)-\mu)v
	-\iO \h(\ph)u_1 v
	\nonumber\\
	& \qquad \pier{\hbox{for every $v \in V $ and a.e. in $(0,T)$,}}\\[2mm]
	%	& \< \beta\dt\ph , v>
%	 +\iO \nabla \ph \cdot \nabla v 
%	+ \iO \xi v
%	+ \iO F_2(\vp)v=
%	\iO \mu v 
%	+\chi\iO \sigma v,
	\label{var:2}&\beta\dt\vp-\Delta\vp+\xi+F_2'(\vp)=\mu+\chi\,\sigma, \quad \hbox{$\xi \in \partial F_1(\ph)$, \, a.e. in $Q$,}
 	\\
 	\label{var:3}
 	& \<	\dt\sigma, v>
 	+ \iO \nabla \sigma\cdot \nabla v
 	=\chi\iO \nabla \vp\cdot \nabla v 
 	- \iO P(\ph)(\sigma+\chi(1-\ph)-\mu)v
 	+ \iO u_2 v	\nonumber\\
	& \qquad \pier{\hbox{for every $v \in V $ and a.e. in $(0,T)$,}}
\end{align}
as well as 
\begin{align}
	\label{var:4}
	\m(0)=\m_0, \quad
	\ph(0)=\ph_0, \quad
	\s(0)=\s_0, \quad \hbox{a.e. in $\Omega$}.
\end{align}
\end{definition}
It is worth noticing that the \pier{homogeneous Neumann boundary conditions \eqref{ss4}
are considered in the condition \eqref{pier2-1} for $\ph$ (cf. the definition of the space $W_0$) and 
incorporated in the variational equalities \eqref{var:1} and \eqref{var:3} for $\mu$ and $\sigma$, 
when using the forms $\iO \nabla \mu \cdot \nabla v  $ and $\iO \nabla \sigma \cdot \nabla v $. 
Moreover, let us point out that, at this level, the control pair} $(u_1,u_2)$ just \pier{yields two fixed 
forcing terms in \eqref{var:1} and \eqref{var:3}. The initial conditions~\eqref{var:4} make sense 
since  \eqref{pier2-1} and \eqref{pier2-2} ensure that $\ph$ and $\mu, \, \sigma$ are continuous 
from $[0,T] $ to $V$ and $H$, respectively.} 

\begin{theorem}[Weak well-posedness]
\label{THM:WEAK}
\pier{Assume that \ref{const:weak}--\anfirst{\ref{h:weak}}} hold.
Moreover, let the initial data $(\m_0,\ph_0,\s_0)$ satisfy
\begin{align}
	\label{weak:initialdata}
	\m_0, \s_0 \in \Lx2, 
	\quad
	\ph_0 \in \Hx1,
	\quad
	F_1(\ph_0) \in \Lx1,
\end{align}
\anfirst{and \juerg{suppose that the} source terms $u_1, u_2$ are such that
\begin{align}
	\label{u:weak}
		(u_1, u_2) \in L^2(Q) \times L^2(Q).
\end{align}}%
Then there exists \pcol{at least a} solution $(\m,\ph,\xi,\s)$ in the sense of Definition \ref{DEF:WEAK}.
Moreover, \pcol{if $u_1 \in L^\infty(Q)$ in addition to \eqref{u:weak}, \elvis{then} the found solution is unique.
Furthermore,} let $(\m_i,\ph_i,\xi_i,\s_i), $ \pier{$i=1,2$,} be two weak solutions to \Statesys\ associated 
with the initial data \pier{$(\m_0^i,\ph_0^i,\s_0^i)$, which
satisfy \eqref{weak:initialdata}, and controls $(u_1^i,u_2^i) \in L^\infty (Q) \times L^2(Q)$},
$i= 1,2 $.Then there is a \anfirst{positive} constant $C_{d}$, depending only on \an{the} data \an{of the system}, such that
\begin{align}
	& \non
	\norma{\a (\m_1-\m_2) + (\ph_1-\ph_2)}_{\L\infty H}
	+ \norma{\nabla(\m_1-\m_2)}	_{\L2 H}
		\\
	& \qquad \non
	+ \norma{\ph_1-\ph_2}_{\L\infty H \cap \L2 V}
	+ \norma{\s_1-\s_2}_{\L\infty H \cap \L2 V}
	\\& \quad \non
	\leq
	C_{d} \Big(\an{
	\norma{\alpha(\m_0^1-\m_0^2) + (\ph_0^1-\ph_0^2)}
	+ \norma{\ph_0^1-\ph_0^2}
	+ \norma{\s_0^1-\s_0^2}}
	\Big)
	\\& \qquad 
	+
	C_{d} \Big(
	 \norma{u_1^1-u_1^2}_{\L2 H}
	+  \norma{u_2^1-u_2^2}_{\L2 H}
	\Big).
	\label{cont:dep:weak}
\end{align} 
\end{theorem}
Before entering the proof, let us remark that the above result is very general and includes also the cases of singular and 
\pier{nonsmooth} potentials, such as the double obstacle potential defined by \eqref{2obspot}.
For the dependencies of the constant $C_{d}$, we invite the reader to follow the proof of 
the estimate~\eqref{cont:dep:weak} below. 

\begin{proof} \pier{For the proof of the existence of a solution, we point out that the arguments are quite  
standard, since similar procedures have already been used in previous contributions. 
\an{Thus, for that part, we} proceed rather formally, just employing 
the Yosida approximation of $\partial F_1$ for our estimates, without recurring to \juerg{finite-dimensional
approximation techniques} like the Faedo--Galerkin scheme.}

Hence, we introduce the Yosida regularization of $ \partial F_1$. 
For $\eps>0$ let $F_{1,\eps}$ denote the Moreau--Yosida
approximation of $F_1$ at the level $\eps$. 
It is well known (see, e.g., \cite{Brezis}) 
that the following conditions are satisfied:
\begin{align}\label{pier3}
&0\leq F_{1,\eps}(r)\le F_1(r) \quad\mbox{for all }\,r\in\erre. \\
&F'_{1,\eps} \mbox{ is Lipschitz continuous on $\,\erre$ \anfirst{with \Lip\ constant $\tfrac 1\eps$,} and $F'_{1,\eps}(0)=0$.}\label{pier4}\\
&|F'_{1,\eps}(r)|\le |(\partial F_1)^\circ (r)| \,\,\mbox{ and }\,\,\lim_{\eps\searrow 0}\,F'_{1,\eps}(r)=(\partial F_1)^\circ (r),
\quad\mbox{for all }\,r\in D(\partial F_1).  \label{pier5} 
\end{align} 
Here,  $(\partial F_1)^\circ $ denotes the minimal section of $\partial F_1$, that is, 
$(\partial F_1)^\circ (r)$ defines the element of $(\partial F_1) (r)$ with minimal \newju{modulus}.

Next, we are going to prove a series of estimates for the solution to problem \eqref{var:1}--\eqref{var:3}, 
where 
$(\partial F_1) (r)$ is replaced by $F'_{1,\eps}$ and the inclusion in  \eqref{var:2} reduces to an equality. Namely, we argue on
\begin{equation}
\label{pier6}\beta\dt\vp-\Delta\vp+F'_{1,\eps} (\vp) +F_2'(\vp)=\mu+\chi\,\sigma \quad \hbox{a.e. in $Q$.}
\end{equation}
For the sake of simplicity, we still denote by $(\m,\ph,\xi,\s)$, with $\xi = F'_{1,\eps} (\vp)$,
the solution to the approximated system in place of $(\m_\eps,\ph_\eps,\xi_\eps,\s_\eps)$; the correct notation will be reintroduced at the end of each estimate.

\noindent
{\sc First Estimate:}
{We add the term $\ph$ to both sides of \eqref{pier6} and test by $\dt \ph$. Then, we take $v=\mu$ in \eqref{var:1}  and $v=\sigma $ in \eqref{var:3}. Moreover, we add the resulting equalities and, with the help of a cancellation, we deduce that
almost everywhere in $(0,T)$ it holds the identity
\begin{align*}
	& \frac 12 \frac d{dt} \Big( \a \norma{\m}^2 + \norma{\ph}_V^2 + 2 \iO F_{1,\eps}(\ph) + \norma{\s}^2  \Big)
 	\\ 
	& \quad\quad
		+  \norma{\nabla \m}^2 
	+ \b \norma{\dt\ph}^2 
	+  \norma{\nabla \s}^2 
	+ \iO P(\ph)(\m-\s)^2 
	\\ 
	& \quad
	\pcol{{}={}}
	 \iO \chi P(\ph)(1-\ph)(\m-\s)
	- \iO \h(\ph)u_1 \m
	+ \iO \chi \s \dt \ph
	\\ & \qquad
	+ \iO ( \ph -F'_2(\ph)) \dt \ph
	+ \chi\iO \nabla \ph \cdot \nabla \s
	+ \iO u_2 \s =:I.
\end{align*}
Note that the last term on the \lhs\ is nonnegative due to \ref{P:weak}. Then, we can integrate the above inequality 
\juerg{over $[0,t]$ for} $t\in (0,T]$, using the initial conditions \eqref{var:4}. We point \newju{out} that the quantity 
$$ \a \norma{\m_0}^2 + \norma{\ph_0}_V^2 + 2 \iO F_{1,\eps}(\ph_0) + \norma{\s_0}^2  \quad 
\hbox{is bounded independently of } \eps, $$
thanks to \eqref{weak:initialdata} and \eqref{pier3}. Next, owing to the boundedness and regularity properties of $P$, $\h$ and $F_2'$,  and by Young's inequality, it is \sfw\ to infer that
\begin{align*}
 \int_0^t I \juerg{(s)\,ds} & 
	\leq \chi \norma{P}_\infty \intQt	(1+|\ph|^2 +|\m|^2+|\s|^2)
	+ \frac12 \norma{\h}_\infty 	\intQt	(|u_1|^2 +|\m|^2) 
	    \\ & \quad	
	+	\frac\b 4 \intQt |\dt\ph|^2  + \frac{\chi^2 }\b \intQt |\s|^2				
			+	\frac\b 4 \intQt |\dt\ph|^2 + c 	\intQt	(1+|\ph|^2 )
			    \\ & \quad	
			+ \frac1 2 	 \intQt |\nabla\s|^2 		
			+ \frac{\chi^2} 2 	 \intQt |\juerg{\nabla\varphi}|^2 					
	+ \frac12 \intQt	(|u_2|^2 +|\s|^2) 								
    \\ & 										
     \leq \iot  \Big( \frac\b 2 \norma{\dt\ph\juerg{(s)}}^2 
	+  \frac12 \norma{\nabla \s \juerg{(s)}}^2 \Big)\juerg{\,ds\,+\,c}\\
	&\quad \juerg{+\, c \iot \bigl(\norma{\ph(s)}_V^2 + \norma{\m(s)}^2+ \norma{\s(s)}^2 + \norma{u_1(s)}^2 + \norma{u_2(s)}^2
	\bigr  )\,ds},
\end{align*}
so that it suffices to apply Gronwall's lemma to conclude that
\begin{align}
	& \norma{\anfirst{\mu_\eps}}_{\L\infty H \cap \L2 V}
	+ \norma{\anfirst{\ph_\eps}}_{\H1 H \cap \L\infty V}
		\nonumber \\ & \quad
	+ \norma{\anfirst{F_{1,\eps}(\ph_\eps})}_{\L\infty {\Lx1}}^{1/2}
	+ \norma{\anfirst{\s_\eps}}_{\L\infty H \cap \L2 V}
	\leq c. \label{pier7}
\end{align}}%

\noindent 
\pier{{\sc Second Estimate:}
Now, owing to \eqref{pier7}, from a comparison of terms in the variational equalities \eqref{var:1} and \eqref{var:3} it follows that
\begin{align}
	\norma{\anfirst{\dt\m_\eps}}_{\L2 {V^*}}
	+ \norma{\anfirst{\dt\s_\eps}}_{\L2 {V^*}}
	\leq c. \label{pier8}
\end{align}
In fact,  arguing for instance on \eqref{var:1}, and taking advantage of \ref{P:weak},\ref{h:weak}, \pcol{\eqref{u:weak} and \eqref{pier7},} \juerg{we have that \pcol{for a.e. $t\in (0,T)$} and for every $v\in V$ it holds}
\begin{align*}
	\Big| \<\alpha \dt\mu \pcol{{}(t){}} , v > \Big|
	&{}\leq
	  \norma{\dt \ph \pcol{{}(t){}}}\,\norma{ v}
	+ \norma{ \nabla \mu \pcol{{}(t){}} }\,\norma{ \nabla v}
	\\ 
	&\quad + c\, 	(\norma{\s\pcol{{}(t){}}}+\norma{\ph\pcol{{}(t){}}}+ 1 
	+\norma{\m\pcol{{}(t){}}})\,\norma{ v}
	+ c \,\norma{u_1\pcol{{}(t){}}}\norma{ v}
	\\ 
	&\pcol{{}\leq 
	c \,( \norma{\dt \ph \pcol{{}(t){}}}+ \norma{u_1\pcol{{}(t){}}}+ 1)\,  \norma{ v}_V}.
\end{align*}
Thus, dividing by $\norma{ v}_V$ and passing to the superior limit, we readily have 
the bound for \pcol{$\norma{\dt\m (t)}_{V^*}$ in terms of $c ( \norma{\dt \ph \pcol{{}(t){}}}+ \norma{u_1\pcol{{}(t){}}}+ 1)$. Then, by squaring and integrating over $(0,T)$, we deduce \eqref{pier8} for $\dt\m$. The corresponding property} for $\dt\s$ can be obtained in a similar way from~\eqref{var:3}.}

\noindent 
{\sc Third Estimate:}
We rewrite \pier{\eqref{pier6} as 
\begin{align}
	\label{ph:elliptic}
	- \Delta \ph + F'_{1,\eps}(\ph) = \m + \chi \s - \b \dt\ph -F'_2(\ph)=: f_\ph \quad 
\end{align}
almost everywhere. Due to \eqref{pier7} and the Lipschitz continuity of $F'_2$, we infer that $f_\ph$ is uniformly bounded in $ \L2 H $.
Hence, we can test \eqref{ph:elliptic} by $F'_{1,\eps}(\ph)$ and integrate by parts in the first term, taking advantage of the homegeneous Neumann boundary condition and obtaining a nonnegative contribution. Thus, by a standard computation it turns out that   $F'_{1,\eps}(\ph)$ is bounded in $ \L2 H $ independently of $\eps$. Then, by comparison in \eqref{ph:elliptic} and thanks to the elliptic regularity theory, we finally derive that}
\begin{align}
	 \pier{\anfirst{\norma{F'_{1,\eps}(\ph_\eps)}}_{\L2 H} + \norma{\anfirst{\ph_\eps}}_{\L2 {W_0}}  \leq c.} \label{pier9}
\end{align}

\noindent 
\pier{{\sc Passage to the limit:}
Denote now by $(\m_\eps,\ph_\eps,\s_\eps)$ the triplet solving the problem \eqref{var:1}, \eqref{pier6}, \eqref{var:3}, \eqref{var:4} with the regularity \eqref{pier2-1}, \eqref{pier2-2}. Then, in view of the estimates \eqref{pier7}, \eqref{pier8}, \eqref{pier9}, \anfirst{which are independent of $\eps$,} by weak and weak-star compactness it turns out that there are $\m,\ph,\s $ and $\xi$ such that 
\begin{align}
& \m_\eps \to \m \quad \hbox{weakly star in } \   \H1 {V^*} \cap \L\infty H \cap \L2 V, \label{pier10-1}
\\
& \ph_\eps \to \ph \quad \hbox{weakly star in } \   \H1 {H} \cap \L\infty V \cap \L2 {W_0},  \label{pier10-2}
\\
&\s_\eps \to \s \quad \hbox{weakly star in } \   \H1 {V^*} \cap \L\infty H \cap \L2 V,  \label{pier10-3}
\\
& F'_{1,\eps}(\ph_\eps) \to \xi \quad \hbox{weakly in } \   \L2 {H}, \label{pier10-4}
\end{align} 
as $\eps \searrow 0$, possibly along a subsequence. By virtue of \eqref{pier10-1}--\eqref{pier10-3} and 
the Aubin--Lions lemma (see, e.g., \cite[Sect.~8, Cor.~4]{Simon}), we deduce \newju{that $\m_\eps\to\m$,  
$\ph_\eps\to \ph$, $\s_\eps\to \s$, all strongly} in $\L2 H$. Then, we can pass to the limit in the variational equalities  \eqref{var:1}, \eqref{var:3} 
and also in \eqref{pier6}, in order to obtain the equality in \eqref{var:2}. The nonlinearities $P(\ph_\eps)$,  $\h(\ph_\eps)$, $F'_{2}(\ph_\eps)$ can be 
easily taken to the limit, because of the Lipschitz continuity of the involved functions and of the strong 
convergence of $\ph_\eps $ to $\ph$ in $\L2 H$. In addition, the inclusion in \eqref{var:2} 
\juerg{results} as a consequence of \eqref{pier10-4} and the
maximal monotonicity of $\partial F_1$, \pcol{since} we can apply, e.g., \cite[Lemma~2.3, p.~38]{Barbu}. Finally, the initial conditions 
\eqref{var:4} can be readily obtained by observing that \eqref{pier10-1}--\eqref{pier10-3} imply weak convergence in $\C 0 H$ 
(actually, even strong for $\ph_\eps $ to $ \ph$).}
%\Brem
%\label{soloL2}
%Notice that for the existence proof we just used the regularity $L^2(Q) $ of the component $u_1$ of the control pair
%$(u_1, u_2)$ specified in \anfirst{\ref{h:weak}}. Hence, the existence of a weak solution can be derived whenever $(u_1, u_2)\in (L^2(Q))^2.$
%\Erem%

\pier{As for uniqueness, it 
suffices to show that \eqref{cont:dep:weak} is fulfilled for weak solutions. In fact, if 
we let $(\m_i,\ph_i,\xi_i,\s_i)$, $i=1,2$, denote two different weak solutions to 
\Statesys\ associated with the same initial data $(\m_0,\ph_0,\s_0)$ and control variables $
(u_1,u_2)$\juerg{, then} we derive that \eqref{cont:dep:weak} 
holds with the \rhs\ equal to zero,  so that $\ph_1=\ph_2$, $\m_1=\m_2$, $\xi_1=\xi_2$, $
\s_1=\s_2$, \juerg{whence the} uniqueness follows.}

\noindent 
\pier{{\sc Continuous dependence estimate:} Now, recalling the notation in} the statement of the theorem, we set, for $i= 1,2 $,
\begin{align*}
	\m&\anfirst{:= \m_1-\m_2, \quad
	\ph:= \ph_1-\ph_2, \quad}
	\xi:= \xi_1-\xi_2, \quad
	\s:= \s_1-\s_2, \quad
	\\ 
	\m_0 &:= \m_0^1-\m_0^2, \quad
	\ph_0 := \ph_0^1-\ph_0^2, \quad
	\s_0 := \s_0^1-\s_0^2, \quad
	\\ 
	u_i &:= u_i^1-u_i^2, \quad
	R_i := P(\ph_i)(\s_i + \chi (1-\ph_i) - \mu_i),
	\quad
	\h_i := \h(\ph_i),
\end{align*}
and consider the difference of the equations \pier{in \eqref{var:1}--\eqref{var:4} 
to infer that
\begin{align}
	 \label{cd:1}& \<\dt(\alpha\mu + \ph), v > 
	+ \iO \nabla \mu \cdot \nabla v
	= \iO \hat R \, v  - \iO ((\h_1-\h_2) u_1^1 + \h_2 u_1)v
	\nonumber\\
	& \qquad \hbox{for every $v \in V $ and a.e. in $(0,T)$,}\\[2mm]
	\label{cd:2}&\beta\dt\vp-\Delta\vp+\xi + (F'_2(\ph_1)-F'_2(\ph_2))=\mu+\chi\,\sigma
	\hbox{ \, a.e.  in $Q$,}
 	\\
 	\label{cd:3}
 	& \<	\dt\sigma, v>
 	+ \iO \nabla \sigma\cdot \nabla v
 	=\chi\iO \nabla \vp\cdot \nabla v 
 	- \iO \hat R\,  v
 	+ \iO u_2 v	\nonumber\\
	& \qquad \hbox{for every $v \in V $ and a.e. in $(0,T)$,}
	\\
	&
	\label{cd:5}
	\m(0)=\m_0, \quad
	\ph(0)=\ph_0, \quad
	\s(0)=\s_0, \quad \hbox{a.e. in $\Omega$},
\end{align}
where 
\begin{align*}
	\hat R : = R_1-R_2 = (P(\ph_1)-P(\ph_2))(\s_1 + \chi (1-\ph_1) - \m_1) + P(\ph_2)(\s- \chi \ph - \m).
\end{align*}
We take $v= \a\m+\ph$ in \eqref{cd:1}, test \eqref{cd:2} by $(\chi^2 + \tfrac 1\a)\ph$, and let $v= \s$ in \eqref{cd:3}. 
Then, we add the resulting equalities and integrate over $(0,t)$ and by parts to obtain that
\begin{align}
	& \frac 12 \IO2 {(\a\m+\ph)}
	+ \a \I2 {\nabla \m}
	+ \frac {\b} 2 (\chi^2 + \tfrac 1\a)\IO2 {\ph} 
	\non \\ & \qquad
	+ (\chi^2 + \tfrac 1\a) \I2 {\nabla \ph}
	+ (\chi^2 + \tfrac 1\a)\intQt \xi \ph
	+ \frac12 \IO2 \s
	+ \I2 {\nabla \s}
	\non \\ & \quad
	= 
	\frac 12 \pcol{\left(\juerg{\norma{\a\m_0+\ph_0}^2}
	+ {\b}(\chi^2 + \tfrac 1\a) \norma{\ph_0}^2
	+  \norma{\s_0}^2 \right)}
		- \intQt \nabla \m \cdot \nabla \ph
	\non \\ & \qquad 
	+ \intQt \hat R (\a\m+\ph -\s)
	-\intQt (\h_1-\h_2) u_1^1  (\a\m+\ph)
			-\intQt  \h_2 u_1 (\a\m+\ph)
	\non \\ & \qquad
	 - (\chi^2 + \tfrac 1\a)\intQt (F'_2(\ph_1)-F'_2(\ph_2))\ph
	 + (\chi^2 + \tfrac 1\a)\intQt \tfrac 1\a (\alpha \m + \ph) \ph  
	\non \\ & \qquad
		- \frac 1 \a (\chi^2 + \tfrac 1\a)\intQt  |\ph|^2  
		+ (\chi^2 + \tfrac 1\a)\chi \intQt \s \ph 
	+\chi \intQt \nabla \ph \cdot \nabla \s + \intQt u_2 \s \label{pier11}
\end{align}
for a.e. $t\in (0,T)$, where we also used that $\m= \tfrac 1\a(\a\m + \ph) - \tfrac 1\a \ph$. 
Observe that all \juerg{of} the terms on the \lhs\ are nonnegative\juerg{; 
in particular, the fifth} is nonnegative thanks to the monotonicity of $\partial F_1$.
Next, we denote by $I_1,...,\pcol{I_{11}}$ the \pcol{eleven terms} on the \rhs\ of \eqref{pier11}, in the above order,
and \juerg{estimate them individually}. Using the Young inequality, we infer that
\begin{align*}
	I_1 + \pcol{I_2 + I_{10}} 
	\leq\  &{}\pcol{c{} \pcol{\left(\norma{\a\m_0+\ph_0}^2
	+ \norma{\ph_0}^2
	+  \norma{\s_0}^2 \right)}} 
	\\ &{}+\frac12 (\an{\chi^2  + \tfrac 1{\a}}) \I2 {\nabla \ph} 
	+ \frac{\a}{2} \I2 {\nabla \m} 
	+ \frac 12 \I2 {\nabla \s},
\end{align*}
and here \pcol{the last three} contributions on the \rhs\ can be absorbed on the \lhs\ of \eqref{pier11}.
We also immediately observe that $\pcol{I_8} \leq 0$.
Moreover, with the help of \ref{F:weak}, \anfirst{\ref{h:weak}}, and recalling \eqref{pier1}, 
\juerg{we deduce from  Young's} inequality  that
\begin{align*}
	\pcol{I_4 + I_5}  &\leq L \, \|u_1^1\|_{L^\infty (Q)} \intQt |\ph||\a\m+\ph|
	+  \norma{\elvis{\h_2}}_\infty \intQt|u_1||\a\m+\ph| 
		 \\ &\leq 
		c \intQt (|{\a\m+\ph}|^2+|\ph|^2+|u_1|^2)
\end{align*}
as well as
\begin{align*}
	&\pcol{I_6+I_7+I_9 + I_{11}}  
	\\
	&\quad \leq L (\chi^2 + \tfrac 1\a)\intQt|\ph|^2
		+ \frac1{\an{2\a}} (\chi^2 + \tfrac 1\a)\intQt (|\a\m+\ph|^2+|\ph|^2) 
		\\
		&\qquad
		+ \frac\chi 2(\chi^2 + \tfrac 1\a) \intQt (|\s|^2+|\ph|^2) 
         + \frac12 \intQt (|u_2|^2+|\s|^2) 
	\\ &\quad \leq  c \intQt (|\ph|^2+|\a\m+\ph|^2 +|\s|^2 + |u_2|^2).
\end{align*}
It remains to estimate $\pcol{I_3}$. Using the boundedness and \Lip\ continuity of $P$, the \anfirst{\Holder} and Young inequalities, 
and the continuous embedding $V \subset L^4 (\Omega)$, we find \juerg{that}
\begin{align*}
	\pcol{I_3}
		& \leq 
	L \intQt |\ph|(|\s_1| +\chi + \chi |\ph_1|+|\m_1|)(|\a\m+\ph| +|\s|)
	\\ & \quad
		+\norma{P}_\infty \intQt (|\s|+ \tfrac1\a |\a \m +\ph| + (\chi +\tfrac1\a )| \ph| )(|\a\m+\ph| +|\s|)
	\\ & \leq
	c \iot \norma{\ph\juerg{(s)}}_4 \an{(\norma{\s_1}_4 + \norma{\ph_1}_4+\norma{\m_1}_4)\juerg{(s)}\,(\norma{\a\m+\ph}+ \norma{\s})}
	\juerg{(s)\,ds}\\ & \quad
	\an{+c \iot \norma{\ph{(s)}}(\norma{\a\m+\ph}+ \norma{\s})(s)\,ds}
	+ c \intQt (|\s|^2+|\a\m+\ph|^2+|\ph|^2)
	\\ & \leq 
	\d \iot \juerg{\norma{\ph(s)}_V^2\,ds}
	+ \cd \iot \juerg{\bigl(}(\an{\norma{\s_1}_V^2 +\norma{\ph_1}_V^2+\norma{\m_1}_V^2})(\norma{\a\m+\ph}^2+ \norma{\s}^2)
	\juerg{\bigr)(s)\,ds}
	\\ & \quad
+ c \intQt (|\s|^2+|\a\m+\ph|^2+|\ph|^2)\,,
\end{align*}
for a positive $\d$ to be chosen, for instance, less than or equal to $\tfrac 1 4 (\chi^2 + \tfrac 1\a)$. Since $(\m_1,\ph_1,\xi_1,\s_1)$ is a weak solution 
to~\Statesys\ in the sense of Definition~\ref{DEF:WEAK}, it follows that the function
\begin{align*}
	t \mapsto (\an{\norma{\s_1(t)}_V^2 +\norma{\ph_1(t)}_V^2+\norma{\m_1(t)}_V^2})
\end{align*}
\juerg{belongs to} $L^1(0,T)$.
Hence, we can collect all the above inequalities and apply the Gronwall lemma to finally derive the estimate \eqref{cont:dep:weak}.}
\end{proof}

Since the control problem introduced above will demand strong regularities, we also prove the existence of strong solutions \pier{(i.e., regularity results for our weak solutions)} 
to the system \Statesys\ under \pier{further} assumptions. In this direction, in addition to \ref{const:weak}--\pier{\anfirst{\ref{h:weak}}}, we postulate that:\pier{%
\begin{enumerate}[label={\bf (S\arabic{*})}, ref={\bf (S\arabic{*})}]
\item 	\label{F:strong:1}
	\quad $F=F_1+F_2$;  $\,F_1:\erre\to [0,+\infty]\,$ is convex and lower semicontinuous with
\hspace*{2mm}
$\,F_1(0)=0$; $F_2\in \newju{C^5}(\erre)$, and $\,F_2'\,$ is Lipschitz continuous on $\erre$.
\item \label{F:strong:2}
	\quad There exists an interval $\,(r_-,r_+)\,$ with $\,-\infty\le r_-<0<r_+\le +\infty\,$ such that
the \hspace*{3mm} restriction of $F_1$ to $\,(r_-,r_+)\,$ belongs to $\,\newju{C^5}(r_-,r_+)$.
%\item \label{F:strong:3}
%	\quad $F_2\in C^4(\erre)$, and $\,F_2'\,$ is globally Lipschitz continuous on $\erre$.
\item \label{F:strong:4}
	\quad It holds $\,\lim_{r\to r_{\pm}} F'(r)=\pm\infty$.
\item \label{Ph:strong}
	\quad  $P, \h\in C^3(\erre)\cap W^{3,\infty}(\erre)$, and $\h$ is positive on $(r_-,r_+)$.
\end{enumerate}}% 
Observe that \ref{Ph:strong} entails that $\anfirst{P,P',P'',\h,\h',\h''}$ are Lipschitz continuous on $\erre$.
Moreover, let us remark that the above setting allows us to include the singular logarithmic potential \eqref{logpot} and the associated quartic approximation \eqref{regpot}, but it excludes the double obstacle potential \eqref{2obspot}, which cannot be considered in the framework of \ref{F:strong:2}--\ref{F:strong:4}.
\anfirst{Furthermore, the prescribed regularity for the potential $F$ entails that its derivative can be defined in the classical manner so that we no longer need considering a selection $\xi$ in the notion of strong solution below.}
Moreover, it will be useful to set\juerg{, for a fixed $R>0$,}
\begin{equation}\label{defUR}
{\cal U}_R:=\left\{\bu =(u_1,u_2) \in L^\infty(Q)^2:\,\|\bu\|_\infty<R\right\}.
\end{equation}
Under these conditions, we have the following result concerning the well-posedness of the state system \pier{\eqref{ss1}--\eqref{ss5}, where \juerg{the} equations and conditions
have to be fulfilled almost everywhere \juerg{in $Q$}.}
\Bthm[Strong well-posedness]
\label{THM:STRONG}
Suppose that the conditions \juergen{\ref{const:weak}}, \ref{F:strong:1}--\ref{Ph:strong}, \juergen{and \eqref{defUR}}
are fulfilled.
Moreover, let the initial data \anfirst{fulfill}
\begin{align}
	\label{strong:initialdata}
	\m_0,\s_0 \in \Hx1 \cap \Lx\infty, \quad \ph_0 \in \an{ W_0 ,}
\end{align}
as well as
\begin{align}
	r_-<\min_{x\in\overline{\Omega}}\,\vp_0(x) \le
	\max_{x\in\overline{\Omega}}\,\vp_0(x)<r_+.
	\label{strong:sep:initialdata}
\end{align}
Then the state system \eqref{ss1}--\eqref{ss5} has for every $\bu=(u_1,u_2)\in \UR$ a unique solution
$(\mu,\vp,\sigma)$ with the regularity
\begin{align}
\label{regmu}
&\mu\in  H^1(0,T;H) \cap C^0([0,T];V) \cap L^2(0,T;W_0)\cap  L^\infty(Q),\\[1mm]
\label{regphi}
&\ph\in W^{1,\infty}(0,T;H)\cap H^1(0,T;V)\cap L^\infty(0,T;W_0) \cap C^0(\overline Q),
\\[1mm]
\label{regsigma}
&\sigma\in H^1(0,T;H)\cap C^0([0,T];V)\cap L^2(0,T;W_0)\cap L^\infty(Q).
\end{align}
Moreover, there \pier{is} a \juerg{constant} $K_1>0$, which depends on $\Omega,T,R,\alpha,\beta$ and the data of the 
system, but not on the choice of $\bu\in \UR $, such that
\pier{\begin{align}\label{ssbound1}
&\|\mu\|_{H^1(0,T;H) \cap C^0([0,T];V) \cap L^2(0,T;W_0)\cap L^\infty(Q)}\nonumber\\[1mm]
&+\,\|\ph\|_{W^{1,\infty}(0,T;H)\cap H^1(0,T;V) \cap L^\infty(0,T;W_0)\cap C^0(\overline Q)}
\nonumber\\[1mm]&+\,\|\sigma\|_{H^1(0,T;H) \cap C^0([0,T];V) \cap L^2(0,T;W_0)\cap L^\infty(Q)}\,\le\,K_1\,.
\end{align}}%
\pier{Furthermore}, there \pier{exist two values} $r_*,r^*$, depending on $\Omega,T,R,\alpha,\beta$ and the data of the 
system, but not on the choice of $\bu\in \UR$, such that
\begin{equation}\label{ssbound2}
r_-  <r_*\le\vp(x,t)\le r^*<r_+ \quad\mbox{for all $(x,t)\in \overline Q$}.
\end{equation}
\pier{Also, there is some constant $K_2>0$, with the same dependencies as $K_1$,} such that
\begin{align}
\label{ssbound3}
	&\max_{i=0,1,2,3}\,\left\|P^{(i)}(\vp)\right\|_{L^\infty(Q)}\,
	+ \max_{i=0,1,2,3}\,\left\|\h^{(i)}(\vp)\right\|_{L^\infty(Q)}\,
	\non\\
	&+\,\max_{i=0,1,2,3,4\an{,5}}\,\left\|
	F^{(i)}(\vp)\right\|_{L^\infty(Q)} \,\le\,K_2\,.
\end{align}
\pier{Finally, for $i= 1,2 $, let $(\m_i,\ph_i,\s_i)$ be a strong solution to \Statesys\ associated with the initial data $(\m_0^i,\ph_0^i,\s_0^i)$ \juerg{satisfying} \eqref{strong:initialdata}--\eqref{strong:sep:initialdata} \juerg{and the}
 control\anfirst{s} $\bu^i= (u_1^i,u_2^i) \in \pier{{}\UR} $.} Then, \pier{there is a \anfirst{positive} constant $C_{D}$, depending only on data, such that}
\begin{align}
	& \non
	\norma{\m_1-\m_2}_{\H1 H \cap \L\infty V \cap \L2 {W_0}}
	+ \norma{\ph_1-\ph_2}_{\H1 H \cap \L\infty V \cap \L2 {W_0}}
	\\
	& \qquad \non
	+ \norma{\s_1-\s_2}_{\H1 H \cap \L\infty V \cap \L2 {W_0}}
	\\& \quad 
	\leq
	C_{D} \Big(
	\norma{\m_0^1-\m_0^2}_V
	+\norma{\ph_0^1-\ph_0^2}_{V}
	+ \norma{\s_0^1-\s_0^2}_{V}
	+ 	\newju{\left \|\bu^1-\bu^2 \right\|_{L^2(0,T;H)^2}}
		\Big).
	\label{cont:dep:strong}
\end{align}
\Ethm
\Brem
(i) The {\em separation property} \eqref{ssbound2} is particularly important for the case of singular 
potentials such as $F_{\rm log}$. Indeed, it guarantees that the phase variable always stays away
from the critical values $\,r_-,r_+$ that \pier{may} correspond to the pure phases. In this way, the singularity 
is no longer an obstacle for the analysis; \pier{indeed, the values of $\ph$ range in some interval in which $F_1$ is smooth.}
%however, the case of pure phases is then excluded, which is
%not desirable from the viewpoint of medical applications.
%\noindent

\noindent (ii) Notice that \eqref{strong:sep:initialdata} entails that \pier{$F^{(i)}(\vp_0)\in C^0(\overline\Omega)$ for $i=0,1,\ldots, \an{5}$.
This condition can be restrictive for singular potentials; for instance, in the case of
$\,F_{\rm log}$ 
% and $\,F_{\rm obs}\,$ (which however does not \anfirst{fulfill} \ref{F:strong:4})
we have $r_\pm=\pm 1$, so that \eqref{strong:sep:initialdata} excludes the}
pure phases (tumor and healthy tissue) as initial data.
\Erem                      
 
%\vspace{2mm}
Notice also that, owing to Theorem 2.1, the control-to-state operator $$\, \S:\bu=(u_1,u_2)\mapsto 
(\mu,\ph,\sigma)\,$$ is well defined as a
mapping between ${\cal U}=L^\infty(Q)^2$ and the Banach space specified by the regularity results~\eqref{regmu}--\eqref{regsigma}.
Actually, the control-to-state operator $\S$ may be well defined just after Theorem~\ref{THM:WEAK}, but the notion of weak solutions \juerg{proposed there} \pier{(cf.~Definition~\ref{DEF:WEAK})} would not \juerg{suffice for the investigation of}
the optimal control problem ($\cal CP$).

\begin{proof}
\pier{Again, we proceed formally, but still using the Yosida approximation $F'_{1,\eps}$, at least in the first part of the proof. 
Of course, we \juerg{take} for granted all the estimates already done in the existence proof for Theorem~\ref{THM:WEAK},
and start now with additional estimates \an{independent of $\eps$}.}
%since we can easily adapt the lines of argument of \cite[Proof of Thm.~2.2]{CSS1}.
\anfirst{To avoid a heavy notation, we proceed as in Theorem \ref{THM:WEAK} and use the simpler notation $(\m,\ph,\s)$ for the variables of the approximated system instead of $(\m_\eps, \ph_\eps,\s_\eps)$, while we will reintroduce the correct notation exhibiting the dependence of $\eps$ at the end of each estimate.}

\noindent
\pier{{\sc First estimate:}
We rewrite the variational equality \eqref{var:1} as
\begin{align}
	\alpha \<\dt\mu, v > 
	+ \iO \nabla \mu \cdot \nabla v
	= \iO f_\mu v \quad \hbox{ for every $v \in V $ and a.e. in $(0,T)$,}
%	\begin{cases}
%	\a \dt \m - \Delta \m =P(\ph)(\sigma+\chi(1-\ph)-\mu) - \h(\ph)u_1 -\dt\ph =:  f_{\m}
%	\quad &\hbox{in $Q$},
%	\\
%	\dn \m=0 \quad &\hbox{on $\Sigma$,}
%	\\ 
%	\m(0)=\m_0 \quad &\hbox{in $\Omega$}.
%	\end{cases}
	\label{parabolic:mu}
\end{align}
where $f_\mu := -\dt\ph + P(\ph)(\sigma+\chi(1-\ph)-\mu) - \h(\ph)u_1 $ is already \juerg{known to be} 
uniformly bounded in $\L 2 H$ by \eqref{pier7}.  As $\m(0)=\m_0$ is now in $H^1(\Omega)$, 
\juerg{it follows from the regularity theory
for parabolic problems (see, e.g., \cite{L1})}  that
\begin{align}
	\norma{\anfirst{\m_\eps}}_{\H1 H \cap \L\infty V \cap \L2 {W_0}} \leq c,  \label{pier12}
\end{align}
and \eqref{parabolic:mu} can be equivalently rewritten as the equation \eqref{ss1} 
along with the Neumann boundary condition $\dn \m=0$ a.e.~on $\Sigma$. Next, recalling also \eqref{pier9} and arguing similarly for the variational equality \eqref{var:3},
rewritten as 
\begin{align*}
 	&\<	\dt\sigma, v>
 	+ \iO \nabla \sigma\cdot \nabla v
 	= - \iO (\chi \Delta \ph + P(\ph)(\sigma+\chi(1-\ph)-\mu) -   u_2) v	\nonumber\\
	& \qquad \pier{\hbox{for every $v \in V $ and a.e. in $(0,T)$,}}
\end{align*} 
we also deduce that
\begin{align}
	\norma{\anfirst{\s_\eps}}_{\H1 H \cap \L\infty V \cap \L2 {W_0}}\leq c .   \label{pier13}
	\end{align}
\juerg{Hence, \eqref{ss3} holds a.e.~in $Q$, and all of the boundary conditions in \eqref{ss4} hold 
a.e. on $\Sigma$}.
}

\noindent
\pier{{\sc Second estimate:}
\juerg{We} formally differentiate \eqref{pier6} with respect to time, obtaining
\begin{align}
	\label{diff:ss2}
	\b \partial_{t}(\dt \ph) - \Delta (\dt\ph) + F''_{1,\eps} (\ph) \dt\ph = \dt\m + \chi \dt\s - F''_{2} (\ph) \dt\ph =:  g_{\ph},
\end{align}
where $g_{\ph}$ is bounded in $\L2 H$ independently of $\eps$, on account of \eqref{pier12}, \eqref{pier13}, \ref{F:strong:1}, and \eqref{pier7} (indeed, $F''_{2}$
is globally bounded on $\erre$). 
Then, multiplying \eqref{diff:ss2} by $\dt\ph$ and integrating over $\Omega$ and by parts\juerg{, we find that} 
\begin{align}
	\label{second:est:reg}
	\frac \b2\frac d{dt} \norma{\dt\ph}^2
	+ \norma{\nabla \dt\ph}^2
	+ \iO F''_{1,\eps}(\ph)|\dt\ph|^2
	 \leq \norma{ g_{\ph}} \norma{\dt\ph},\juerg{\quad \mbox{a.e. in $(0,T)$,}}
\end{align}
where the third term on the \lhs\ is nonnegative owing to the monotonicity of $F'_{1,\eps}$. Now, we aim to integrate \eqref{second:est:reg} with respect to time. 
Note that taking $t=0$ in \eqref{pier6} produces 
\begin{align*}
	\dt\ph(0)= \tfrac 1\b 
		\big( 
		\Delta\vp_0 - F'_{1,\eps}(\vp_0) +\mu_0+\chi\,\sigma_0 - F'_{2}(\vp_0)		\big),
\end{align*}
\juerg{where the} \rhs\ is bounded in $H$ by virtue of \an{\eqref{pier5},} \eqref{strong:initialdata},\eqref{strong:sep:initialdata},  and \ref{F:strong:1}. 
Hence, we can integrate \eqref{second:est:reg} over $[0,t]$, with $t\in (0,T]$, 
\juerg{to} conclude that
\begin{align}
	\norma{\anfirst{\ph_\eps}}_{W^{1,\infty}(0,T;H) \cap \H1 V}\leq c . \label{pier14} 
\end{align}
}%

\noindent 
\pier{{\sc Third estimate:}
We come back to the elliptic equation \eqref{ph:elliptic} and observe that now we have at hand that $f_\ph $ is bounded in $\L\infty H$.
Then, arguing similarly as in the proof of \eqref{pier9}, \juerg{using} monotonicity and elliptic regularity theory, we easily infer that
\begin{align}
	\norma{F'_{1,\eps}(\anfirst{\ph_\eps})}_{\L\infty H} + \norma{\anfirst{\ph_\eps}}_{\L\infty {W_0}}  \leq c,
	\label{pier15-1} 
\end{align}
so that the \juerg{continuity of the} embedding $W_0 \subset C^0(\ov \Omega)$ entails that 
\begin{align}
	\norma{\anfirst{\ph_\eps}}_{L^\infty (Q)} \leq c .  \label{pier15-2} 
\end{align}}%

\noindent
\pier{{\sc Fourth estimate:}
Next, we \pcol{consider} the parabolic equation \eqref{ss1}, \an{written as }
%(cf.~\eqref{parabolic:mu})
$$
\alpha\dt\mu -\Delta\mu = -\dt\ph + P(\ph)(\sigma+\chi(1-\ph)-\mu) - \h(\ph)u_1 =: f_{\m},  
$$
and observe that now, thanks to \eqref{pier14}, we have that $\dt\ph$, and consequently $f_\m$, are bounded in $ \L\infty H$.
Moreover, we recall \eqref{strong:initialdata} and note that $\m_0 \in \Lx\infty$, in particular. Thus, we can apply the regularity result~\cite[Thm.~7.1, p.~181]{LSU} to 
show that
\begin{align}
		\norma{\anfirst{\m_\eps}}_{L^\infty(Q)} \leq c.
			\label{pier16} 
\end{align}
With similar arguments we can easily obtain the same property for the nutrient variable.
In fact, \an{it suffices to rewrite} \eqref{ss3} as a parabolic equation with forcing term 
$$f_\s = -\chi\Delta\vp+P(\ph)(\sigma+\chi(1-\ph)-\mu)+u_2 , $$ 
\an{and notice that} \eqref{pier15-1} allows us to infer that $\Delta\vp$, and thus $f_\s$, are bounded in $ \L\infty H$. Hence, we can apply the same argument to conclude that
\begin{align}
		\norma{\anfirst{\s_\eps}}_{L^\infty(Q)} \leq c.
		\label{pier17} 
\end{align}}%

\pier{Now, we collect the estimates \eqref{pier12}--\eqref{pier13}, \eqref{pier14}--\eqref{pier17} and point out that they still hold for the real solution $(\m,\ph,\s)$ when passing to the limit as $\eps \searrow 0$, because of the weak or weak star lower semicontinuity of norms. Then, we realize that indeed the global estimate \eqref{ssbound1} in the statement has been proved, with the observation that $L^\infty(Q)$ for $\ph$ is replaced by $C^0 (\ov Q)$, since this continuity property is actually ensured by
$$ \ph \in W^{1,\infty}(0,T;H)\cap L^\infty(0,T;W_0) $$
and the compact embedding $W_0 \subset C^0 (\ov \Omega)$ (see, e.g., \cite[Sect.~8, Cor.~4]{Simon}).}

\noindent
\pier{{\sc Separation property:} At this point, the equation \eqref{ss2} holds for the limit 
functions, with the datum $F'= F_1' + F_2'$ as in \ref{F:strong:1}--\ref{F:strong:4} and 
with the \rhs\ bounded in  $L^\infty(Q)$. Thus, there exists a positive constant $C_*$ for 
which 
\begin{equation} 
\label{pier18}
\norma{\mu+\chi\s}_{L^\infty(Q)}\leq C_* .
\end{equation}
Moreover, the condition \eqref{strong:sep:initialdata} for the initial $\ph_0$ and the growth assumption \ref{F:strong:4} imply the existence of some 
constants $r_{*}$ and $r^{*}$ such that $ r_- < r_{*}\leq r^{*} < r_+ $
and
\begin{gather}
	\label{separation_first}
	 r_{*} {{}\leq{}} \infess_{x\in \Omega} \ph_0(x), \quad r^{*} {{}\geq{}} \supess_{x\in \Omega} \ph_0(x),
	\\
	F'(r)
	 + C_* \leq 0 
	\quad \forall r \in (r_-,r_{*}), \quad 
	F'(r)
	 - C_* \geq 0 
	\quad \forall r \in (r^{*},r_+).
	\label{separation_second}
\end{gather}
Then, let us multiply \eqref{ss2} by $v = (\ph - r^*)^+ - (\ph - r_*)^-$, where the
standard positive $(\,\cdot\,)^+$ and negative $(\,\cdot\,)^- $ parts are used here. 
Then, \an{we} integrate over $Q_t= \Omega \times (0,t)$, for $t\in (0,T],$
and with the help of \eqref{pier18} deduce that
\begin{align*}
	&\frac \b2 \an{\norma{v(t)}^2} + \I2 {\nabla v} \\
	&\quad ={} \int_{Q_t \cap \{\ph < r_{*}\}} (F'(\ph) - \mu - \chi\s) (r_{*}- \ph ) 
	+ \int_{Q_t \cap \{\ph > r^{*}\}} ( \mu + \chi\s - F'(\ph)) (\ph - r^{*}) \\
	&\quad \leq{} \int_{Q_t \cap \{\ph < r_{*}\}} (F'(\ph) + C_*) (r_{*}- \ph ) 
	+ \int_{Q_t \cap \{\ph > r^{*}\}} (C_* - F'(\ph)) (\ph - r^{*}),
\end{align*}
where we also applied \eqref{separation_first} to have that $v(0)=0$.
Note that the \rhs\ above is nonpositive due to 
\eqref{separation_second}, so that $v =0$ almost everywhere,
which in turn implies~that 
%$\ph \leq r^*$  {in} $Q$. In a similar
%manner, we easily conclude that $\ph \geq r_*$ almost everywhere {in} $Q$
%by testing \eqref{sev_proof} by $-(\ph-r_*)^-$.
%Thus, we have just shown that
%\begin{align}
%	\label{seventh_estimate}
%	r_* \leq {\ph} \leq  r^* \quad \aeQ.
%\end{align}
\begin{align}
	\label{seventh_estimate}
	r_* \leq {\ph} \leq  r^* \quad \hbox{a.e. in } Q.
\end{align}
Then, \eqref{ssbound2} is proven, and, at this point, the assumptions~\ref{F:strong:1}--\ref{Ph:strong} enable us to directly deduce~\eqref{ssbound3}.}

On account of the above regularity, the separation property,  and the assumptions \ref{F:strong:1}--\newju{\ref{Ph:strong}},
 we are now in a position to show the refined continuous dependence estimate given by \eqref{cont:dep:strong}. 
In this direction, we employ the notation introduced in the proof of Theorem \ref{THM:WEAK} and consider the system of the differences \eqref{cd:1}--\eqref{cd:5}. Notice that we now have that $F$ is differentiable, so that $\xi + (F'_2(\ph_1)-F'_2(\ph_2) \anfirst{)}= F'(\ph_1)-F'(\ph_2)$.
Moreover, let us remark that due to the separation property \eqref{ssbound2} and to the reinforced assumptions \ref{F:strong:1}--\ref{F:strong:4}, it follows that \elvis{$F'$} is \Lip\ continuous in \juerg{the range of the 
occurring arguments}. 

\noindent
{\sc First estimate:}
We test \eqref{cd:1} by $\mu$, \eqref{cd:2} \juerg{-- to which we add the term $\ph$ on both sides --} by $\dt \ph$, \juerg{as well as} \eqref{cd:3} by $\s$\pcol{. Then, we sum up and integrate over $Q_t$ and by parts. With the help of the cancellation of two terms we} deduce that
\begin{align*}
		& \frac \a 2 \IO2 {\m} 
		+ \I2 {\nabla \m}
		+ \b \I2 {\dt \ph}
		+\frac 12 \norma{\ph\juerg{(t)}}_V^2 
		+ \frac 12 \IO2 {\s}
		+ \I2 {\nabla \s}
		\\ & \quad = 
		 \pcol{\frac 1 2 \left(\a \norma{\m_0}^2
		 + \norma{\ph_0}^2_V
		 + \norma{\s_0}^2 \right)}
		+ \intQt \hat R (\m -\s)
		- \intQt (\h_1-\h_2) u_1^1 \m
		\next\quad 
		- \intQt \h_2 u_1\m
		- \intQt (F'(\ph_1)-F'(\ph_2) \pcol{{}-\ph} ) \dt\ph
		+ \intQt \m \dt \ph
				 \next \quad
		 + \chi \intQt \s \dt\ph
%		 + \intQt \ph \dt \ph 
		 +\chi \intQt \an{\nabla \ph\cdot\nabla \s}  
		 + \intQt u_2 \s= : I_1+ ... + \pcol{I_{9}}.
\end{align*}
Using the Young and \Holder\ inequalities, the \Lip\ continuity and boundedness of $P$ along with the strong regularity~\pcol{\eqref{ssbound1}} of the solutions, we infer that
\pcol{%
\begin{align*}
	\pcol{I_2}
		& \leq 
	L \intQt |\ph|(|\s_1| +\chi + \chi |\ph_1|+|\m_1|)(|\m| +|\s|)
	\\ & \quad
		+\norma{P}_\infty \intQt (|\s|+ \chi | \ph| + |\m| )(|\m| +|\s|)
	\\ & \leq
	 c \,\left(\norma{\s_1}_{L^\infty(Q)}^2 + 1 +\norma{\ph_1}_{L^\infty(Q)}^2+\norma{\m_1}_{L^\infty(Q)}^2\right) \iot \norma{\juerg{\varphi(s)}}^2\juerg{\,ds} 
	 \\ &\quad + 
	c \intQt (\an{|\m|^2+|\ph|^2+|\s|^2}).
\end{align*}
Recalling \eqref{defUR}, we have that
\begin{align*}
	&I_3  \leq  \intQt |\h(\ph_1)-\h(\ph_2)|\, | u_1^1| \, | \m	|
	\\ &\leq 
	L \iot \norma{\juerg{\ph(s)}} \,\norma{u_1^1\juerg{(s)}}_\infty \,\norma{\m\juerg{(s)}}\juerg{\,ds}\leq \frac{LR}2  \intQt (|\ph|^2+|\m|^2).
\\ 
%	&\quad \leq	
%	\anfirst{\frac 1 4 \I2 {\nabla\m}}
%	+ c \iot \juerg{( \norma{\m(s)}^2 + \norma{u_1^1(s)}^2\, \norma{\ph(s)}^2_V )\,ds}.
\end{align*}
Moreover, it is easy to see that 
\begin{align*}
%	\\ 
&  I_\anfirst{4}\leq 
	c \intQt (|\m|^2+|u_1|^2),
	\ \quad	I_\anfirst{5}  
	\leq  
	\frac \beta 4  \I2 {\dt\ph}
	+ c\I2 \ph ,
	\\
&	I_\anfirst{6}+I_7  \leq 
	\frac \beta 4 \I2 {\dt\ph}
	\,+\, c \intQt (|\m|^2+|\s|^2),
	\\ 
&	I_{8}   \leq 
	\frac 1 2 \I2 {\nabla \s}
	\anfirst{\,+\, c \intQt |\nabla \ph|^2 , \ \quad  I_{9}	\leq c \intQt (|\s|^2+|u_2|^2).}
\end{align*}
Hence, we can collect the inequalities and} apply the Gronwall lemma to infer that the differences 
\newju{satisfy}
\begin{align}
	& \norma{\m}_{\L\infty H \cap \L2 V}
	+	\norma{\ph}_{\H1 H \cap \L\infty V}
	+	\norma{\s}_{\L\infty H \cap \L2 V} \non
	\\  & \quad
	\leq \,c\, (\norma{\m_0}+ \norma{\ph_0}_V
	+\norma{\s_0} + \norma{u_1}_{\L2 H } + \norma{u_2}_{\L2 H } ). \label{pierbase}
\end{align}

\anfirst{
\noindent
{\sc Second estimate:}
By exploiting the ellipticity of \pcol{equation~\eqref{cd:2}}, the \Lip\ continuity of $F'$, along with the above estimates, it is \sfw\ to derive that
\begin{align}
	& \norma{\ph}_{ \L2 {W_0}}
	\leq c \,(\norma{\m_0}+ \norma{\ph_0}_V
	+\norma{\s_0} + \norma{u_1}_{\L2 H } + \norma{u_2}_{\L2 H } ). \label{pier19}
	\end{align}
}

\anfirst{
\noindent
{\sc Third estimate:}
We argue in a similar way as in \eqref{parabolic:mu} and 
rewrite \eqref{cd:1} as a parabolic \pcol{variational equality} for $\m=\m_1-\m_2$ with source term given by 
\begin{align*}
	f_\mu&  := -\dt\ph + (P(\ph_1)-P(\ph_2))(\s_1 + \chi (1-\ph_1) - \m_1) 
	+ P(\ph_2)(\s- \chi \ph - \m) \\ & \quad
	- (\h_1-\h_2) u_1^1 + \h_2 u_1.
\end{align*}
\juerg{On account of \pcol{\eqref{ssbound1} and} the above estimates,} we easily deduce that
\begin{align*}
	\norma{f_\m}_{\L2 H} \leq c \,(\norma{\m_0}+ \norma{\ph_0}_V
	+\norma{\s_0} + \norma{u_1}_{\L2 H } + \norma{u_2}_{\L2 H } ).
\end{align*}
Hence, observing that $\m(0)=\m_0=\m_0^1-\m_0^2$ is in $H^1(\Omega)$, \juerg{we can readily infer from the \pcol{parabolic} regularity theory \an{(see, e.g., \cite{L1})}}
that
\begin{align}
	& \norma{\m}_{\H1 H \cap \L\infty V \cap \L2 {W_0}} \non
	\\ & \quad 
	\leq c\, (\pcol{\norma{\m_0}_V}+ \norma{\ph_0}_V
	+\norma{\s_0} + \norma{u_1}_{\L2 H } + \norma{u_2}_{\L2 H } ). \label{pier20}
	\end{align}}% 

\anfirst{
\noindent
{\sc Fourth estimate:}
Arguing in a similar \pcol{manner for the equality}~\eqref{cd:3}, we infer that the 
\pcol{\rhs\ can be rewritten as $\iO f_\sigma v$, with}
\begin{align*}
	f_\s: = - \chi \Delta \ph 
	- (P(\ph_1)-P(\ph_2))(\s_1 + \chi (1-\ph_1) - \m_1) 
	- P(\ph_2)(\s- \chi \ph - \m)
	+ u_2\pcol{.}
\end{align*}
\pcol{Note that $\norma{f_\s}_{\L 2 H}$ is again already bounded as above and 
$\s(0)=\s_0=\s_0^1-\s_0^2$ is in $H^1(\Omega)$}, so that
\begin{align}
	& \norma{\s}_{\H1 H \cap \L\infty V \cap \L2 {W_0}}  \non
	\\ & \quad
	\leq c \,(\norma{\m_0}+ \norma{\ph_0}_V
	+\pcol{\norma{\s_0}_V} + \norma{u_1}_{\L2 H } + \norma{u_2}_{\L2 H } ). 
	\label{pier21}
\end{align}
\pcol{By collecting \eqref{pierbase}--\eqref{pier21}, we finally conclude} the proof of the assertion.}%
\end{proof}

\section{Existence of a Minimizer}\label{MINIMIZER}
\setcounter{equation}{0}

\juerg{Now that} the well-posedness results for system \Statesys\ have been addressed, we can deal with a corresponding optimal control problem,
where the source terms $u_1$ and $u_2$ act as controls. In this direction, we require that the cost functional $\J$ is defined by \eqref{cost} and that \juerg{the following assumptions are fulfilled}:
\begin{enumerate}[label={\bf (C\arabic{*})}, ref={\bf (C\arabic{*})}]
\item \label{ass:control:const}	
	\quad ${\anfirst{b_1}},{\anfirst{b_2}}$ are nonnegative constants, and ${\anfirst{b_0}}$ \,is positive.
\item \label{ass:control:targets}
	\quad \juerg{$\hat \vp_Q\in L^2(Q)$ and $\,\hat\vp_\Omega\in \Ldue$.}
\item \label{ass:control:Uad}
	\quad $\uad=\left\{\bu=(u_1,u_2)\in {\cal U}: \underline u _i\le u_i\le\hat u_i\mbox{\, a.e. in 
$Q$}, \,\,\,i=1,2\right\},$\\[1mm]
\hspace*{5mm} where \,$\underline u_i, \hat u_i\in L^\infty(Q)$\, satisfy \,$\underline u_i\le
\hat u_i$\,  a.e. in $\,Q$, $i=1,2$. 
\end{enumerate}

\vspace{2mm}\noindent
Notice that $\uad$ is a nonempty, closed and convex subset of $\,{\cal U}=L^\infty(Q)^2$. In the 
following, it will sometimes be \an{necessary} to work with a bounded open superset of $\uad$.
We therefore once and for all fix some $R>0$ such that
${\cal U}_R \supset \uad$, \juerg{where $\UR$ is} defined by \eqref{defUR}.
The first result for (${\cal CP}$) concerns the existence of a minimizer, \juerg{where the proof readily} follows from the direct method of calculus of variations,
along with weak and weak star compactness arguments.
\begin{theorem}[Existence of minimizers]
\label{THM:EXISTENCEOPTCONTROL}
Assume that \juergen{\ref{const:weak}}, \ref{F:strong:1}--\ref{Ph:strong}, \juerg{\eqref{strong:initialdata},
\eqref{strong:sep:initialdata},}
and \ref{ass:control:const}--\ref{ass:control:Uad} hold true.
Then the minimization problem (${\cal CP}$) admits a minimizer. 
\end{theorem}

\begin{proof}
At first, let us notice that $\J$ is nonnegative, so that we can pick a minimizing sequence $\{\bu_n\}_{n\in\enne}
 \subset \uad$
with the corresponding sequence of states $\{(\m_n,\ph_n,\s_n)\}_{n\in\enne}$ to have
\begin{align*}
	\lim_{n\to\infty}\J ((\m_n,\ph_n,\s_n), \bu_n) = \inf_{\bf v \in \uad} \J(\S(\bf v), \bf v) \geq 0\pcol{,}
\end{align*}
\pcol{where $\S(\bf v)$ denotes the state corresponding to the control $\bf v$.}
Furthermore, by combining the estimates \eqref{ssbound1}--\eqref{ssbound3}, which are uniform with respect to $n$,
with the structure of $\uad$, it is a standard matter (upon extracting a subsequence that we do not relabel) to infer the existence of limits 
$\ov \bu \in \uad$ and $\ov \m,\ov \ph, \ov \s$ such that\an{, as $n \to \infty$,}
\begin{align*}
	\bu_n & \to \ov \bu \quad \hbox{weakly star in $L^\infty(Q)^2$},
	\\
	\m_n & \to \ov \m \quad \hbox{weakly star in $H^1(0,T;H) \cap \pcol{\L \infty V} 
	 \cap L^2(0,T;W_0)\cap \anfirst{L^\infty(Q)}$},	
	\\
	\ph_n & \to \ov \ph \quad \hbox{weakly star in $W^{1,\infty}(0,T;H)\cap H^1(0,T;V) \cap L^\infty(0,T;W_0)\cap 
	\juerg{L^\infty(Q)}$},	
	\\
	\s_n & \to \ov \s \quad \hbox{weakly star in $H^1(0,T;H) \cap \pcol{\L \infty V} \cap L^2(0,T;W_0)\cap \anfirst{L^\infty(Q)}$}\anfirst{,}
\end{align*}
\juerg{and, by \pcol{the compactness of the embedding $W_0 \subset C^0(\ov \Omega)$,
also that~(see, e.g., \cite[Sect.~8, Cor.~4]{Simon}),}}
\begin{align*}
	\ph_n \to \bph \quad \hbox{strongly in $C^0(\ov Q)$.}
\end{align*} 

It is then a standard matter to pass to the limit in the \pcol{formulation \eqref{ss1}--\eqref{ss5} written for $(\m_n,\ph_n,\s_n)$ and $ \bu_n$} to conclude that
$(\ov \m,\ov \ph,\ov \s)= \S(\ov \bu)$, and then to 
exploit the weak lower semicontinuity of $\J$ to derive that $\ov \bu$ is a minimizer for the optimization problem~($\cal CP$).
\end{proof}

\section{Differentiability Properties of the Solution Operator}
\label{SEC:DIFF}
\setcounter{equation}{0}
 
We now discuss the Fr\'echet differentiability of  $\S$, considered 
as a mapping between suitable Banach spaces. To show such a result,
it is favorable to employ the implicit function theorem, because, if applicable, it yields
that the control-to-state operator automatically inherits the differentiability order from that of
the involved nonlinearities. The technique employed here is an adaptation of that used recently in \cite{ST} in a similar
context. For the reader's convenience, we give the details of the argument. 
For this, some functional analytic preparations are in order. We first define the linear spaces
\juerg{
\begin{align}
\X\,&:=\,X\times \widetilde X\times X, \mbox{\,\,\, where }\nonumber\\
X\,&:=\,H^1(0,T;H)\cap {\C0 V }\cap L^2(0,T;W_0)\cap L^\infty(Q), \nonumber\\
\widetilde X\,&:=W^{1,\infty}(0,T;H)\cap H^1(0,T;V)\cap L^\infty(0,T;W_0)\cap C^0(\overline Q),
\label{calX}
\end{align}
}
which are Banach spaces when endowed with their natural norms. Next, we introduce the linear space
\begin{align}
&\Y\,:=\,\bigl\{(\mu,\ph,\sigma)\in \calX: \,\alpha\dt\mu+\dt\ph-\Delta\mu\in \Liq, 
\,\,\,\beta\dt\ph-\Delta\ph-\mu\in \Liq,\nonumber\\
& \hspace*{14mm}\dt\sigma-\Delta\sigma+\chi\Delta\ph\in\Liq\bigr\},
\label{calY}
\end{align}
which becomes a Banach space when endowed with the norm
\begin{align}
\|(\mu,\ph,\sigma)\|_{\Y}\,:=\,&\|(\mu,\ph,\sigma)\|_\calX
\,+\,\|\alpha\dt\mu+\dt\ph-\Delta\mu\|_{\Liq}
\,+\,\|\beta\dt\ph-\Delta\ph-\mu\|_{\Liq}\nonumber\\
&+\,\|\dt\sigma-\Delta\sigma+\chi\Delta\ph\|_{\Liq}\,.
\label{normY}
\end{align}

Finally, we fix constants $\tilde r_-,\tilde r_+$ such that
\begin{equation}
\label{openphi}
r_-<\tilde r_-<r_*\le r^*<\tilde r_+<r_+,
\end{equation}  
with the constants introduced in \pcol{\bf (S2)} and \eqref{ssbound2}. We then consider the set
\begin{equation}
\label{Phi}
\Phi\,:=\,\bigl\{\anfirst{(\mu,\ph,\sigma)}\in \Y: \tilde r_-<\ph(x,t) <\tilde r_+ \,
\mbox{ for all }\,(x,t)\in\overline Q\bigr\},
\end{equation}
which is obviously an open subset of the space $\Y$.  

We first prove an auxiliary result for the linear initial-boundary value problem
\begin{align}\label{aux1}
&\alpha\dt\mu+\dt\ph-\Delta\mu\,=\,\lambda_1\left[P(\bvp)(\sigma-\chi\ph-\mu)+P'(\bvp)(\bsigma+\chi(1-\bvp)-\bmu)\ph
- \h'(\bvp)\,\uebar\,\ph\right]\nonumber\\
&\hspace*{37mm}-\lambda_2 \h(\bvp)\anfirst{h_1}+\lambda_3 f_1\quad\mbox{in }\,Q,\\[1mm]
\label{aux2}
&\beta\dt\ph-\Delta\ph-\mu\,=\,\lambda_1\left[\chi\,\sigma-F''(\bvp)\ph\right]+\lambda_3 f_2 \quad\mbox{ in \,$Q$},
\\[1mm]
\label{aux3}
&\dt\sigma-\Delta\sigma+\chi\Delta\ph\,=\,\lambda_1\left[-P(\bvp)(\sigma-\chi\ph-\mu)-
P'(\bvp)(\bsigma+\chi(1-\bvp)-\bmu)\ph \right]\nonumber\\
&\hspace*{37mm}+\lambda_2 \anfirst{h}_2+\lambda_3 f_3\quad\mbox{ in \,$Q$},\\[1mm]
\label{aux4}
&\dn\mu\,=\,\dn\ph\,=\,\dn\sigma\,=\,0 \quad\mbox{ on \,$\Sigma$},\\[1mm]
\label{aux5}
&\mu(0)\,=\,\lambda_4 \mu_0,\quad \ph(0)\,=\,\lambda_4\ph_0, \quad \sigma(0)\,=\,\lambda_4
\sigma_0,\,\,\,\mbox{ in }\,\Omega,
\end{align}
which for $\lambda_1 = \lambda_2 = 1$ and $\lambda_3 = \lambda_4 = 0$ coincides  with the linearization 
of the state equation at $((\uebar,\uzbar),(\bmu,\bvp,\bsigma))$. We will need this slightly more general 
version later for the application of the implicit function theorem. \juerg{To shorten the exposition, we introduce the Banach space of the initial data
satisfying \eqref{strong:initialdata}, which is given by}
\begin{align}
\label{defN}
\juerg{{\cal N}\,:=}&\juerg{\{(\mu_0,\varphi_0,\sigma_0):\mu_0,\sigma_0\in \Huno\cap L^\infty(\Omega), \,\,\varphi_0\in W_0\},}
\end{align}
\juerg{endowed with its natural norm.}
\Blem \label{LEM:FRE}
Suppose that $\lambda_1,\lambda_2,\lambda_3,\lambda_4 \in\{0,1\}$ are given
and that the assumptions \ref{const:weak}, \ref{F:strong:1}--\ref{Ph:strong}, \juergen{and 
\eqref{defUR}--\eqref{strong:sep:initialdata}} are fulfilled. Moreover, let $((\uebar,\uzbar),
 (\bmu,\bvp,\bsigma))\linebreak
\in{\cal U}_R\times \Phi$ be arbitrary. Then the linear initial-boundary value problem \eqref{aux1}--\eqref{aux5}
has for every $\anfirst{(h_1,h_2)}\in \Liq^2$ and $(f_1,f_2,f_3)\in {L^\infty(Q)\times(H^1(0,T;H)\cap L^\infty(Q))
\times L^\infty(Q)}$ a unique solution
\,$(\mu,\ph,\sigma)\in {\mathcal Y}$. 
Moreover, the linear mapping $$\,(\anfirst{(h_1,h_2)},(f_1,f_2,f_3),(\mu_0,\ph_0,\sigma_0))
\mapsto (\mu,\ph,\sigma)\,$$ is continuous from $\,\Liq^2\times \bigl(L^\infty(Q)\times (H^1(0,T;H)\cap \Liq)
\times\Liq\bigr)\times \juerg{{\mathcal N}}$
into $\,\Y$.
\Elem                 

\Bdim
We use a standard Faedo--Galerkin approximation. To this end, let $\{\lambda_k\}_{k\in\enne}\,$ and 
$\,\{e_k\}_{k\in\enne}\,$ denote the eigenvalues and associated eigenfunctions of the eigenvalue problem
$$-\Delta y+y=\lambda y\quad\mbox{in }\,\Omega,\quad \dn y=0\quad\mbox{on }\,\Gamma,$$
where the latter are normalized by $\,\an{\|e_k\|}=1$ \pcol{for all $k\in \enne$}. Then $\,\{e_k\}_{k\in\enne}\,$ forms a complete orthonormal   
system in $\,H\,$ which is also dense in $\,V$. We put $\,V_n:={\rm span}\,\{e_1,\ldots,e_n\}$, $n\in\enne$, noting
that $\,\bigcup_{n\in\enne}V_n\,$ is dense in $\,V$. We look for functions
of the form
$$\mu_n(x,t)=\sum_{k=1}^n \rho_k^{(n)}(t)e_k(x), \quad \ph_n(x,t)=\sum_{k=1}^n v_k^{(n)}(t)e_k(x), \quad 
\sigma_n(x,t)=\sum_{k=1}^n w_k^{(n)}(t)e_k(x),$$
that satisfy the system
\begin{align}
\label{gs1}
&(\alpha\dt\mu_n(t),v)+(\dt\ph_n(t),v)+(\nabla\mu_n(t),\nabla v)\,=\,(z_{n1}(t),v),
\\
\label{gs2}
&(\beta\dt\vp_n(t),v)+(\nabla\vp_n(t),\nabla v)-(\mu_n(t),v)\,=\,(z_{n2}(t),v),
\\
\label{gs3}
&(\dt\sigma_n(t),v)+(\nabla\sigma_n(t),\nabla v)-\chi(\nabla\ph_n(t),\nabla v)\,=\,(z_{n3}(t),v),
\\
\label{gs4}
&\mu_n(0)=\lambda_4 \Pn \mu_0, \quad \vp_n(0)=\lambda_4 \Pn \ph_0,\quad \sigma_n(0)=\lambda_4 \Pn\sigma_0,
\end{align} 
for all $v\in V_n,$ and almost every $\,t\in (0,T),$
where $\,\Pn\,$ denotes the $H^1(\Omega)$--orthogonal projection onto $\,V_n$, and where
\begin{align}
\label{zn1}
&z_{n1}\,=\, \lambda_1\left[P(\bvp)(\sigma_n-\chi\ph_n-\mu_n)+P'(\bvp)(\bsigma+\chi(1-\bvp)-\bmu)\ph_n -\h'(\bvp)\,\uebar\,\ph_n\right]\nonumber\\
&\hspace*{11.5mm}-\lambda_2 \h(\bvp)\,\anfirst{h_1} 
+\lambda_3 f_1,\\
\label{zn2}
&z_{n2}\,=\,\lambda_1\left[\chi\,\sigma_n-F''(\bvp)\ph_n\right]+\lambda_3 f_2,\\
\label{zn3}
&z_{n3}\,=\,\lambda_1\left[-P(\bvp)(\sigma_n-\chi\ph_n-\mu_n)-P'(\bvp)(\bsigma+\chi(1-\bvp)-\bmu)\ph_n\right]\nonumber\\
&\hspace*{11.5mm}+\lambda_2 \anfirst{h_2}+\lambda_3 f_3\,.
\end{align}
Insertion of $v=e_k$, for $k\in\enne$, in
\eqref{gs1}--\eqref{gs3},
and substitution for the second summand in \eqref{gs1} by means of \eqref{gs2}, then lead to an initial value
problem for an explicit  linear system of ordinary differential equations for the unknowns $\,\rho_1^{(n)},\ldots,\rho_n^{(n)},v_1^{(n)},\ldots,v_n^{(n)},\linebreak w_1^{(n)},\ldots,w_n^{(n)}$, in which all of the coefficient
functions belong to $L^\infty(0,T)$. Hence, by virtue of Carath\'eodory's theorem, there exists a unique solution
in $W^{1,\infty}(0,T;\erre^{3n})$ that specifies the unique solution $(\mu_n,\ph_n,\sigma_n)\in W^{1,\infty}(0,T;
\pcol{W_0})^3$ to the system \eqref{gs1}--\eqref{gs4}, for $n\in\enne$.

We now derive some a priori estimates for the Faedo--Galerkin approximations. 
In this procedure, $C_i>0$,
$i\in\enne$, will denote constants that are independent of $n\in\enne$ and the data $((f_1,f_2,f_3),(\mu_0,\ph_0,
\sigma_0))$, while the constant $M>0$ is given by
\begin{align}
\label{defM}
M:=&\lambda_2\,\|\anfirst{(h_1,h_2)}\|_{L^\infty(Q)^2}\,+\,\lambda_3\,{\|(f_1,f_2,f_3)\|_{L^\infty(Q)\times(H^1(0,T;H)\cap \Liq)\times\Liq}}\nonumber\\[0.5mm]
&+\,\lambda_4\,\|(\mu_0,\ph_0,\sigma_0)\|_{\juerg{\mathcal N}}.
\end{align}
Moreover, $\,(\bmu,\bvp,\bsigma)\in \Phi$, and thus it follows that $\juergen{\bmu,}\bsigma\in L^\infty(Q)$ and $\juergen{\bvp},\h(\bvp), \h'(\bvp),P(\bvp),\linebreak
P'(\bvp), F''(\bvp)\in
C^0(\overline Q)$. Hence, there is some constant $C_1>0$ such that, for a.e. $(x,t) \in  Q$ and for all $n\in\enne$,
\begin{align}
\big(|z_{n1}|+|z_{n2}|+|z_{n3}|\big)(x,t)\,&\le\,C_1\,\big(\lambda_1(|\mu_n|+|\vp_n|+|\sigma_n|)(x,t) +
 \lambda_2(|\anfirst{h_1}|+|\anfirst{h_2}|)(x,t)\nonumber\\
& \quad\quad\,+ \lambda_3(|f_1|+|f_2|+|f_3|)(x,t)\big)\nonumber\\[0.5mm]
 &\,\,\le\,C_1\,\bigl(\lambda_1(|\mu_n|+|\vp_n|+|\sigma_n|)(x,t)\,+\,M\bigr).
\label{basic}
\end{align}

\vspace{1mm}\noindent
{\sc First estimate:}
We insert $v=\mu_n(t)$ in \eqref{gs1}, $v=\dt\ph_n(t)$ in \eqref{gs2}, $v=\sigma_n(t)$ in \eqref{gs3}
and add the resulting equations, whence a cancellation of two terms occurs. Then, we add to both sides of the 
resulting equation the same term $\,{\frac 12\frac d{dt}\an{\|\ph_n(t)\|^2}=}(\ph_n(t),\dt\ph_n(t))$. Integration over $[0,\tau]$, where $\tau\in (0,T]$, then
yields the identity
\begin{align}
\label{John}
&\frac 1 2\,\bigl(\alpha\|\mu_n(\tau)\|^2+\|\ph_n(\tau)\|_V^2 {+\|\sigma_n(\tau)\|^2\bigr)}\,+{\int_{Q_\tau} \bigl(|\nabla\mu_n|^2+|\nabla\sigma_n|^2\bigr)}
\,+\,\beta\int_{Q_\tau} |\dt\ph_n|^2\nonumber\\
&\quad={\,\frac{\lambda_4^2} 2\left(\alpha\|\Pn\mu_0\|^2+\|\Pn\ph_0\|_V^2+\|\Pn\sigma_0\|^2\right)}
\,+\int_0^\tau (\mu_n(t),z_{n1}(t))\,dt \,+{\int_0^\tau(\an{z_{n3}(t),\sigma_n(t)})\,dt}\nonumber\\
&\qquad{+\int_0^\tau(\an{z_{n2}(t)+\ph_n(t),\dt\ph_n(t)})\,dt\,+\,\chi\int_0^\tau(\nabla\ph_n(t),\nabla\sigma_n(t))\,dt
\,=:\,\sum_{i=1}^5 \an{I}_i,}
\end{align}                                  
with obvious notation. We estimate the terms on the right-hand side individually. First observe that\juerg{,
for all $n\in\enne$,}
\begin{equation}
\label{inibound}
{\an{I_1}\,\le\,C_2\,\lambda_4^2
\,\|(\Pn\mu_0,\Pn\ph_0,\Pn\sigma_0)\|_{V\times V\times V}^2\,\le\,C_2\,\lambda_4^2\,
\|(\mu_0,\ph_0,\sigma_0)\|_{V\times
V\times V}^2\,\le\,C_2\,M^2.}
\end{equation}
Moreover, by virtue of \eqref{basic} and Young's inequality, \pcol{we deduce that}
\begin{align}
\label{Paul}
\an{I_2 + I_3}\,\le\,C_3\,M^2\,+\,C_4\int_{Q_\tau}\left(|\mu_n|^2+|\ph_n|^2+|\sigma_n|^2\right).
\end{align}
Likewise, \pcol{it results that}
\begin{align}
\label{George}
\an{I_4}\,\le\,&\,\frac{\beta}2\int_{Q_\tau}|\dt\ph_n|^2\,+\,\frac{C_5}{\beta}\,M^2\,+\,\frac{C_6}{\beta}
\int_{Q_\tau}\left(|\mu_n|^2+|\ph_n|^2+|\sigma_n|^2\right).	
\end{align}
Finally, we have that
\begin{equation}
\label{Ringo}
\an{I_5}\,\le\,\an{\frac {\chi^2}2\int_{Q_\tau}|\nabla\ph_n|^2+\frac 12\int_{Q_\tau}|\nabla\sigma_n|^2\,}.
\end{equation}
Combining the estimates \eqref{John}--\eqref{Ringo}, we have shown that
\begin{align*}
&{\frac 12\left(\alpha\|\mu_n(\tau)\|^2+\|\ph_n(\tau)\|_V^2+\|\sigma_n(\tau)\|^2\right)\,+\,\int_{Q_\tau}
\bigl(|\nabla\mu_n|^2\,+\,\frac 12\,|\nabla\sigma_n|^2\bigr)\,+\,\frac{\beta}2\int_{Q_\tau}|\dt\ph_n|^2
}\nonumber\\
&\quad{\le\,C_7\,M^2\,+\,C_8\int_0^\tau\bigl(\|\mu_n(t)\|^2
+\|\ph_n(t)\|_V^2+\|\sigma_n(t)\|^2\bigr)\,dt\,.}
\end{align*}
{Therefore, invoking
Gronwall's lemma, we conclude that}
\begin{align}
& \non
\|\mu_n\|_{L^\infty(0,T;H)\cap L^2(0,T;V)}\,+\,\|\ph_n\|_{H^1(0,T;H)\cap L^\infty(0,T;V)}
\\ & \quad
+\,\|\sigma_n\|_{L^\infty(0,T;H)\cap L^2(0,T;V)}
\,\le\,C_9 M \quad\an{\forall\,n\in\enne}.                                             
\label{esti1}
\end{align}  
\vspace{1mm}\noindent
{\sc Second estimate:} 
Next, we insert $v=\dt\mu_n(t)$ in \eqref{gs1} and integrate 
over $[0,\tau]$, where $\tau\in (0,T]$, to obtain the identity
\begin{align*}
&{\frac 12\,\|\nabla\mu_n(\tau)\|^2\,+\,\alpha\int_0^\tau\|\dt\mu_n(t)\|^2\,dt}\\
&\quad{=\,\frac {\lambda_4^2}2\,\|\nabla \Pn\mu_0\|^2 +\int_0^\tau(\an{z_{n1}
(t),\dt\mu_n(t)})\,dt
\,-\int_0^\tau(\an{\dt\ph_n(t),\dt\mu_n(t)})\,dt.}
\end{align*}
{Applying Young's inequality appropriately \pcol{and making} use of 
\eqref{basic} and \eqref{esti1}, we conclude the estimate}
\begin{align}\label{esti2}
&\|\mu_n\|_{H^1(0,T;H)\cap L^\infty(0,T;V)}\,\le\,C_{10}\,M\quad\pcol{\forall\,n\in\enne}.\end{align}

\vspace{1mm}\noindent
{\sc Third estimate:} \quad At this point, we insert $v=-\Delta\mu_n(t)$ in \eqref{gs1} and $v=-\Delta\ph_n(t)$
in \eqref{gs2}, add, and integrate over $[0,\tau]$, where $\tau\in (0,T]$. We then obtain that
\begin{align*}
&\frac \alpha 2\,\|\nabla\mu_n(\tau)\|^2\,+\,\frac{\beta} {2}\,\|\nabla\ph_n({\tau})\|^2\,
+\int_0^\tau\|\Delta\mu_n(t)\|^2\,dt
\,+\int_0^\tau\|\Delta\ph_n(t)\|^2\,dt \\
&\quad=\,\frac{\alpha\lambda_4^2} 2\|\nabla \Pn\mu_0\|^2\,+\,\frac {\beta\lambda_4^2}2 \|\nabla \Pn\ph_0\|^2\,-\int_0^\tau(\an{z_{n1}(t) \pcol{{}- \partial_t \ph_n(t)},\Delta\mu_n(t)})\,dt\\
&\qquad \,\,-\int_0^\tau(\an{\mu_n(t)+z_{n2}(t),\Delta\ph_n(t)})\,dt\,,
\end{align*}
whence, using \pcol{\eqref{basic}, \eqref{esti1}, \eqref{esti2}} and Young's inequality, 
\begin{equation}
\label{Hugo}
\int_0^T\bigl(\|\Delta\mu_n(t)\|^2\,+\,\pcol{\|\Delta\ph_n(t)\|^2}\bigr)\,dt\,\le\,C_{11}M^2\quad\forall\,n\in\enne.
\end{equation}
Classical elliptic estimates, using \eqref{aux4} and \eqref{esti1}--\eqref{Hugo}, then yield that 
\begin{equation}\label{esti3}
\|\mu_n\|_{L^2(0,T;\Hdue)}\,+\,\|\ph_n\|_{L^2(0,T;\Hdue)}\,\le\,C_{12} M \quad\forall\,n\in\enne.
\end{equation}
With the estimate \eqref{esti3} at hand, we may (by first taking $v=\dt\sigma_n(t)$ in \eqref{gs3} and then $v=-\Delta\sigma_n(t)$) infer by similar reasoning that also
\begin{equation}\label{esti4}
\|\sigma_n\|_{H^1(0,T;H)\cap L^\infty(0,T;V)\cap L^2(0,T;\Hdue)}\,\le\, C_{13} M \quad\forall\,n\in\enne.
\end{equation}

At this point, we can conclude from standard weak and weak star compactness arguments the existence of a triple
$(\mu,\ph,\sigma)$ such that, possibly only on a subsequence which is again indexed by $\,n$, 
\begin{align*}
&\mu_n\to\mu,\quad \ph_n\to\ph,\quad \sigma_n\to\sigma,\\[0.5mm]
&\mbox{all weakly star in }\,H^1(0,T;H)\cap L^\infty(0,T;V)\cap L^2(0,T;\pcol{W_0}).
\end{align*}

Standard arguments, which need no repetition here, then show that $(\mu,\ph,\sigma)$ is a strong solution to the system \eqref{aux1}--\eqref{aux5}. \juerg{In particular, it turns out that $\dn\mu=\dn\ph=\dn\sigma=0$ almost everywhere on
$\Sigma$.} 
Moreover, recalling \eqref{esti1}--\eqref{esti4}, and invoking the weak sequential lower
semicontinuity of norms, we conclude that there is some $C_{14}>0$ such that
\begin{align}
\label{esti5}
\|(\mu,\ph,\sigma)\|_{\anfirst{\Z}}\,\le\,C_{14}M ,
\end{align}
\pcol{where $\Z$ is defined in \eqref{defZ}.}

\andrea{%
\Brem
\label{REM:L2}
A careful inspection of the estimates \eqref{John}--\eqref{esti4} shows that the \juergen{term} $M$ appearing on the \rhs\ of the above estimates can be 
replaced by 
\begin{align*}
\juergen{\widetilde	M:=}&\lambda_2\,\|(h_1,h_2)\|_{\L2 H^2}\,+\,\lambda_3\,\juergen{\|(f_1,f_2,f_3)\|_{L^2(0,T;H)^3}}\,+\,\lambda_4\,\|(\mu_0,\ph_0,\sigma_0)\|_{\pcol{V^3}},
\end{align*}
\juergen{since\pcol{, in particular,} only the $L^2(Q)$ norms of the increments $h_i$, $i=1,2$, and of the source terms $f_i$, $i=1,2,3$, enter} 
the computations.
This, along with \eqref{esti5}, entails that the linearized variables $(\mu,\ph,\sigma)$ (which correspond to the choices $\lam_1=\lambda_2=1,\lambda_3=\lambda_4=0$)
\juergen{satisfy}
\begin{align}
	\label{est:rem}
	\|(\mu,\ph,\sigma)\|_{\anfirst{\Z}}\,\leq c\, \norma{\bh}_{\L2 H^2},
\end{align}
\juergen{with some} positive constant $c$ (cf. also Remark \ref{REM:FRE})\pcol{, being 
$\bh = (h_1, h_2)$.}
\Erem}

\juerg{We now derive further estimates for $(\mu,\ph,\sigma)$.}
%At first, we observe that $\ph$ solves the
%parabolic problem
%$$\juerg{\beta\dt\ph-\Delta\ph=\mu+\lambda_1(\chi\sigma-F''(\overline \ph)\ph)+\lambda_3 f_2=:f_\ph,}$$
%\juerg{with zero Neumann boundary condition and initial datum $\ph_0\in W_0\subset L^\infty(\Omega)$, 
%where we know from \eqref{ssbound3} and \eqref{esti5} that
%$\|f_\ph\|_{L^\infty(0,T;H)}\le C_{15}M$. It thus follows from \cite[Theorem~7.1]{LSU} that
%$\ph\in\Liq$ and}
%\begin{equation}
%\label{uwe1}
%\juerg{\|\ph\|_{L^\infty(Q)}\,\le\,C_{16}M.}
%\end{equation}
\juerg{\pcol{In the next one we argue formally, noting that it} can be carried out rigorously on the level of
the Faedo--Galerkin approximations. Indeed, we differentiate \eqref{aux2} formally with respect to time to obtain the
identity}
\begin{equation}
\label{uwe2}
\juerg{\beta\dt(\dt\ph)-\Delta\dt\vp\,=\,\dt\mu+\lambda_1\bigl(\chi\dt\sigma-F'''(\bvp)\dt\bvp\,\ph
-F''(\bvp)\dt\ph\bigr)+\lambda_3\dt f_2=:g_\ph.}
\end{equation},
\juerg{Testing \eqref{uwe2} formally by $\dt\ph$ and integrating formally by parts, we find that}
\begin{equation}
\label{uwe3}
\juerg{\frac {\beta}2\|\dt\ph(t)\|^2+\int_{Q_t}|\nabla\dt\ph|^2=\frac{\beta}2\|\dt\ph(0)\|^2+
\int_{Q_t}g_\ph\,\dt\ph\,.}
\end{equation}
\juerg{Now observe that, owing to \eqref{esti5}, \eqref{ssbound1}, \eqref{ssbound3}, \pcol{and \eqref{defM}} we have that}
\begin{equation}\label{uwe4}
\pcol{\int_{Q_t}g_\ph\,\dt\ph \leq  \|g_\ph\|_{L^2(0,T;H)}\,\|\dt\ph\|_{L^2(0,T;H)}\le\,C_{\pcol{15}}M^2.}
\end{equation} 
\pcol{As for \eqref{uwe4}, we point out that the term $- \lambda_1 F'''(\bvp)\dt\bvp\,\ph$, which is part of $g_\ph$, is bounded in $L^2(0,T;H)$ since $V\subset L^4 (\Omega) $ with continuous embedding, and consequently it follows that $\norma{\dt\bvp}_{\L 2{L^4 (\Omega)}} $ is under control and 
$\norma{\ph}_{\L \infty {L^4 (\Omega)}} \leq C_{16}\, M$, whence
\begin{align*}
&\norma{- \lambda_1 F'''(\bvp)\dt\bvp\,\ph}_{L^2(0,T;H)} \\
&\quad \le \lambda_1 
\norma{F'''(\bvp)}_{L^\infty (Q)} \norma{\dt\bvp}_{\L 2{L^4 (\Omega)}} \norma{\ph}_{\L \infty {L^4 (\Omega)}} \le  C_{17}\, M. 
\end{align*}
Next,} \juerg{writing
\eqref{aux2} for $t=0$ \pcol{and recalling \eqref{aux5}}, we have that}
\begin{equation*} 
\juerg{\dt\ph(0)=\beta^{-1}\bigl(\pcol{\lambda_4(\Delta\ph_0+\mu_0+\lambda_1\chi\sigma_0-\lambda_1 F''(\ph_0)\ph_0)}+\lambda_3\,f_2(0)\bigr),}
\end{equation*}
\juerg{and it follows from \eqref{strong:initialdata}, 
\pcol{\eqref{strong:sep:initialdata}, 
\eqref{defM}} that}
\begin{equation}\label{uwe5}
\juerg{\|\dt\ph(0)\|\le C_{18}M.} 
\end{equation}
\juerg{Combining \eqref{uwe3}--\eqref{uwe5}, and invoking Young's inequality and Gronwall's lemma, we thus can conclude that}
\begin{equation}\label{uwe6}
\juerg{\|\ph\|_{W^{1,\infty}(0,T;H)\cap H^1(0,T;V)}\le C_{19}M.}
\end{equation}
%\juerg{From now on, our argumentation will again be rigorous.}
\juerg{\pcol{Now, in view of \eqref{uwe6} and \eqref{esti5}, a comparison of terms} in \eqref{aux2} and standard elliptic estimates yield that}
\begin{equation}\label{uwe7}
\juerg{\|\ph\|_{L^\infty(0,T;\Hdue)}\le C_{20}M,}
\end{equation}
\juerg{and the compactness of the embedding $(W^{1,\infty}(0,T;H)\cap L^\infty(0,T;\Hdue))\subset C^0(\overline Q)$
(see, e.g., \cite[Sect.~8, Cor.~4]{Simon})
then shows that also}
\begin{equation}\label{uwe8}
\juerg{\|\ph\|_{C^0(\overline Q)}\le C_{21}M.}
\end{equation}
\juerg{At this point, we observe that, by bringing the term $\dt\ph$ to the \rhs, equation \eqref{aux1} can be rewritten as
a linear parabolic equation for $\,\mu\,$ whose \rhs\ is already known to be bounded 
in $L^\infty(0,T;H)$ by an expression of the form
$\,C_{22}M$. Since $\mu$ satisfies zero Neumann boundary conditions and $\mu_0\in \Huno\cap L^\infty(\Omega)$, we 
can apply the classical result of \cite[Theorem~7.1]{LSU} to conclude that $\mu\in\Liq$ and}
\begin{align}\label{uwe9}
\juerg{\|\mu\|_{\Liq}\le\,C_{23}M.}
\end{align}
\juerg{Similar reasoning on equation \eqref{aux3}, invoking the $L^\infty(0,T;H)$-bound  for $\Delta\ph$ implied by
\eqref{uwe7}, shows that
also $\sigma\in\Liq$ and}
\begin{equation}\label{uwe10}
\juerg{\|\sigma\|_{\Liq}\le C_{24}M.}
\end{equation}
\pcol{About the linear dependence of the right-hand sides of \eqref{uwe9} and \eqref{uwe10} on the constant $M$ that is specified in \eqref{defM}, we point out that this dependence is
a consequence of the linearity of the problem \eqref{aux1}--\eqref{aux5}. {Indeed, 
e.g., if} we choose a full set of data for which $M =1$ and prove the above estimates, \elvis{then we  obtain all of the }
bounds \eqref{esti5} and \eqref{uwe6}--\eqref{uwe10} with some 
particular constants -- fully determined -- and without specification of $M$. 
Next, we can take a generic set of data for which the constant $M \, (\not=0) $ in 
\eqref{defM} is different from~$1$. Thus, pick the corresponding solution $(\mu, \varphi, \sigma)$ of \eqref{aux1}--\eqref{aux5} 
and divide all components of the triplet $(\mu, \varphi, \sigma)$ by $M$; then,  the 
scaled triplet  $(\mu/M, \varphi/M, \sigma/M)$ solves another problem in which the 
data are all divided by $M$  and satisfy \eqref{defM} with \an{constant} $1$. 
Hence, the \elvis{previously} found universal estimates work also for $(\mu/M, \varphi/M, \sigma/M)$ 
with the same constants as before. As a consequence, it is straightforward to finally 
obtain \eqref{uwe9} and \eqref{uwe10} (and previous estimates as well),
simply \elvis{by} multiplying by $M.$}  

\pcol{\Brem
\label{ANDR}
We point out that the above linearity argument has been implicitly used in the paper \cite{CSS1} (see, in particular, \cite[(3.10) and~(3.14)]{CSS1}). On the other hand, 
distinct approaches may be possible; in particular, we aim to mention the analysis \an{performed in \cite{ST}} where, for slightly more regular initial data, the authors can obtain \elvis{continuity for all components
of the solutions $(\mu, \varphi, \sigma)$, and boundedness estimates analogous to \eqref{uwe9} and \eqref{uwe10} are satisfied}.
\Erem}%

\pcol{At this point, we can combine  all the estimates \eqref{esti5}, \eqref{uwe6}--\eqref{uwe10} and deduce that}
\begin{align}
\label{esti6}
\juerg{\|(\mu,\ph,\sigma)\|_{{\cal X}}\,\le\,C_{25}M,} 
\end{align} 
\pcol{${\cal X}$ being defined in \eqref{calX}. Moreover,}
\juerg{it is readily seen from the equations \eqref{aux1}--\eqref{aux3} 
that $(\mu,\ph,\sigma)\in {\cal Y}$ \pcol{(cf.~\eqref{calY})} and that}
\begin{equation}
\juerg{\|(\mu,\ph,\sigma)\|_{{\cal Y}}\,\le\,C_{26}M.}
\end{equation}

The existence of a solution with the asserted properties is thus shown. It remains to
prove the uniqueness. To this end, let $(\mu_i,\ph_i,\sigma_i)\in\Y$, $i=1,2$, be two solutions to 
the system. Then $(\mu,\ph,\sigma):=(\mu_1,\ph_1,\sigma_1)-
(\mu_2,\ph_2,\sigma_2)$ solves the system \eqref{aux1}--\eqref{aux5} 
with zero initial data, where the terms $\lambda_2\anfirst{h}_i$,
$i=1,2$, and 
$\lambda_3f_i$, $i=1,2,3$, on the  right-hand sides do not occur. {By the definition of $\,{\cal Y}$
(recall \eqref{calX} and \eqref{calY}), and since} $(\mu,\ph,\sigma)\in\Y$,
all of the generalized partial derivatives occurring in \eqref{aux1}--\eqref{aux3} belong to $L^2(Q)$.
Therefore, we may repeat -- now for the continuous
problem -- the a priori estimates performed for the Faedo--Galerkin approximations that led us to the
estimate \eqref{esti1}. We then find analogous estimates for
$(\mu,\ph,\sigma)$, where this time the constant $\,M\,$ from \eqref{defM} equals zero. Thus, $(\mu,\ph,\sigma)
=(0,0,0)$. With this, the uniqueness is shown, which finishes the proof of the assertion.
\Edim

\vspace{2mm}
\juerg{\pcol{Having proved Lemma~\ref{LEM:FRE}}, we are in a position to prepare for the application of the implicit function
theorem.  For this purpose, let us consider two auxiliary linear initial-boundary value problems. The first,}
\begin{align}
\label{sysG11} 
&\alpha\dt\mu+\dt\ph-\Delta\mu\,= \, f_1&&\mbox{ in \,$Q$},\\[0.5mm]
\label{sysG12} 
&\beta\dt\ph-\Delta\ph-\mu\,= \, f_2 &&\mbox{ in \,$Q$},
\\[0.5mm]
\label{sysG13}
&\dt\sigma-\Delta\sigma+\chi\Delta\ph\,=\,f_3&&\mbox{ in \,$Q$},\\[0.5mm]
\label{sysG14}
&\dn\mu\,=\,\dn\ph\,=\,\dn\sigma\,=\,0 &&\mbox{ on \,$\Sigma$},\\[0.5mm]
\label{sysG15}
&\mu(0)\,=\,\ph(0)\,=\,\sigma(0)\,=0 &&\mbox{ in }\,\Omega,
\end{align}
is obtained from \eqref{aux1}--\eqref{aux5} for $\lambda_1=\lambda_2=\lambda_4=0, \, \lambda_3=1$.
Thanks to Lemma \ref{LEM:FRE}, this system has for each $(f_1,f_2,f_3) \in L^\infty(Q)\times\bigl(
H^1(0,T;H)\cap\Liq\bigr)\times\Liq$ a unique solution
$(\mu,\ph,\sigma) \in {\cal Y}$, and the associated linear mapping
\begin{equation}\label{defG1}
{\cal G}_1:{\left(L^\infty(Q)\times\bigl(H^1(0,T;H)\cap\Liq\bigr)\times\Liq\right)}\to {\cal Y}; \,\,(f_1,f_2,f_3)\mapsto (\mu,\ph,\sigma),
\end{equation}
is continuous. The second system reads
\begin{align}\label{sysG21}
&\alpha\dt\mu+\dt\ph-\Delta\mu\,=\,0 &&\mbox{ in \,$Q$},\\[0.5mm]
\label{sysG22}
&\beta\dt\ph-\Delta\ph-\mu\,=\,0 &&\mbox{ in \,$Q$},
\\[0.5mm]
\label{sysG23}
&\dt\sigma-\Delta\sigma+\chi\Delta\ph\,=\,0&&\mbox{ in \,$Q$},\\[0.5mm]
\label{sysG24}
&\dn\mu\,=\,\dn\ph\,=\,\dn\sigma\,=\,0 &&\mbox{ on \,$\Sigma$},\\[0.5mm]
\label{sysG25}
&\mu(0)\,=\,\mu_0,\quad \ph(0)\,=\,\ph_0, \quad \sigma(0)\,=\,
\sigma_0 &&\mbox{ in }\,\Omega,
\end{align}
and results from \eqref{aux1}--\eqref{aux5} for $\lambda_1=\lambda_2=\lambda_3=0,\lambda_4=1$.
For each $(\mu_0,\ph_0,\sigma_0) \in \juerg{{\mathcal{N}}}$, it also enjoys a unique solution $(\mu,\ph,\sigma) \in {\cal Y}$, and the associated mapping
\begin{equation}\label{defG2}
{\cal G}_2: \juerg{\mathcal{N}} \to {\cal Y}; \,\,(\mu_0,\ph_0,\sigma_0)
\mapsto (\mu,\ph,\sigma),
\end{equation}
is linear and continuous as well.
%\noindent 
In addition, we define on the open set ${\cal A}:=({\cal U}_R\times\Phi)
\subset ({\cal U}\times\Y)$ the nonlinear mapping
\begin{align}
\label{defG3}
&{\cal G}_3:\an{{\cal A}}\to {\left(L^\infty(Q)\times\bigl(H^1(0,T;H)\cap\Liq\bigr)\times\Liq\right)};
\nonumber\\[0.5mm]
&((u_1,u_2),(\mu,\ph,\sigma))\mapsto (f_1,f_2,f_3),
\,\,\mbox{ where}\nonumber\\[0.5mm]
&(f_1,f_2,f_3)=(P(\ph)(\sigma+\chi(1-\ph)-\mu)-\h(\ph)u_1, \chi\,\sigma-F'(\ph),\nonumber\\
&\hspace*{25mm}-P(\ph)(\sigma+\chi(1-\ph)-\mu)+u_2)\,. 
\end{align}
The solution $(\mu,\ph,\sigma)$ to the nonlinear state equation \eqref{ss1}--\eqref{ss5} is the sum
of the solution to the system \eqref{sysG11}--\eqref{sysG15}, where $(f_1,f_2,f_3)$ is chosen as above (with $(\mu,\varphi,\sigma)$ considered as known), and of the solution to the system
\eqref{sysG21}--\eqref{sysG25}. Therefore, the state vector $(\mu,\ph,\sigma)$ associated with the control vector $(u_1,u_2)$ is the unique solution to the nonlinear equation
\begin{equation} \label{nonlineq}
(\mu,\ph,\sigma) = {\cal G}_1 \big({\cal G}_3((u_1,u_2),(\mu,\ph,\sigma)\big) +  {\cal G}_2(\mu_0,\varphi_0,\sigma_0).
\end{equation}
Let us now define  the nonlinear mapping  $\,{\cal F}:{\cal A}\to \Y$,
\begin{align}
\label{defF}
 {\cal F}((u_1,u_2),(\mu,\ph,\sigma))\,=\,{\cal G}_1({\cal G}_3
((u_1,u_2),(\mu,\ph,\sigma))+{\cal G}_2(\mu_0,\ph_0,\sigma_0) - (\mu,\ph,\sigma).
\end{align} 
With $ {\cal F}$, the state equation can be shortly written as
\begin{equation} \label{nonlineq2}
{\cal F}((u_1,u_2)\anfirst{,}(\mu,\ph, \sigma))=(0,0,0).
\end{equation}
This equation just means that $(\mu,\ph,\sigma)$ is a solution to the state system \eqref{ss1}--\eqref{ss5} such that    $((u_1,u_2),(\mu,\ph,\sigma))\in{\cal A}$. From Theorem \juergen{2.3} we 
know that such a solution exists for every $(u_1,u_2)\in{\cal U}_R$. A fortiori, any such solution automatically enjoys the separation property  \eqref{ssbound2} and is uniquely determined.  

We are going to apply the implicit function theorem to the equation  \eqref{nonlineq2}. To this end,
we need the differentiability of the involved mappings.

Observe that, owing to the differentiability properties of the
involved Nemytskii operators (see, e.g., \cite[Thm.~4.22, \pcol{p.~229}]{Fredibuch}), the mapping $\,{\cal G}_3\,$ is twice continuously Fr\'echet differentiable \newju{into the space $L^\infty(Q)\times L^\infty(Q)\times L^\infty(Q)$}, and for the first partial derivatives at any point 
$\,((\uebar,\uzbar),(\bmu,\bvp,\bsigma))\in {\cal A}$, and for all $(u_1,u_2)\in{\cal U}$
and $(\mu,\ph,\sigma)\in \Y$, we have the identities
\begin{align}
\label{Freu}
&D_{(u_1,u_2)}{\cal G}_3((\uebar,\uzbar),(\bmu,\bvp,\bsigma))(u_1,u_2)\,=\,
\anfirst{(-\h(\bvp)u_1},0, u_2),\\[0.5mm]
\label{Frey}
&D_{(\mu,\ph,\sigma)}{\cal G}_3((\uebar,\uzbar),(\bmu,\bvp,\bsigma))(\mu,\ph, \sigma)\nonumber\\
&=\,
\bigl(P(\bvp)(\sigma-\chi\ph-\mu)+P'(\bvp)(\bsigma+\chi(1-\bvp)-\bmu)\ph-\h'(\bvp)\,\uebar\,\ph,\,\,\chi\sigma-F''(\bvp)\ph,
\nonumber\\
& \quad\,-P(\bvp)(\sigma-\chi\ph-\mu)-P'(\bvp)(\bsigma+\chi(1-\bvp)-\bmu)\ph\bigr)\,.
\end{align}  
\newju{It remains to show that the second component of ${\cal G}_3$ is also twice continuously Fr\'echet differentiable
into the space $H^1(0,T;H)$. But this follows exactly as in \cite[Sect.~2,~\elvis{(2.70)} ff.]{ST}. We thus can refer the reader
to \cite{ST} for this argument.}

\vspace{2mm}
\newju{At this point, we may  apply the chain rule, which yields}
 that ${\cal F}$ 
is twice continuously Fr\'echet differentiable from ${\cal U}_R\times\Phi$ into $\Y$, with the first-order partial derivatives
\begin{align}
D_{(u_1,u_2)}{\cal F}((\uebar,\uzbar),(\bmu,\bvp,\bsigma))\,&=\,{\cal G}_1\circ D_{(u_1,u_2)}{\cal 
 G}_3((\uebar,\uzbar),(\bmu,\bvp,\bsigma)),
\label{DFu}
\\[1mm]
D_{(\mu,\ph,\sigma)}{\cal F}((\uebar,\uzbar),(\bmu,\bvp,\bsigma))\,&=\,{\cal G}_1\circ
D_{(\mu,\ph,\sigma)}{\cal G}_3((\uebar,\uzbar),(\bmu,\bvp,\bsigma))-I_{\Y},
\label{DFy}
\end{align}
where $\,I_{\Y}\,$ denotes the identity mapping on $\,\Y$.

At this point, we introduce for convenience abbreviating denotations, namely,
\begin{align*}
\bu&:=(u_1,u_2),\quad \overline\bu:=(\uebar,\uzbar), \quad {\bf y}:=(\mu,\vp,\sigma),\quad
{\bf {\overline y}}:=(\bmu,\bvp,\bsigma),\nonumber\\
{\bf y}_0&:=(\mu_0,\ph_0,\sigma_0), \quad \mathbf{0}:=(0,0,0).
\end{align*}
With these denotations, we want to prove the differentiability of the control-to-state mapping $\bu \mapsto  {\bf y}$  defined implicitly by the equation $\,{\cal F}(\bu, {\bf y})=\mathbf{0}$, using the implicit function
theorem. Now let $\overline\bu\in  {\cal U}_R$ be given and $\overline{\bf y}={\cal S}(\overline\bu)$. We need to show that the linear
and continuous operator $\,D_{{\bf y}}{\cal F}(\overline\bu,\overline{\bf y})$ is a topological isomorphism from $\Y$ into itself. 

To this end, let ${\bf v}\in\Y$ be arbitrary. Then the identity $\,D_{\bf y}{\cal F}(\overline\bu,
\overline{\bf y})({\bf y})={\bf v}\,$ just means that $\,{\cal G}_1\left(D_{\bf y}{\cal G}_3(\overline\bu,
\overline{\bf y})(\bf y)\right)-{\bf y}={\bf v}$, which is equivalent to saying that   
\begin{equation*}
{\bf w}\,:=\,{\bf y}+{\bf v}={\cal G}_1\left(D_{\bf y}{\cal G}_3(\overline\bu,
\overline{\bf y})(\bf w)\right)-{\cal G}_1\left(D_{\bf y}{\cal G}_3(\overline\bu,
\overline{\bf y})(\bf v)\right).
\end{equation*}
The latter identity means that ${\bf w}$ is a solution to \eqref{aux1}--\eqref{aux5} for $\lambda_1=\lambda_3=1,
\lambda_2=\lambda_4=0$, with the specification $(f_1,f_2,f_3)=-{\cal G}_1\left(D_{\bf y}{\cal G}_3(\overline\bu,
\overline{\bf y})(\bf v)\right)\in\Y$. By \juergen{Lemma 4.1}, such a solution ${\bf w}\in\Y$ exists and is uniquely 
determined. We thus can infer that $D_{{\bf y}}{\cal F}(\overline\bu,\overline{\bf y})$ is surjective. At the same time, taking ${\bf v}=\mathbf{0}$, we see that the equation $D_{{\bf y}}{\cal F}(\overline\bu,\overline{\bf y})({\bf y})=\mathbf{0}$
means that ${\bf y}$ is the unique solution to \eqref{aux1}--\eqref{aux5} for $\lambda_1=1,\lambda_2=\lambda_3=
\lambda_4=0$. Obviously, ${\bf y}=\mathbf{0}$, which implies that $D_{{\bf y}}{\cal F}(\overline\bu,\overline{\bf y})$ is
also injective and thus, by the open mapping principle, a topological isomorphism from $\Y$ into itself. 

At this point, we may employ the implicit function theorem (cf., e.g., \cite[Thms. 4.7.1 and 5.4.5]{cartan1967} or \cite[10.2.1]{Dieu}) to conclude    
that the mapping $\S$ is twice continuously Fr\'echet differentiable from ${\cal U}_R$ into $\Y$ and that the
first Fr\'echet derivative $\,D\S(\ubar)\,$ of $\,\S\,$ at $\overline\bu\in{\cal U}_R$ is given by the formula
\begin{equation}
\label{DS1}
D{\cal S}(\overline\bu)\,=\,-D_{\bf y}{\cal F}(\overline\bu,\overline{\bf y})^{-1}\circ D_{\bu}
{\cal F}(\overline\bu,\overline{\bf y}).
\end{equation}
Now let $\anfirst{{\bf h}=(h_1,h_2)}\in{\cal U}$ be arbitrary and ${\bf y}=(\mu,\ph,\sigma)=D{\cal S}(\overline\bu)({\anfirst{\bf h}})$. 
Then, $$ D_{\bf y}{\cal F}(\overline\bu,\overline{\bf y})({\bf y})=- \an{D}_{\bu}
{\cal F}(\overline\bu,\overline{\bf y})({\anfirst{\bf h}}),$$ which is obviously equivalent to saying that
$$
{\bf y}={\cal G}_1\left(D_{\bf y}{\cal G}_3(\overline\bu,\overline{\bf y})({\bf y})\right)
+{\cal G}_1(\juerg{-}\h(\bvp)\anfirst{h_1},0,\anfirst{h_2}).
$$
This, in turn, means that ${\bf y}$ is the unique solution to the problem \eqref{aux1}--\eqref{aux5} for
$\lambda_1=\lambda_2=1,\lambda_3=\lambda_4=0$.

In summary, we have shown the following result.

\Bthm[\Frechet\ differentiability of $\S$]
\label{THM:FRECHET}
Suppose that the conditions \ref{const:weak}, \ref{F:strong:1}--\juergen{\ref{Ph:strong},}  
and \eqref{defUR} are fulfilled.
Moreover, let the initial data $(\m_0,\ph_0,\s_0)$ verify \eqref{strong:initialdata} and \eqref{strong:sep:initialdata},
and let $\overline \bu=(\uebar,\uzbar)\in {\cal U}_R$ be arbitrary and $(\bmu,\bvp,\bsigma)={\cal S}(\overline\bu)$. Then the control-to-state 
operator $\S$ is twice continuously Fr\'echet differentiable at $\,\overline\bu\,$ as a mapping from $\,{\cal U}\,$ into 
$\,{\cal Y}$. Moreover, for every $\anfirst{ \bh=(h_1,h_2)}\in {\cal U}$, the Fr\'echet derivative $\,D \S(\overline\bu)\in 
{\cal L}({\cal U},\Y)\,$ of $\,\S\,$ at
$\,\overline\bu\,$ is given by the identity $\,D\S(\overline\bu)\anfirst{(\bh)}=(\mu,\ph,\sigma)$, where 
$(\mu,\ph,\sigma)$ is the 
unique solution to the linear system \eqref{aux1}--\eqref{aux5} with $\lambda_1=\lambda_2=1,\lambda_3=\lambda_4=0$.
\Ethm

\anfirst{Motivated by the forthcoming analysis, \juerg{we} now present a stability estimate for the solutions to the linearized system. \juerg{In abuse of notation, we will} denote the \juergen{linearized} variable associated with $\ph$ \juerg{by} $\xi$, which up to now was devoted to indicate a selection of the subdifferential $\partial F_1$ evaluated at some point. Since the optimal control problem demands to work with strong solutions, we no longer have any selection to work with, so that from now on the variable $\xi$ will play a different role (cf. Theorem \ref{THM:STRONG}).} \newju{We have the
following result, where we recall the definition \eqref{defV} of the space ${\cal V}$.}
\begin{theorem}
\label{THM:CONTDEP:LIN}
Suppose that the conditions \ref{const:weak}, \ref{F:strong:1}--\an{\ref{Ph:strong}}, \ref{ass:control:Uad} and 
\eqref{defUR} are fulfilled.
Moreover, let the initial data $(\m_0,\ph_0,\s_0)$ \juergen{satisfy} \an{\eqref{strong:initialdata}--\eqref{strong:sep:initialdata}},
and let $\ov \bu^i=(\uebar^i,\uzbar^i)\in \UR$ be arbitrary and $(\bmu_i,\bvp_i,\bsigma_i)={\cal S}(\overline\bu^i)$, \anfirst{$i= 1,2 $.}
Furthermore, let $(\eta_i,\xi_i,\theta_i)$ denote the associated solutions to the linearized system (i.e.\pcol{, the} system \eqref{aux1}--\eqref{aux5} with $\lambda_1=\lambda_2=1,\lambda_3=\lambda_4=0$).
Then the mapping $D\S$ is \Lip\ continuous on $\UR$ in the sense that there exists a positive constant $\anfirst{K_d}$ such that for all $\bh \in \juergen{\cal U}$ we have
\begin{align}
	\label{fre:cd:estimate}
\juergen{	\norma{\big(D\S(\ov \bu^1)-D\S(\ov \bu^2)\big)\anfirst{(\bh)}}_\newju{{\cal V}} \leq \anfirst{K_d}\, \norma{\ov\bu^1-\ov\bu^2}_{\L2 H^2}\norma{\bh}_{\L2 H^2}.}
\end{align}
\end{theorem}
\begin{proof}
Due to Theorem \ref{THM:FRECHET}, the proof of \eqref{fre:cd:estimate} reduces to showing
 that there exists a constant $c>0$
such that
\begin{align}
	\label{est:cd:lin}
	\norma{(\eta_1,\xi_1,\th_1)-(\eta_2,\xi_2,\th_2)}_\newju{{\cal V}} \leq 	c \,
	\norma{\ov \bu^1-\ov \bu^2}_{\L2 H ^2}\norma{\bh}_{\L2 H ^2},
\end{align}
which is the estimate we are going to check.
Moreover, let us notice that \eqref{est:rem} in Remark \ref{REM:L2} guarantees the existence of a positive constant $c$ such 
that
\begin{align}
	\label{est:lin:cont}
	\norma{(\eta_i,\xi_i,\theta_i)}_\Z \leq c \,\norma{\bh}_{\L2 H ^2}, \quad \hbox{ $i=1,2$.}
\end{align}
Next, we set
\begin{align*}
	 \bm &= \bm_1-\bm_2, \quad 
	\bph = \bph_1-\bph_2, \quad
	 \bs= \bs_1-\bs_2, \quad
	\\  \eta &= \eta_1-\eta_2, \quad 
	\xi = \xi_1-\xi_2, \quad
	 \th= \th_1-\th_2,
\end{align*}
and observe that the triple $( \eta,  \xi, \th)$ solves the system obtained from taking the difference 
between the linearized system written for $(\eta_1,\xi_1,\theta_1)$
and $(\eta_2,\xi_2,\theta_2)$, which reads 
\begin{align}
	\label{cd:fre:1}
	& \a \dt  \eta + \dt  \xi - \Delta  \eta =  f_1 +  f_2 && \hbox{in $Q$,}
	\\
	\label{cd:fre:2}
	& \b \dt \xi - \Delta  \xi + (F''(\bph_1)-F''(\bph_2))\xi_1 + F''(\bph_2) \xi =  \eta + \chi  \th && \hbox{in $Q$,}
	\\
	\label{cd:fre:3}
	& \dt  \th - \Delta  \th \anold{+ \chi \Delta  \xi } =-  f_1 && \hbox{in $Q$,}
	\\
	\label{cd:fre:4}
	& \dn  \eta =\dn  \xi= \dn  \th = 0 && \hbox{on $\Sigma$,}
	\\
	\label{cd:fre:5}
	&  \eta(0)= \xi(0)=  \th(0)=0 && \hbox{in $\Omega$,}
\end{align}
where \pcol{now $f_1$ and $f_2$ are specified by}
\begin{align*}
		 f_1 & = 
		(P(\bph_1)-P(\bph_2))(\th_1-\chi \xi_1-\eta_1)
		+ P(\bph_2)( \th - \chi  \xi -  \eta)
		\\ & \quad
		+ (P'(\bph_1)-P'(\bph_2))\xi_1(\bs_1+\chi(1-\bph_1)-\bm_1)
		+ P'(\bph_2) \xi(\bs_1+\chi(1-\bph_1)-\bm_1)
		\next
		+ P'(\bph_2)\xi_2( \bs - \chi  \bph -  \bm),
		\\ 
		 f_2 & = 
		- (\h'(\bph_1)-\h'(\bph_2)) \xi_1 \ov u_1^1
		- \h'(\bph_2)\,  \xi \, \ov u_1^1
		- \h'(\bph_2)\,\xi_2 \, ({\ov u_1^1-\ov u_1^2})
		\anfirst{- (\h (\bph_1)-\h (\bph_2))h_1}.
\end{align*}
Moreover, let us recall that due to \eqref{cont:dep:strong}\anfirst{, we have the stability estimate}
\begin{align}
	\norma{\bm_1-\bm_2}_{\newju{Z}}
	+ \norma{\bph_1-\bph_2}_{\newju{Z}} + \norma{\bs_1-\bs_2}_{\newju{Z}}
	\,	\leq\,
	c \,\norma{\ov \bu^1-\ov \bu^2}_{\L2 H ^2}.
	\label{cons}
\end{align}

We now aim at deriving some a priori estimates for the differences $(\eta,\xi,\th)$ which will \andrea{entail \eqref{est:cd:lin}.
\juergen{Prior to this}, let us premise a general fact that will be employed several times later on. To this end,
let $f: \erre \to \erre$ be a regular, bounded and \Lip\ continuous function
with \Lip\ continuous and bounded derivative $f'$ (in what follows 
the role of $f$ will be played by $P$, $\h$, $F$, and possibly their derivatives).
Let $\ph_1$ and $\ph_2$ be the second components of different strong solutions $(\m_1,\ph_1,\s_1)$ and $(\m_2,\ph_2,\s_2)$ associated \juergen{with controls $\bu^1,\bu^2\in {\cal U}_R$ according to Theorem \ref{THM:STRONG}. Then it follows from the 
continuity of the \pcol{embeddings $V\subset L^p(\Omega)$ and $W_0 \subset W^{1,p} (\Omega)$, $p \in [1,6]$,} and the estimate \eqref{ssbound1} that, a.\,e.
in $(0,T)$,}
\begin{align}
\juergen{\|f(\ph_1)-f(\ph_2)\|_\anold{p}}&\,\juergen{\le\,c\,}	\norma{f(\ph_1)-f(\ph_2)}_V\non\\[0.2mm]
	&  = \,c\,\norma{f(\ph_1)-f(\ph_2)} + c\,\norma{f'(\ph_1) \nabla \ph_1-f'(\ph_2)\nabla \ph_2}\non
	\\[0.2mm]
	 &	\leq 	\non
	c \,\norma{\ph_1-\ph_2}
	+ c \,\norma{(f'(\ph_1)-f'(\ph_2))\nabla \ph_1 + f'(\ph_2)\nabla (\ph_1-\ph_2) }
	\\[0.2mm]
	& \leq 
	\non
	c \,\big( \norma{\ph_1-\ph_2} +\pcol{\norma{\ph_1-\ph_2}_4\,\norma{\nabla \ph_1}_4}
			 + \norma{\nabla (\ph_1-\ph_2)} \big)
	\\[0.2mm] 
	& \leq 
	\juergen{c \,\big(1+K_1\big) \,\norma{\ph_1-\ph_2}_V.}
	\label{normaVdiff}
\end{align}}%
In addition, for the sake of a shorter exposition, in the upcoming estimates 
\pcol{we will} avoid to explicitly write
the integration variable
$\,s\,$ in the time integrals.

\noindent
{\sc First estimate:}
We test \eqref{cd:fre:1} by $ \eta$, \eqref{cd:fre:2}, to which we add $\xi$ on both sides, by $\dt  \xi$,
and \eqref{cd:fre:3} by $ \th$, add the resulting equalities and integrate over time, obtaining
\begin{align*}
	&\frac\a2 \IO2 { \eta }
	+ \I2 {\nabla  \eta }
	+ \b\I2 {\dt  \xi}
	+ \frac 12 \norma{ \xi(t)}^2_V
	+ \frac 12 \IO2 { \th }
	+ \I2 {\nabla  \th}
	\\  & \quad
	= \intQt  f_1 ( \eta -  \th) + \intQt  f_2  \eta 
	- \intQt (F''(\bph_1)-F''(\bph_2))\xi_1 \dt  \xi
	- \intQt  F''(\bph_2) \xi\dt  \xi
	\\ & \qquad
	+ \intQt  ( \eta + \chi  \th +  \xi)\dt  \xi
	+ \chi \intQt \nabla  \xi \cdot \nabla  \th \,= : \,I_1+...+I_6.
\end{align*}
Using the Young and \Holder\ inequalities, the \Lip\ continuity of $P, \, P'$ and $\h'$, the continuous embedding $V \subset \Lx4$, the uniform bounds \eqref{ssbound1}--\eqref{ssbound3} for $(\bm_i,\bph_i,\bs_i)$, $i=1,2$,
{as well as \eqref{est:lin:cont},} \eqref{cons}, and \eqref{normaVdiff}, we infer that
\begin{align*}
	I_1 & \leq 
	c \intQt (| \eta|^2+| \xi|^2+|  \th|^2)
	+ c \iot \norma{P(\bph_1)-P(\bph_2)}^2_4\,(\norma{\th_1}^2_\anfirst{4}+\norma{\xi_1}^2_\anfirst{4}+\norma{\eta_1}^2_\anfirst{4})
\juergen{\,ds}
	\next
	+ c \, (\norma{\bs_1}_{L^\infty(Q)}^2+\norma{\bph_1}_{L^\infty(Q)}^2+\norma{\bm_1}_{L^\infty(Q)}^2+1)\iot \norma{P'(\bph_1)-P'(\bph_2)}^2_4\,\norma{\xi_1}_{4}^2\juergen{\,ds} 
	\next
	+ c\, (\norma{\bs_1}_{L^\infty(Q)}^2+\norma{\bph_1}_{L^\infty(Q)}^2+\norma{\bm_1}_{L^\infty(Q)}^2+1) \I2 {\xi}
	\next
	+ c  \iot (\norma{ \bs_1-\bs_2}_4^2+\norma{ \bph_1-\bph_2}_4^2+\norma{ \bm_1-\bm_2}_4^2) \,\norma{\xi_2}_{4}^2
	\juergen{\,ds}
	\\ & \leq
	c \intQt (| \eta|^2+| \xi|^2+|  \th|^2)
	\next 
	+ c\, (\norma{\th_1}^2_{\L\infty V}+\norma{\xi_1}^2_{\L\infty V}+\norma{\eta_1}^2_{\L\infty V})\iot \norma{\bph_1-\bph_2}^2_V\juergen{\,ds} 
	\next
	+ c \,\norma{\xi_1}_{\L\infty V}^2\,  (\norma{\bs_1}_{L^\infty(Q)}^2+\norma{\bph_1}_{L^\infty(Q)}^2+\norma{\bm_1}_{L^\infty(Q)}^2+1)\iot \norma{\bph_1-\bph_2}^2_V \juergen{\,ds}               
	\next
	+ c\, (\norma{\bs_1}_{L^\infty(Q)}^2+\norma{\bph_1}_{L^\infty(Q)}^2+\norma{\bm_1}_{L^\infty(Q)}^2+1) \I2 {\xi}
	\next
	+ c \, \norma{\xi_2}_{\L\infty V}^2 \iot (\norma{ \bs_1-\bs_2}_V^2+\norma{ \bph_1-\bph_2}_V^2+\norma{ \bm_1-\bm_2}_V^2)
	\juergen{\,ds}
%\end{align*}
%\begin{align*}	
   \\
  & \leq
	c \intQt (| \eta|^2+| \xi|^2+|  \th|^2)
	+ c \,\norma{\ov \bu^1 -\ov \bu^2}_{\L2 H ^2}^2\, \norma{\bh}_{\L2 H ^2}^2.
\end{align*}
\pcol{Next}, by similar computations, we deduce that, \juergen{for any $\delta>0$ (to be chosen later),}
\begin{align*}
	I_2 & \leq 
	c \intQt | \eta|^2
	+ c\, \norma{\ov u_1^1}^2_{L^\infty(Q)}\iot \norma{\h'(\bph_1)-\h'(\bph_2)}^2_4\, \norma{\xi_1}^2_{4}\juergen{\,ds}
	+ c\,\norma{\ov u_1^1}^2_{L^\infty(Q)} \I2 { \xi} \non\\
	&\quad + 2\d \iot \norma{\eta}_{\pcol{V}}^2\juergen{\,ds}
		+ \cd \iot \juergen{\norma{\ov \bu^1- \ov \bu^2}^2}\,\norma{\xi_2}^2_{4}\juergen{\,ds}
	+ \cd \iot \norma{\h(\bph_1)-\h(\bph_2)}_4^2\, \norma{h_1}^2\juergen{\,ds}
	\\ & 
	\leq
	(c+2\d) \intQt | \eta|^2
	+ c \,\norma{\xi_1}^2_{\L\infty V}\,\norma{\ov u_1^1}^2_{L^\infty(Q)}\iot \norma{\bph_1-\bph_2}^2_V
	\juergen{\,ds}               
	\next
	+ c\,\norma{\ov u_1^1}^2_{L^\infty(Q)} \I2 { \xi}
	+ 2\d \I2 {\nabla \eta}
	\next
	+ \cd \,\norma{\xi_2}^2_{\L\infty V }\iot \juergen{\norma{\ov \bu^1- \ov \bu^2}^2\,ds}
	+ \cd \iot \norma{\bph_1-\bph_2}^2_V\,\juergen{ \norma{h_1}^2\,ds}
	\\ &\leq 
	2\d \intQt | \nabla \eta|^2
	+\cd \intQt | \eta|^2
	+c  \intQt | \xi|^2
	+ \cd \,\norma{\ov \bu^1 -\ov \bu^2}_{\L2 H ^2}^2\,\norma{\bh}_{\L2 H^2}^2\,.
\end{align*}
\juergen{The terms involving the potentials can be easily handled by invoking} the separation principle \eqref{ssbound2}, which entails the \Lip\ continuity of $F$ and of its derivatives\juergen{. Using this, \eqref{est:lin:cont},} 
 \eqref{normaVdiff}, and the Young inequality\juergen{, we obtain that}
\begin{align*}
	I_3+I_4 & \leq 
	\d \I2 {\dt  \xi}
	+ \cd\, \norma{\xi_1}^2_{\L\infty V} \iot\norma{\bph_1-\bph_2}^2_V\juergen{\,ds} 
	+ \cd \I2 \xi
	\\ & \leq
	\d \I2 {\dt  \xi}
	+ \cd  \,\norma{\ov \bu^1 -\ov \bu^2}_{\L2 H ^2}^2\,\norma{\bh}_{\L2 H ^2}^2	
	+ \cd \I2 { \xi}.
\end{align*}
Finally, we have that
\begin{align*}
	I_5+I_6 & \leq 
		\d \intQt (|{\dt  \xi}|^2+| {\nabla  \th}|^2)
		+ \cd \intQt ( | \eta|^2+ | \theta|^2+ |  \xi|^2+ |\nabla  \xi|^2).
\end{align*}
\juergen{At this point, we collect the above estimates and adjust $\d\in (0,1)$ small enough. Gronwall's lemma 
then yields that}
\begin{align}\label{paolauno}
	& \norma{ \eta }_{\L\infty H \cap \L2 V}
	+ \norma{ \xi }_{\H1 H \cap \L\infty V}
	+ \norma{ \theta }_{\L\infty H \cap \L2 V} \non
	\\[0.3mm] 
	&\quad 	\leq c\, \norma{\ov \bu^1 -\ov \bu^2}_{\L2 H ^2}\,\norma{\bh}_{\L2 H ^2}.
\end{align}

\noindent
{\sc Second estimate:}
{Estimate \eqref{paolauno} entails} that the term $\dt  \xi $ is bounded in $\L2 H$ by the \juergen{expression}
 on the \rhs\ of the above estimate. \newju{Hence, equation \eqref{cd:fre:2} can be expressed as an elliptic equation 
for $\xi$ whose \rhs\ is bounded in $\L2 H$
by the same {expression. Therefore}, we easily get from elliptic regularity theory that}
%\begin{align}\label{paoladue}
%	\norma{ \eta }_{\H1 H \cap \L\infty V \cap \L2 {W_0}}
%	\leq c \,\norma{\ov \bu^1 -\ov \bu^2}_{\L2 H ^2}\,\norma{\bh}_{\L2 H ^2}.
%\end{align}
%
%\noindent 
%{\sc Third estimate:}
%Expressing \eqref{cd:fre:2} as an elliptic equation for $ \xi$, we easily get from elliptic regularity theory that
\begin{align}\label{paolatre}
	\norma{ \xi}_{\L2 {W_0}}\leq c \,\norma{\ov \bu^1 -\ov \bu^2}_{\L2 H ^2}\,\norma{\bh}_{\L2 H ^2}.
\end{align}
\newju{ This concludes the proof of Theorem \ref{THM:CONTDEP:LIN}.}
\end{proof}

\andrea{
\Brem
\label{REM:FRE}
Let us point out that \anfirst{Theorem \ref{THM:FRECHET}} establish\anfirst{es} the \Frechet\ differentiability of $\S$ 
\juergen{as a mapping from $L^\infty(Q)^2$ into $\Y$, a space of very regular functions. However,
the \Frechet\ differentiability can also be directly obtained as a mapping from $\L2 H^2$ into a space of less regular 
functions. Indeed, a closer look at the proof of \cite[Theorem~2.5]{S} reveals that the line of argumentation employed 
there can straightforwardly be adapted to our present situation, yielding that in the notation used there it holds that}
\begin{align}
	\label{freI:estimate}
	\norma{(\bm^\bh-\bm-\eta^\bh, \bph^\bh - \bph - \xih,\bs^\bh-\bs-\th^\bh )}_\Z \leq c\, \norma{\bh}_{\L2H^2}^2,
\end{align}
\juergen{which in turn} entails that the control-to-state operator $\S$ is \Frechet\ differentiable as a mapping 
from $\L2 H ^2$ into $\Z \,\juergen{\supset}\, \Y$.
We have chosen to not follow this \juergen{approach because, although it will suffice to handle
the first-order necessary conditions pointed out below in Section \ref{SEC:FOC},
it would not allow us to deal with} the second-order sufficient conditions established in Section \ref{SEC:SOC}.
\Erem         
}

Motivated by the forthcoming analysis on the second-order sufficient optimality conditions, let us now also 
explicitly identify the second-order \Frechet\ derivative of the control-to-state operator $\S$, \juergen{which exists according to Theorem \ref{THM:FRECHET}. To this end, let $\ov\bu\in{\cal U}_R$ be given. For arbitrary increments $\bh,\bk \in \Uh$, we set}
\begin{align*}
	(\bm,\bph,\bs) \an{:=} \S(\ov \bu), \quad (\eta^\bh,\xih,\theta^\bh) \an{:=} D\S(\ov \bu)\anfirst{(\bh)}, \quad (\eta^\bk,\xik,\theta^\bk)\an{: =} D\S(\ov \bu)\anfirst{(\bk)},
\end{align*}
{where it is known that $(\eta^\bh,\xih,\theta^\bh),(\eta^\bk,\xik,\theta^\bk)\in{\mathcal{Y}}$. Now, if one 
adapts the argumentation of the proof of \cite[Theorem~5.16,~pp.~288--289]{Fredibuch} to the present situation, starting 
from the identities \eqref{nonlineq}--\eqref{DFy}, then one concludes that the second-order Fr\'echet derivative 
\juergen{$D^2{\cal S}(\ov\bu)(\bk)(\bh)$} can be evaluated using the system \eqref{bilin:1}--\eqref{bilin:5} introduced below. 
Since we intend to give an independent proof, we do not give the details
of the argument, here. We begin our analysis with the following result. }
\begin{theorem}
\label{THM:BILIN}
Assume that \ref{const:weak}, \ref{F:strong:1}--\juergen{\ref{Ph:strong}, 
\pcol{\eqref{defUR}--%
%, \eqref{strong:initialdata}, and 
\eqref{strong:sep:initialdata} are fulfilled}}. Then the following initial-boundary value problem, \juergen{which
will be referred to as the ``bilinearized system'', admits a unique solution $(\nu,\psi,\rho)\in{\cal Y}$:}
%\begin{align*}
%	(\nu,\psi,\rho) \in \Z,
%\end{align*}
%and \anfirst{$(\nu,\psi,\rho)$} verifies
\begin{align}
\label{bilin:1}
&\alpha\dt\nu+\dt\psi-\Delta\nu\,= g_1 + g_2 &&\mbox{ in \,$Q$},\\[1mm]
\label{bilin:2}
&\beta\dt\psi-\Delta\psi-\nu = \chi \rho  - F''(\bvp)\psi - F^{(3)}(\bph) \xih \xik &&\mbox{ in \,$Q$},
\\[1mm]
\label{bilin:3}
&\dt\rho-\Delta\rho+\chi\Delta\psi\,=\, - g_1 &&\mbox{ in \,$Q$},\\[1mm]
\label{bilin:4}
&\dn\nu\,=\,\dn\psi\,=\,\dn\rho\,=\,0 &&\mbox{ on \,$\Sigma$},\\[1mm]
\label{bilin:5}
&\nu(0)\,=\, \psi(0)\,=\, \rho(0)\,=\,0  &&\mbox{ in }\,\Omega,
\end{align}
\Accorpa\Bilin {bilin:1} {bilin:5}
where \pcol{$g_1$ and $g_2$ are defined as}
\begin{align}
	g_1 \non & = P(\bph) (\rho - \chi \psi -\nu)
		+ P'(\bph)\,\xik\,(\theta^\bh-\chi\xih-\eta^\bh)
		\\ & \quad \non
		+ P''(\bph)\,\xik\,\xih\,(\bs + \chi (1-\bph) - \bm) 
		+ P'(\bph)\, \psi\,(\bs + \chi (1-\bph) - \bm) 
		\\ & \quad \label{g1}
		+ P'(\bph)\,\xih \,(\theta^\bk - \chi \xik - \eta^\bk),
		\\
		g_2 & = - \h''(\bph)\,\xik\,\xih\, \ov u_1
		\anold{- \h '(\bph) \xih k_1}
		- \h'(\bph)\,\psi\, \ov u_1
		- \h'(\bph)\,\xik\, h_1.
		\label{g2}
\end{align}
\newju{Moreover, there is some constant $D>0$, which only depends on the data of the system and $R$, such that}
\begin{align}\label{bsderi}
\newju{\|(\nu,\psi,\rho)\|_{\cal V}\,\le\,D\,\|\bh\|_{L^2(0,T;H)^2}\,\|\bk\|_{L^2(0,T;H)^2} \quad\forall\,
\bh,\bk\in{\cal U}.}
\end{align}
\end{theorem}
\Bdim
{Writing $(\mu,\varphi,\sigma)$ in place of $(\nu,\psi,\rho)$, we see that the system \eqref{bilin:1}--\eqref{g2} is 
of the form \eqref{aux1}--\eqref{aux5} with the specifications $\lambda_1=\lambda_3=1, \lambda_2=\lambda_4=0$, and}
\begin{align*}
\juergen{f_1}\,&\juergen{=\,P'(\bph)\,\xik\,(\theta^\bh-\chi\xih-\eta^\bh)
				+ P''(\bph)\,\xik\,\xih\,(\bs + \chi (1-\bph) - \bm)}\\ 
				&\quad\,\juergen{+\, P'(\bph)\,\xih \,(\theta^\bk - \chi \xik - \eta^\bk)
		- \h''(\bph)\,\xik\,\xih\, \ov u_1
				\anold{- \h '(\bph) \xih k_1}
				- \h'(\bph)\,\xik\, h_1,           }\\
\juergen{f_2}\,&\juergen{=\,-F^{(3)}(\ov\varphi)\,\xih\,\xik,}\\
\juergen{f_3}\,&\juergen{=\,-\,P'(\bph)\,\xik\,(\theta^\bh-\chi\xih-\eta^\bh)\,
				-\, P''(\bph)\,\xik\,\xih\,(\bs + \chi (1-\bph) - \bm)}\\ 
				&\quad\,\juergen{-\, P'(\bph)\,\xih \,(\theta^\bk - \chi \xik - \eta^\bk).}
\end{align*}
{Since $(\eta^\bh,\xih,\theta^\bh),(\eta^\bk,\xik,\theta^\bk)\in{\mathcal{Y}}$, it is easily seen that
$(f_1,f_2,f_3)\in\Liq\times(H^1(0,T;H)\linebreak\cap\Liq)\times\Liq$, and so the existence result follows from Lemma 4.1.} 

%{\sc preliminary estimates:}
%By using the Young and \Holder\ inequalities as well as the regularity of the fixed state variables $(\bm,\bph,\bs)$
%and of the linearized variables $(\eta^\bh,\xi^\bh,\th^\bh)$ and $(\eta^\bk,\xi^\bk,\th^\bk)$, we have
%\begin{align*}
%	\norma{g_1}_{\L2 H} & \leq
%	c (\norma{\rho}_{\L2 H} +\norma{\psi}_{\L2 H} +\norma{\nu}_{\L2 H} )
%	\next 
%	+ c \norma{\xik}_{L^\infty(Q)}(\norma{\th^\bh}_{L^\infty(Q)} +\norma{\xi^\bh}_{L^\infty(Q)}+\norma{\eta^\bh}_{L^\infty(Q)})
%	\next 
%	+ c \norma{\xik}_{L^\infty(Q)}\norma{\xih}_{L^\infty(Q)} %(\norma{\bs}_{L^\infty(Q)}+\norma{\bph}_{L^\infty(Q)}+\norma{\bm}_{L^\infty(Q)}+1)
%	\next 
%	+ c \norma{\psi}_{\L2 H} (\norma{\bs}_{L^\infty(Q)}+\norma{\bph}_{L^\infty(Q)}+\norma{\bm}_{L^\infty(Q)}+1)
%	\next
%	+ c \norma{\xih}_{L^\infty(Q)}(\norma{\th^\bk}_{L^\infty(Q)} +\norma{\xi^\bk}_{L^\infty(Q)}+\norma{\eta^\bk}_{L^\infty(Q)})
%	\\ & \leq
%	c (\norma{\nu}_{\L2 H} +\norma{\psi}_{\L2 H} +\norma{\rho}_{\L2 H} +1),
%	\\
%	\norma{g_2}_{\L2 H} & \leq
%	c \norma{\xik}_{L^\infty(Q)}	\norma{\xih}_{L^\infty(Q)}\norma{\ov u_1}_{L^\infty(Q)}
%	+ c \norma{\psi}_{\L2 H} \norma{\ov u_1}_{L^\infty(Q)}
%	\next
%	+ c \norma{\xik}_{L^\infty(Q)}\norma{h_1}_{L^\infty(Q)}
%	\\ & \leq 
%	c (\norma{\psi}_{\L2 H}  +1).
%\end{align*}	
%With this information at disposal we can now easily proceed by performing some a priori estimates.
%
%\noindent

\newju{It remains to show \eqref{bsderi}. To this end, we first add the term $\psi$ on both sides of 
\eqref{bilin:2} and then test \eqref{bilin:1} by $\nu$, \eqref{bilin:2} by $\dt \psi$, and 
\eqref{bilin:3} by $\rho$. Adding the resulting equalities, and integrating over time and by parts, we infer that}
\begin{align}\label{vino}
	&\newju{\frac\a2 \IO2 \nu 
	+ \I2 {\nabla \nu}
	+ \b\I2 {\dt \psi}
	+ \frac 12 \norma{\psi(t)}^2_V
	+ \frac 12 \IO2 \rho 
	+ \I2 {\nabla \rho}}\non
	\\  &\newju{ \quad
	= \intQt g_1(\nu-\rho) + \intQt (\chi \rho  - F''(\bvp)\psi - F^{(3)}(\bph) \xih \xik + \psi)\dt\psi}\non
	\next \quad\newju{
	+ \chi \intQt \nabla \psi \cdot \nabla \rho +\intQt g_2 \nu =: I_1+I_2+I_3+I_4,}
\end{align}
\newju{with obvious meaning. Now, using the global estimates \eqref{ssbound1} and \eqref{ssbound3}, the continuity
of the embedding $V\subset L^4(\Omega)$, as well as H\"older's 
inequality and the estimate \eqref{est:rem} for the linearized variables, we can easily check that}
%\newju{so that}
\begin{align}\label{vinouno}
\newju{\an{I_1}\,}&\newju{\le\,c\intQt(|\nu|^2+|\psi|^2+|\rho|^2)\,+\,c\int_0^t\|\xi^\bk\|_4^2 
\,\bigl(\|\eta^\bh\|_4^2+\|\xi^\bh\|_4^2+\|\theta^\bh\|_4^2\bigr)\,ds}\non\\
&\newju{\quad + \,c\int_0^t \|\xi^\bk\|_4^2\,\|\xi^\bh\|_4^2\,ds\,+
c\int_0^t\|\xi^\bh\|_4^2 
\,\bigl(\|\eta^\bk\|_4^2+\|\xi^\bk\|_4^2+\|\theta^\bk\|_4^2\bigr)\,ds}\non\\
&\newju{\le\,c\intQt(|\nu|^2+|\psi|^2+|\rho|^2)\,+\,c\,\|\bh\|^2_{L^2(0,\pcol{T};H)\an{^2}}\,\|\bk\|^2_{L^2(0,\pcol{T};H)\an{^2}}\,.}
\end{align}
\newju{Arguing similarly, where we also invoke Young's \pcol{inequality}, we obtain 
%for every $\delta>0$ 
that}
\begin{equation}\label{vinodue}                
\newju{\an{I_2}\,\le\,\pcol{\frac \beta 2}\intQt|\dt\psi|^2\,+\,\pcol c\intQt(|\psi|^2+|\rho|^2)\,+\,\pcol c\,
\|\bh\|^2_{L^2(0,\pcol{T};H)^2}\,\|\bk\|^2_{L^2(0,\pcol{T};H)^2},}
\end{equation}
\newju{as well as}
\begin{equation}\label{vinotre}
\newju{\an{I_3}\,\le\,\frac 12 \intQt|\nabla\rho|^2\,+\,\frac {\chi^2}2\intQt|\nabla\psi|^2.}
\end{equation}
\newju{Finally, using also that $\ov u_1$ is bounded, we see that}
\begin{align}\label{vinoquattro}
\newju{\an{I_4}}\,&\newju{\le\,c\intQt(|\nu|^2+|\psi|^2)\,\pcol{{}+\,
c\int_0^t \|\xi^\bk\|_4^2\,\|\xi^\bh\|_4^2\,ds}} \non\\
&\quad\newju{{}+ c\int_0^t\|\nu\|_4(\|\xi^\bh\|_4\|k_1\|+\|\xi^\bk\|_4\|h_1\|)\,ds}\non\\
&\newju{\le\,(c+\an{2}\delta)\intQt(|\nu|^2+|\psi|^2)\,+\,\an{2}\delta\intQt|\nabla\nu|^2\,+\,c_\delta\,\|\bh\|^2_{L^2(0,\pcol{T};H)^2}\,\|\bk\|^2_{L^2(0,\pcol{T};H)^2}.}
\end{align}
\newju{Combining \eqref{vino}--\eqref{vinoquattro}, and choosing $\delta>0$ appropriately small, we thus obtain from
Gronwall's lemma that} 
\begin{align}
&\newju{\|\nu\|_{L^\infty(0,T;H)\cap L^2(0,T;V)}\,+\,\|\psi\|_{H^1(0,T;H)\cap L^\infty(0,T;V)}\,+\,
\|\rho\|_{L^\infty(0,T;H)\cap L^2(0,T;V)}}\non\\
&\newju{\quad \le\,c \,\|\bh\|_{L^2(0,\pcol{T};H)^2}\,\|\bk\|_{L^2(0,\pcol{T};H)^2}.}
\end{align}                                                                                          
\newju{Having shown this, we may bring $\beta\dt\psi$ to the \rhs\ of \eqref{bilin:2}. Elliptic regularity theory
and the estimates shown above then yield that also
}
\begin{align*}
\newju{\norma{\psi}_{\L2 {W_0}}\leq c \,\|\bh\|_{L^2(0,\pcol{T};H)^2}\,\|\bk\|_{L^2(0,\pcol{T};H)^2},}
\end{align*}
\newju{which concludes the proof of the assertion.}
\Edim

%\todo{
%\begin{remark}
%Add a remark concerning the differences w.r.t. J\"urgen's approach
%\end{remark}}
\juergen{We now provide the announced independent proof for the form of the second-order derivative.}  
\begin{theorem}
\label{THM:FRECHET:II}
Assume that \ref{const:weak}, \ref{F:strong:1}--\ref{Ph:strong}, \pcol{\eqref{defUR}--%
%\eqref{strong:initialdata}, and 
\eqref{strong:sep:initialdata}} are fulfilled, and let $\ov\bu\in{\cal U}_R$ be given. Then the
second Fr\'echet derivative $D^2\S(\ov \bu) \in {\cal L}(\Uh, {\cal L}(	\Uh,\juergen{\Y}))$
is \pcol{given, for every $\bh=(h_1,h_2),\bk=(k_1,k_2) \in \Uh$, by} the identity 
$\,D^2\S(\ov \bu)\juergen{(\bk)(\bh)}=(\nu,\psi,\rho)$,
where $(\nu,\psi,\rho)$ is the unique solution to the bilinearized system \juergen{\eqref{bilin:1}--\eqref{g2} 
introduced in} Theorem \ref{THM:BILIN}.
\end{theorem}

\begin{proof}
\juergen{By virtue of Theorem \elvis{4.4}, $D^2{\cal S}(\ov\bu)$ exists as an element of ${\cal L}(\Uh, {\cal L}(	\Uh,\juergen{\Y}))$.
Now, the embedding of ${\cal Y}$ in $\newju{\cal V}$ is continuous. Therefore, 
${\cal S}$ is also twice continuously Fr\'echet 
differentiable between $\Uh$ and $\newju{\cal V}$, and the expressions for the derivatives coincide. It thus suffices 
to work in the space
$\newju{\cal V}$. To  this end, we recall that $\Uh_R$ is open in $\Uh$. Hence, there is some $K>0$ such that
$\ov\bu+\bk\in\Uh_R$ whenever $\|\bk\|_{\Uh}\le K$. In the following, we always tacitly assume that the occurring 
increments $\bk$ satisfy this condition. }

To prove the claim, we proceed in a direct \pcol{way}, by showing that there \juergen{exist some $\eps>0$ and $\widehat C>0$
such that}
\begin{align}
	& \non
	\norma{D\S(\ov \bu + \bk) - D\S(\ov \bu) - D^2\S(\ov \bu)(\bk)}_{{\cal L}(\Uh, \newju{\cal V})}
	\\ & \quad
	= \sup_{\norma{\bh}_{\Uh}=1} \Big\|\big(D\S(\ov \bu + \bk) - D\S(\ov \bu) - D^2\S(\ov \bu)(\bk) \big)
	\anfirst{(\bh)}\Big\|_{\newju{\cal V}}
	\leq \widehat C\, \norma{\bk}_{\L2 H ^2}^{1+\eps}.
	\label{D2:formal}
\end{align}

\juergen{At this point, we introduce some additional notation: for arbitrary $\bh,\bk \in \Uh$, we define the linearized variables}
\begin{align*}
	(\eta^\bh,\xi^\bh,\theta^\bh):=D{\cal S}(\ov\bu)(\bh), \quad(\ov \eta^\bh,\ov \xi^\bh,\ov \theta^\bh) := D\S(\ov\bu + \bk)\anfirst{(\bh)}.
	\end{align*}
\juergen{Notice that \eqref{est:rem} implies that}
\begin{equation}\label{elvis}
\juergen{\|(\eta^\bh,\xi^\bh,\theta^\bh)\|_{\cal Z}\le c\,\|\bh\|_{L^2(0,T;H)^2}, \quad
\|(\ov \eta^\bh,\ov \xi^\bh,\ov \theta^\bh)\|_{\cal Z}\le c\,\|\bh\|_{L^2(0,T;H)^2}\,,}
\end{equation}	
\juergen{where, here and in the remainder of this proof, $c>0$ denote constants that may depend on the data, but
not on the choice of $\bk\in{\cal U}$ with $\ov\bu+\bk\in{\cal U}_R$.}
  	
Next, we fix some $\bh \in\Uh$ with $\norma{\bh}_{\Uh}=1$ and introduce the auxiliary variables
\begin{align*}
	\zeta = \ov \eta ^\bh - \eta^\bh - \nu, 
	\quad
	\phi = \ov \xi^\bh - \ov \xi - \psi,
	\quad
	\omega = \ov \th^\bh - \th^\bh - \rho,
\end{align*}
where $(\nu,\psi,\rho)$ stands for the unique solution to the bilinearized system as obtained from Theorem \ref{THM:BILIN}.
With this notation, we realize that \eqref{D2:formal} reduces to 
\begin{align}
	\label{fre:formal:II}
	\norma{(\zeta,\phi,\omega)}_{\newju{\cal V}}^2
	\leq \widehat C\,\norma{\bk}_{\L2 H ^2}^{2+\eps},
\end{align}
which is the estimate we are going to show for some $\eps>0$.
To this end, we first observe that these new variables solve the initial-boundary value problem
\begin{align}
\label{freII:1}
&\alpha\dt\zeta+\dt\phi-\Delta\zeta\,= \anold{\Lambda_1} + \Lambda_2 &&\mbox{ in \,$Q$},\\[1mm]
\label{freII:2}
&\beta\dt\phi-\Delta\phi-\zeta = \chi \omega  - F''(\bvp)\phi + \Lambda_3 &&\mbox{ in \,$Q$},
\\[1mm]
\label{freII:3}
&\dt\omega-\Delta\omega+\chi\Delta\phi\,=\,\anold{-\Lambda_1} &&\mbox{ in \,$Q$},\\[1mm]
\label{freII:4}
&\dn\zeta\,=\,\dn\phi\,=\,\dn\omega\,=\,0 &&\mbox{ on \,$\Sigma$},\\[1mm]
\label{freII:5}
&\zeta(0)\,=\, \phi(0)\,=\, \omega(0)\,=\,0  &&\mbox{ in }\,\Omega,
\end{align}
where 
\begin{align*}
%\juergen{g} & \juergen{=P(\bph)(\omega - \chi \phi -\zeta)+ P'(\bph)(\bs + \chi (1-\bph) - \bm)\phi},
%\\ 
	\Lambda_1 & = 
		\anold{P(\bph)(\omega - \chi \phi -\zeta)
		+ (P(\bph^\bk)-P(\bph))\big( (\ov\th^\bh -\th^\bh)- \chi (\ov\xi^\bh -\xi^\bh )- (\ov\eta^\bh-\eta^\bh)\big)}
		\\ 
		& \quad
		+ [P(\bph^\bk)-P(\bph)-P'(\bph)\xik](\th^\bh- \chi \xih - \eta^\bh)
		\\ 
		& \quad
		+ [P'(\bph^\bk)-P'(\bph)-P''(\bph)\xik](\bs + \chi (1-\bph) - \bm)\xih
		\\ & \quad
		+ \,P'(\bph)(\omega - \chi \phi - \zeta)\xih
		+ P'(\bph)(\bs + \chi (1-\bph) - \bm)\phi
		\\ & \quad 
		+ (P'(\bph^\bk)-P'(\bph))\big( (\ov\s^\bk -\bs)- \chi (\ov\ph^\bk -\bph )- (\ov\m^\bk-\bm)\big)\xih
		\\ & \quad 
		+ (P'(\bph^\bk)-P'(\bph))(\bs + \chi (1-\bph) - \bm)(\ov \xi^\bh - \xih)		
		\\ & \quad 
		+ P'(\bph)\big( (\ov\s^\bk -\bs)- \chi (\ov\ph^\bk -\bph )- (\ov\m^\bk-\bm)\big)(\ov \xi^\bh - \xih)
		\\ & \quad
		+ (P'(\bph^\bk)-P'(\bph))\big( (\ov\s^\bk -\bs)- \chi (\ov\ph^\bk -\bph )- (\ov\m^\bk-\bm)\big)(\ov \xi^\bh - \xih),
		\\ \Lambda_2 & = 
		- h_1 \anold{[}\h(\bph^\bk) - \h(\bph) - \h'(\bph)\xik\anold{]}
		- [\h'(\bph^\bk)- \h'(\bph)-\h''(\bph)\xik]\xih \ov u_1 
		\anold{- \h'(\bph)\phi \ov u_1}
		\\ & \quad
		%- \h'(\bph)\xih k_1
		-(\h'(\bph^\bk)-\h'(\bph))(\ov \xi^\bh- \xih)\ov u_1
		{-(\h'(\bph^\bk)-\h'(\bph)) \xih k_1}
		\\ & \quad
		- \h'(\bph)(\ov \xi^\bh - \xih)k_1
		-  (\h'(\bph^\bk)-\h'(\bph))(\ov \xi^\bh - \xih)k_1,
		\\ 
		\Lambda_3 & =
		- [F''(\bph^\bk) - F''(\bph) - F^{(3)}(\bph)\xik]\xih
		-(F''(\bph^\bk)-F''(\bph))(\ov \xi^ \bh- \xih),
\end{align*}
\pcol{and $(\bm, \bph, \bs)= \calS (\ov\bu) $, $(\bm^\bk, \bph^\bk, \bs^\bk)= \calS \an{(\ov\bu + \bk)} $.}

Before estimating in detail, let us recall that, owing to the continuous dependence results obtained in Theorem \ref{THM:STRONG} and \anfirst{in} Theorem \ref{THM:CONTDEP:LIN}, and recalling that $\norma{\bh}_{\Uh}=1$, we have
that
\begin{align*}
	& \norma{(\bm^\bk,\bph^\bk,\bs^\bk)-(\bm,\bph,\bs)}_{\Z}
	+ 	\norma{(\ov \eta^\bh,\ov \xi^\bh,\ov \th^\bh)-(\eta^\bh,\xih,\th^\bh)}_{\newju{\cal V}}\non
	\\ & \quad
	\leq c \,(1+ \norma{\bh}_{\L2 H ^2})\,\norma{\bk}_{\L2 H ^2}
	\leq c \,\norma{\bk}_{\L2 H ^2}.
\end{align*}
\newju{In particular,} \pcol{we point out that} 
\begin{equation}\label{yeah}
\newju{\|\bph^\bk-\bph\|_{L^\infty(0,T;V)}+\|\ov\xi^\bh-\xih\|_{L^\infty(0,T;V)}\,\le\,c \,\norma{\bk}_{\L2 H ^2}.}
\end{equation}
\newju{We will use this estimate in the following at several places without further reference.}

Moreover, let us recall that $(\bm,\bph,\bs)$ and $(\eta^\bh,\xih,\th^\bh)$, as solutions to \Statesys\ and \eqref{aux1}--\eqref{aux5} with $\lambda_1=\lambda_2=1,\lambda_3=\lambda_4=0$, enjoy the 
bounds provided by \eqref{ssbound1}--\eqref{ssbound3} and 
\begin{align*}
	\norma{(\eta^\bh,\xih,\th^\bh)}_{\Y} \leq c,
\end{align*} 
respectively. {Notice also that the estimate \eqref{normaVdiff} is valid for the arguments $\varphi_1=\ov\varphi^\bk,
\varphi_2=\ov\varphi$. }In addition, we also owe to Taylor's formula with integral remainder for \newju{$P^{(i)}$} and $\h^{(i)}$ for $i=0,1$, which yields that
\begin{align}
	\label{Taylor:1}
\newju{	P^{(i)}(\bph^\bk)-P^{(i)}(\bph)-P^{(i+1)}(\bph)\xik }&\newju{ = P^{(i+1)}(\bph) {(\bph^\bk - \bph - \xik)} + 
{{ Q}^{\bk}_{i+2}} (\bph^\bk-\bph)^2,}
	\\[0.5mm] 
		\label{Taylor:2}
	\h^{(i)}(\bph^\bk)- \h^{(i)}(\bph)-\h^{(i+1)}(\bph)\xik & = \h^{(i+1)}(\bph) {(\bph^\bk - \bph - \xik)} + {{ R}^{\bk}_{i+2}} (\bph^\bk-\bph)^2,
\end{align}
for $i=0,1$, where the remainders $\,{{ Q}^{\bk}_i},{{ R}^{\bk}_i}\,$ are uniformly bounded and defined by
\begin{align*}
\newju{Q^{\bk}_{i+2} = \!\int_0^1\!\! P^{(i+2)}(\bph+s (\bph^\bk-\bph)) (1-s) \,ds,}\quad
{{ R}^{\bk}_{i+2}} = \!\int_0^1\!\! \h^{(i+2)}(\bph+s (\bph^\bk-\bph)) (1-s) \,ds.
\end{align*}
Namely, there exists a positive constant $\anfirst{R^*}$ such that, for all $\bk\in\Uh$ with $\ov\bu+\bk\in{\cal U}_R$,
\begin{align*}
\newju{\|Q_2^\bk\|_{L^\infty(Q)}+
	\norma{{{ Q}^{\bk}_3} }_{L^\infty(Q)}}+\norma{{{ R}^{\bk}_2} }_{L^\infty(Q)} +\norma{{{ R}^{\bk}_3} }_{L^\infty(Q)} \leq \anfirst{R^*}.
\end{align*}
Moreover, owing to Remark \ref{REM:FRE}, we also have \anfirst{that}
\begin{align*}
	\norma{(\bm^\bk-\bm-\eta^\bk, \bph^\bk - \bph - \anfirst{\xik},\bs^\bk-\bs-\th^\bk )}_\Z \leq c \norma{\bk}_{\L2H\anfirst{^2}}^2.
\end{align*}\anold{%
Thus, we test \eqref{freII:1} by $\zeta$, \eqref{freII:2}, to which we add on both sides the term $\phi$, by $\dt \phi$, and 
\eqref{freII:3} by $\om$, add the resulting equations, and integrate over time 
and by parts to obtain that
\begin{align}
	\non 
	&\frac\a2 \IO2 \zeta 
	+ \I2 {\nabla \zeta}
	+ \b\I2 {\dt \phi}
	+ \frac 12 \norma{\phi(t)}^2_V
	+ \frac 12 \IO2 \om 
	+ \I2 {\nabla \om}
	\\  & \quad\non
	=  \intQt \Lambda_1(\zeta-\om) + \intQt\Lambda_2 \zeta
	+ \intQt (\chi \omega  - F''(\bvp)\phi + \Lambda_3 +\phi ) \dt \phi
	\\ & \qquad
	+ \chi \intQt \nabla \phi \cdot \nabla \om  =: \an{I}_1+\an{I}_2+\an{I}_3+\an{I}_4.
	\label{fre:proof}
\end{align}
Using the Young and H\"older inequalities, the above properties, the continuous embedding $V\subset \Lx6$, 
\an{\eqref{fre:cd:estimate}, the general property \eqref{normaVdiff} with $p=4$ and $p=6$,
and \eqref{yeah},} 
 \newju{as well as Taylor\anfirst{'s} formulae \anfirst{\eqref{Taylor:1}--\eqref{Taylor:2}}}, we infer that
\begin{align*}
	\an{I}_1 & \leq 
	c \!\intQt \!(|\om|^2+|\phi|^2+|\zeta|^2)
	+ c \iot \norma{\diff}_V^2 \big(  \norma{\ov\th^\bh -\th^\bh}_4^2  +  \norma{\ov\xi^\bh -\xi^\bh}_4^2 +  \norma{\ov\eta^\bh-\eta^\bh}_4^2 \big)\juergen{\,ds}
	\next
	+ c (\norma{\th^\bh }_{L^\infty(Q)}^2+\norma{\xi^\bh }_{L^\infty(Q)}^2+\norma{\eta^\bh }_{L^\infty(Q)}^2)\iot (\juergen{\norma{\bk}^4_{L^2(0,s;H)^2}} + \norma{\bph^\bk-\bph}^4_{4})\juergen{\,ds}
	\next
	+ c (\norma{\bs }_{L^\infty(Q)}^2+\norma{\bph}_{L^\infty(Q)}^2+\norma{\bm}_{L^\infty(Q)}^2+1)\norma{\xih}_{L^\infty(Q)}^2\!\!\iot\! (\juergen{\norma{\bk}^4_{L^2(0,s;H)^2}} + \norma{\bph^\bk-\bph}^4_{4})
	\juergen{\,ds}
	\next
	+	c (\norma{\xi^\bh }_{L^\infty(Q)}^2+1)\intQt (|\om|^2+|\phi|^2+|\zeta|^2)
	\next
	+ c \norma{\xih}_{L^\infty(Q)}^2 \iot \norma{\diff}_{V}^2 \big(  \norma{\ov\s^\bk -\bs}_4^2  +  \norma{\ov\ph^\bk -\bph}_4^2 +  \norma{\ov\m^\bk-\bm}_4^2 \big)\juergen{\,ds}
%\end{align*}
%\begin{align*}
\next
+ c(\norma{\bs }_{L^\infty(Q)}^2+\norma{\bph}_{L^\infty(Q)}^2+\norma{\bm}_{L^\infty(Q)}^2+1) \iot \norma{\diff}_V^2 \norma{\ov \xi^\bh - \xih}_4^2\juergen{\,ds}
	\next
	+ c \iot \big(  \norma{\ov\s^\bk -\bs}_4 ^2 +  \norma{\ov\ph^\bk -\bph}_4^2 +  \norma{\ov\m^\bk-\bm}_4^2 \big)\norma{\ov \xi^\bh - \xih}_4 ^2\juergen{\,ds}
	\next
	+ c \iot \norma{\diff}_V^2 \big(  \norma{\ov\s^\bk -\bs}_6^2 +  \norma{\ov\ph^\bk -\bph}_6^2 +  \norma{\ov\m^\bk-\bm}_6^2 \big)\norma{\ov \xi^\bh - \xih}_6 ^2\juergen{\,ds}
	\\ &\leq
	c \iot (\norma{\om}^2 + \norma{\phi}^2 + \norma{\z}^2)\pcol{\,ds}	+ c \,\norma{\bk}^4_{\L2 H^2} \,.
\end{align*}
\juergen{Similar computations show that also}
\begin{align*}
	\an{I}_2 & \leq
	c \I2 \z
	+c\norma{h_1}_{L^\infty(Q)}^2 \iot (\juergen{\norma{\bk}^4_{L^2(0,s;H)^2}} + \norma{\diff}_4^4)
		\juergen{\,ds}
	\next
	+ c \,\norma{\xih}_{L^\infty(Q)}^2\norma{\ov u_1}_{L^\infty(Q)}^2\iot  
	(\juergen{\norma{\bk}^4_{L^2(0,s;H)^2}} + \norma{\diff}_4^4)\juergen{\,ds}
	\next
	%+	c \I2 \z 
	\pcol{{}+ c\, \norma{\ov u_1}_{L^\infty(Q)}^2 \I2 \phi
	%\next
	+ c \,\norma{\ov u_1}_{L^\infty(Q)}^2\iot \norma{\diff}_{4}^2\norma{\ov \xi^\bh - \xih}_4^2
	\juergen{\,ds}}
	\next
	+ \d \iot \norma{\z}_{\pcol V}^2\juergen{\,ds}
	+ \cd \,\norma{\xih}^2_{L^\infty(Q)} \iot \norma{\diff}^2_{\pcol 4}\norma{k_1}^2\juergen{\,ds}
	\next
	+ \d \iot \norma{\z}_{\pcol V}^2\juergen{\,ds}
	+ \cd \iot \norma{\ov \xi^\bh - \xih}^2_4\norma{k_1}^2\juergen{\,ds}
	\next
	+ \d \iot \norma{\z}_{\pcol V}^2\juergen{\,ds}
	+ \cd \iot \norma{\diff}^2_{\pcol 6}\norma{\ov \xi^\bh - \xih}^2_6\norma{k_1}^2\juergen{\,ds}
	\\ & \leq
	\an{3}\d \I2 {\nabla\zeta}
	+ \cd \I2 \zeta
	+ c \intQt |\phi|^2 + c \norma{\bk}^4_{\L2 H^2},
\end{align*}
for a positive $\d$ yet to be determined.
\pcol{As for $\an{I}_3$, we use \eqref{ssbound3} and, for the first term in $\Lambda_3$, we can argue as in \eqref{Taylor:1}--\eqref{Taylor:2} and \elvis{the} following lines. Then,} we have that
\begin{align*}
	\an{I}_3 & \leq 
	\d \I2 {\dt\phi}
	+ \cd \intQt(|\om|^2 + |\phi|^2) 
	+ \cd \I2 {\Lambda_3 }
	\\ & \leq
	\d \I2 {\dt\phi}
	+ \cd \intQt(|\om|^2 + |\phi|^2) 
	+ c \norma{\xih}^2_{L^\infty(Q)} \iot  (\juergen{\norma{\bk}^4_{L^2(0,s;H)^2}} + 
	\norma{\diff}_4^4)\juergen{\,ds}
	\next
	+ \cd \iot  \norma{\diff}^2_V  \norma{\ov \xi ^\bh - \xih}^2_4\juergen{\,ds} 
	\\ & \leq
	\d \I2 {\dt\phi}
	+ \cd \intQt(|\om|^2 + |\phi|^2)
	+ \cd \,\norma{\bk}^4_{\L2 H^2},
\end{align*}
and, as for $\an{I}_4$, we \pcol{see that}
\begin{align*}
	\an{I}_4 \leq \frac 12 \I2 {\nabla \om}
	+ \frac {\chi^2}2 \I2 {\nabla \phi},
\end{align*}
where the first term can be absorbed on the \lhs\ of \eqref{fre:proof}.
\newju{Therefore, collecting the above estimates and choosing $\d>0$} small enough, we obtain from Gronwall's lemma that}
\an{\begin{align*}
&\newju{\|\zeta\|_{L^\infty(0,T;H)\cap L^2(0,T;V)}+\|\phi\|_{H^1(0,T;H)\cap L^\infty(0,T;V)}+\|\omega\|
_{L^\infty(0,T;H)\cap L^2(0,T;V)}}\\
&\quad\newju{\le\,c\,\|\bk\|_{L^2(0,T;H)^2}^2\,.}
\end{align*}}
\newju{With this estimate shown, we can move \pcol{the} term $\beta\dt\phi$ to the \rhs\ of \eqref{freII:2}, and it readily follows
from elliptic regularity theory that also}
\begin{equation*}
\newju{\|\phi\|_{L^2(0,T;W_0)}\,\le\,c\,\|\bk\|_{L^2(0,T;H)^2}^\an{2}\,.}
\end{equation*}
\newju{In conclusion, the estimate \eqref{fre:formal:II} is valid with $\varepsilon=2$. This finishes the proof of the assertion.}
\end{proof}

The last result we present in this section will be useful later on to handle the second-order sufficient optimality conditions and establishes the \Lip\ continuity of $D^2\S$ in a suitable sense.
\begin{theorem}
\label{THM:CONTDEP:BILIN}
Suppose that the conditions \ref{const:weak}, \ref{F:strong:1}--\newju{\ref{Ph:strong}}, \ref{ass:control:Uad} and \eqref{defUR} are fulfilled.
Moreover, let the initial data $(\m_0,\ph_0,\s_0)$ satisfy \an{\eqref{strong:initialdata}--\eqref{strong:sep:initialdata}},
and let $\ov \bu^i=(\uebar^i,\uzbar^i)\in \UR$ be controls with the associated states $(\bmu_i,\bvp_i,\bsigma_i)={\cal S}(\overline\bu^i)$, \anfirst{$i= 1,2 $.}
Furthermore, let $\bh,\bk \in \Uh$ be admissible increments with the corresponding linearized \an{variables} $(\eta_i^\bh,\xi_i^\bh,\theta_i^\bh)= D\S(\ov\bu_i)(\bh),(\eta_i^\bk,\xi_i^\bk,\theta_i^\bk)= D\S(\ov\bu_i)(\bk)$, and bilinearized variables $(\nu_i,\psi_i,\rho_i)= D^2\S(\ov\bu_i)(\bh)(\bk)$, $i=1,2$.
\newju{Then it holds that}
\begin{align}
	\non & \norma{\big(D^2\S(\ov \bu^1)-D^2\S(\ov \bu^2)\big)(\bh)(\bk)}_{\newju{\cal V}} 
	\\ & \quad
	\leq {K}\, \norma{\ov\bu^1-\ov\bu^2}_{\L2 H^2}\norma{\bh}_{\L2 H^2}\norma{\bk}_{\L2 H^2},
	\label{fre:cd:estimate:bil}
\end{align}
\newju{with a constant $K>0$ that does not depend on the choice of $\,\bh,\bk\in {\cal U}_R$.}
\end{theorem}
\begin{proof}
By virtue of Theorem \ref{THM:FRECHET:II}, \eqref{fre:cd:estimate:bil} directly follows once we show
the existence of some $c>0$ such that
\begin{align}
	\label{est:cd:bilin}
	\norma{(\nu_1,\psi_1,\r_1)-(\nu_2,\psi_2,\r_2)}_{\newju{\cal V}} \leq 	c \,
	\norma{\ov \bu^1-\ov \bu^2}_{\L2 H ^2}\norma{\bh}_{\L2 H ^2}\norma{\bk}_{\L2 H ^2},
\end{align}
which is the estimate we are going to check. To this end, we set
\begin{align*}
	 \bm &= \bm_1-\bm_2, \quad 
	\bph = \bph_1-\bph_2, \quad
	 \bs= \bs_1-\bs_2, \quad
	\\  \nu &=\nu_1-\nu_2, \quad 
	\psi = \psi_1-\psi_2, \quad
	 \r= \r_1-\r_2,
\end{align*}
and observe that the triple of differences $( \nu,\psi,\r)$ solves the system
\begin{align}
	\label{cd:fre:1:bil}
	& \a \dt  \nu + \dt  \psi - \Delta  \nu =  f_1 +  f_2 && \hbox{in $Q$,}
	\\
	\label{cd:fre:2:bil}
	& \b \dt \psi - \Delta  \psi - \nu  = \chi  \r +f_3 && \hbox{in $Q$,}
	\\
	\label{cd:fre:3:bil}
	& \dt  \r - \Delta  \r +\chi \Delta  \psi=  -  f_1 && \hbox{in $Q$,}
	\\
	\label{cd:fre:4:bil}
	& \dn  \nu =\dn  \psi= \dn  \r = 0 && \hbox{on $\Sigma$,}
	\\
	\label{cd:fre:5:bil}
	&  \nu(0)= \psi(0)=  \r(0)=0 && \hbox{in $\Omega$,}
\end{align}
with 
\begin{align*}
		 f_1 & = 
		(P(\bph_1)-P(\bph_2))(\r_1 - \chi \psi_1 -\nu_1)
		+ P(\bph_2) (\r - \chi \psi - \nu)
		\next
		+ (P'(\bph_1)-P'(\bph_2))\xi_1^\bk (\th_1^\bh - \chi \xi_1^\bh - \eta_1^\bh)
		+ P'(\bph_2)(\xi^\bk_1-\xi^\bk_2) (\th_1^\bh - \chi \xi_1^\bh - \eta_1^\bh)
		\next	
		+ P'(\bph_2)\xi^\bk_2 \big( (\th_1^\bh -\th_2^\bh) - \chi (\xi_1^\bh -\xi_2^\bh )- (\eta_1^\bh-\eta_2^\bh) \big)
		\next
		+ (P''(\bph_1)-P''(\bph_2))\xi_1^\bk\xi_1^\bh (\bs_1 + \chi (1-\bph_1) - \bm_1)	
		\next
		+ P''(\bph_2)(\xi_1^\bk-\xi_2^\bk)\xi_1^\bh (\bs_1 + \chi (1-\bph_1) - \bm_1)			
		\next
		+ P''(\bph_2)\xi_2^\bk(\xi_1^\bh -\xi_2^\bh )(\bs_1 + \chi (1-\bph_1) - \bm_1)				
		\next
		+ P''(\bph_2)\xi_2^\bk \xi_2^\bh(\bs - \chi \bph - \bm)	
		+ (P'(\bph_1)-P'(\bph_2))\psi_1 (\bs_1 + \chi (1-\bph_1) - \bm_1)	
		\next
		+ P'(\bph_2)\psi (\bs_1 + \chi (1-\bph_1) - \bm_1)	
		+ P'(\bph_2)\psi_2 (\bs - \chi \bph - \bm)	
		\next
		+ (P'(\bph_1)-P'(\bph_2))\xi^\bh_1 	(\th_1^\bk - \chi \xi_1^\bk - \eta_1^\bk)
		+ P'(\bph_2)(\xi^\bh_1-\xi^\bh_2) 	(\th_1^\bk - \chi \xi_1^\bk - \eta_1^\bk)
		\next
		+ P'(\bph_2)\xi^\bh_2 	\big( (\th_1^\bk-\th_2^\bk) - \chi (\xi_1^\bk-\xi_2^\bk) - (\eta_1^\bk-\eta_2^\bk)\big)
		,
%		\\ 
\end{align*}
\begin{align*}
		 f_2 & = 
		 - (\h''(\bph_1)-\h''(\bph_2)) \xi^\bk_1\xi^\bh_1 \ov u_1^1
		 - \h''(\bph_2)( \xi^\bk_1- \xi^\bk_2)\xi^\bh_1 \ov u_1^1
		 - \h''(\bph_2)\xi^\bk_2(\xi^\bh_1 -\xi^\bh_2)\ov u_1^1		 
		 \next
		 - \h''(\bph_2)\xi^\bk_2 \xi^\bh_2 (\ov u_1^1-\ov u_1^2)
		 - (\h'(\bph_1)-\h'(\bph_2))\xi^\bh_1 k_1
		 - \h'(\bph_2)(\xi^\bh_1 -\xi^\bh_2) k_1
		 \next
		 -  (\h'(\bph_1)-\h'(\bph_2)) \psi_1 \ov u_1^1
		 -  \h'(\bph_2) \psi \ov u_1^1
		 -  \h'(\bph_2) \psi_2 (\ov u_1^1-\ov u_1^2)
		 \next
		 -  (\h'(\bph_1)-\h'(\bph_2)) \xi^\bk_1 h_1 
		 -  \h'(\bph_2)( \xi^\bk_1- \xi^\bk_2) h_1 
		,
		\\ 
		 f_3 & = -
		 \pcol{{} (F''(\bph_1)-F''(\bph_2))\psi_1
		 - F''(\bph_2) \psi
		 -(F^{(3)}(\bph_1)-F^{(3)}(\bph_2))\xi^\bh_1\xi^\bk_1}	 
		 \next
		 -  \pcol{{}F^{(3)}(\bph_2)(\xi^\bh_1-\xi^\bh_2)\xi^\bk_1		
		 -F^{(3)}(\bph_2)\xi^\bh_2(\xi^\bk_1-\xi^\bk_2).}
\end{align*}
Moreover, due to \eqref{cont:dep:strong} and \eqref{fre:cd:estimate}, we have \eqref{cons} as well as
\begin{align}
	\label{est:cd:bil:proof}
	\non
	& \norma{\eta^\bh_1-\eta^\bh_2}_\newju{L^\infty(0,T;H) \cap L^2(0,T;V)}
	+ \norma{\xi^\bh_1-\xi^\bh_2}_{\H1 H \cap \L\infty V \cap \L2 {W_0}}
	\\[0.5mm]
	& 	+ \norma{\th^\bh_1-\th^\bh_2}_\newju{L^\infty(0,T;H)\cap L^2(0,T;V)}
	\,	\leq\,
	c \,\norma{\ov \bu^1-\ov \bu^2}_{\L2 H ^2}\norma{\bh}_{\L2 H ^2}
\end{align}
and the corresponding estimate for $\bk$. \newju{Also, owing to \pcol{Theorems~\an{\ref{THM:STRONG},} \ref{THM:FRECHET} and~\ref{THM:BILIN}}, 
for $i=1,2$ \pcol{it is clear that} $(\bmu_i,\bvp_i,\bsigma_i)$, $(\eta_i^\bh,\xi_i^\bh,\theta_i^\bh)$,
$(\eta_i^\bk,\xi_i^\bk,\theta_i^\bk)$,
 and also $(\nu_i,\psi_i,\rho_i)$ belong to $\,\Y$}. Moreover, for both $(\eta_i^\bh,\xi_i^\bh,\theta_i^\bh)$ and $(\eta_i^\bk,\xi_i^\bk,\theta_i^\bk)$, $i=1,2$,
we have \eqref{est:lin:cont} with the corresponding increment,
and from Theorem \ref{THM:BILIN}, that 
\begin{align}
	\label{est:lin:cont:bil}
	\norma{(\nu_i,\psi_i,\r_i)}_\newju{\cal V} \leq c \,\norma{\bh}_{\L2 H ^2}\norma{\bk}_{\L2 H ^2}, \quad \hbox{ $i=1,2$.}
\end{align}
We are now ready to proceed by arguing in a similar fashion as in Theorem \ref{THM:CONTDEP:LIN} \pcol{in order to check} \eqref{est:cd:bilin}.

\noindent
{\sc First estimate:}
To begin with, we multiply \eqref{cd:fre:1:bil} by $ \nu$, \eqref{cd:fre:2:bil}, to which we add to both sides the term $ \psi$, by $\dt  \psi$,
and \eqref{cd:fre:3:bil} by $ \r$. Then we add the resulting equalities and integrate $Q_t$ and by parts to obtain
\begin{align*}
	&\frac\a2 \IO2 { \nu }
	+ \I2 {\nabla  \nu }
	+ \b\I2 {\dt  \psi}
	+ \frac 12 \norma{ \psi(t)}^2_V
	+ \frac 12 \IO2 { \r }
	+ \I2 {\nabla  \r}
	\\  & \quad
	= \intQt  f_1 ( \nu -  \r) + \intQt  f_2  \nu 
	+ \intQt  ( \chi  \r +  \psi + f_3)\dt  \psi
	\next \quad 
	+ \chi \intQt \nabla  \psi \cdot \nabla  \r \,= : \,I_1+I_2+I_3+I_4.
\end{align*}
Then, using the Young and \Holder\ inequalities, the boundedness and the \Lip\ continuity of $P^{(i)}$ and $\h^{(i)}$ for $i=1,2$,
the continuous inclusion $V\subset \Lx p $ for $p \in [1,6]$, \eqref{normaVdiff}, along with the uniform bounds indicated above and the stability estimates \eqref{cons} and \eqref{est:cd:bil:proof}, we infer that
\begin{align*}
	I_1 & \leq 
	c \intQt (| \nu|^2+| \psi|^2+|  \r|^2)
%	\next
+	\newju{{}c\,\|\bph\|^2_{L^\infty(0,T;V)}\iot( \norma{\r_1}_4^2 +
	  \norma{\psi_1}_4^2+\norma{\nu_1}_4^2) \,ds}
	\next
	+ c \norma{\xi_1^\bk}_{\L\infty V}^2 (\norma{\th_1^\bh}_{\L\infty V}^2+\norma{\xi_1^\bh}_{\L\infty V}^2+\norma{\eta_1^\bh}_{\L\infty V}^2)  \iot \norma{\bph}^2_{\pcol 6} \,ds
	\next
	+ c (\norma{\th_1^\bh}_{\L\infty V}^2+\norma{\xi_1^\bh}_{\L\infty V}^2+\norma{\eta_1^\bh}_{\L\infty V}^2) \iot \norma{\xi^\bk_1-\xi^\bk_2}^2_4	\,ds
	\next
	+ c \norma{\xi_2^\bk}_{\L\infty V}^2\iot (\norma{\th_1^\bh -\th_2^\bh}^2_4 +  \norma{\xi^\bh_1-\xi^\bh_2}^2_4 +\norma{\eta_1^\bh-\eta_2^\bh}^2_4 )\,ds
	\next
	+ c  \norma{\xi_1^\bk}_{\L\infty V}^2 \norma{\xi_1^\bh}_{\L\infty V}^2( \norma{\bs_1}_{L^\infty(Q)}^2 +  \norma{\bph_1}_{L^\infty(Q)}^2+\norma{\bm_1}_{L^\infty(Q)}^2 +1) \iot \norma{\bph}^2_{\pcol 6} \,ds
	\next
	+ c  \norma{\xi_1^\bh}_{\L\infty V}^2( \norma{\bs_1}_{L^\infty(Q)}^2 +  \norma{\bph_1}_{L^\infty(Q)}^2+\norma{\bm_1}_{L^\infty(Q)}^2 +1)\iot \norma{\xi_1^\bk-\xi_2^\bk}^2_4\,ds
		\next
	+ c \norma{\xi_2^\bk}_{\L\infty V}^2( \norma{\bs_1}_{L^\infty(Q)}^2 +  \norma{\bph_1}_{L^\infty(Q)}^2+\norma{\bm_1}_{L^\infty(Q)}^2 +1)\iot \norma{\xi_1^\bh-\xi_2^\bh}^2_4\,ds
	\next
	+ c \norma{\xi_2^\bk}_{\L\infty V}^2 \norma{\xi_2^\bh}_{\L\infty V}^2 \iot (\norma{\bs}^2_6 + \norma{\bph}_6^2 + \norma{\bm}_6^2 ) \,ds
	\next
	+ c   \norma{\psi_1}_{\L\infty V}^2( \norma{\bs_1}_{L^\infty(Q)}^2 +  \norma{\bph_1}_{L^\infty(Q)}^2+\norma{\bm_1}_{L^\infty(Q)}^2 +1)   \iot \norma{\bph}^2_{\pcol 4} \,ds
		\next
	+ c  ( \norma{\bs_1}_{L^\infty(Q)}^2 +  \norma{\bph_1}_{L^\infty(Q)}^2+\norma{\bm_1}_{L^\infty(Q)}^2 +1)  \I2 \psi
			\next
	+ c \norma{\psi_2}_{\L\infty V}^2   \iot (\norma{\bs}^2_4 + \norma{\bph}_4^2 + \norma{\bm}_4^2 ) \,ds
	\next
	+ c \norma{\xi_1^\bh}_{\L\infty V}^2 (\norma{\th_1^\bk}_{\L\infty V}^2+\norma{\xi_1^\bk}_{\L\infty V}^2+\norma{\eta_1^\bk}_{\L\infty V}^2)  \iot \norma{\bph}^2_{\pcol 6} \,ds	
	\next
	+ c (\norma{\th_1^\bk}_{\L\infty V}^2+\norma{\xi_1^\bk}_{\L\infty V}^2+\norma{\eta_1^\bk}_{\L\infty V}^2) \iot \norma{\xi^\bh_1-\xi^\bh_2}_4^2	\,ds
	\next
	+ c \norma{\xi_2^\bh}_{\L\infty V}^2\iot (\norma{\th_1^\bk-\th_2^\bk}^2_4 +  \norma{\xi^\bk_1-\xi^\bk_2}_4^2 +\norma{\eta_1^\bk-\eta_2^\bk}_4^2 )\,ds
	\\ & \leq
	c \intQt (|\nu|^2+| \psi|^2+|  \r|^2) 
	 + c \norma{\ov \bu^1-\ov \bu^2}_{\L2 H ^2}^2\norma{\bh}_{\L2 H ^2}^2\norma{\bk}_{\L2 H ^2}^2.
\end{align*}
Besides, similar computations \pcol{on the second term allow us to find out} that, for any $\delta>0$,
\newpage
\begin{align*}
	I_2 & \leq 
	c  \intQt | \nu|^2
	+ c\, \norma{\ov u_1^1}^2_{L^\infty(Q)} \norma{\xi_1^\bk}^2_{\L\infty V} \norma{\xi_1^\bh}^2_{\L\infty V} \iot \norma{\bph}^2_{\pcol 6} \,ds
	%\end{align*}
	%\begin{align*}
	\next
	+ c  \norma{\ov u_1^1}^2_{L^\infty(Q)}\norma{\xi_1^\bh}^2_{\L\infty V} \iot \norma{\xi^\bk_1-\xi^\bk_2}^2_4	\,ds
	\next
	+ c  \norma{\ov u_1^1}^2_{L^\infty(Q)}\norma{\xi_2^\bk}^2_{\L\infty V} \iot \norma{\xi^\bh_1-\xi^\bh_2}^2_4	\,ds
	\next
	+ \d \iot \norma{\nu}^2_{\pcol V}\,ds
	+ \cd \norma{\xi_2^\bk}^2_{\L\infty V} \norma{\xi_2^\bh}^2_{\L\infty V}\iot \norma{\ov u_1^1-\ov u_1^2}^2\,ds
 	\next
 	+ \cd \norma{ \xi^\bh_1 }^2_{\L\infty V} \iot \norma{\bph}_{\pcol 6}^2
	\norma{k_1}^2\,ds
%	+ \d \iot \norma{\nu}^2_4\,ds
	+ \cd  \iot \norma{\xi^\bh_1- \xi^\bh_2}_4^2\norma{k_1}^2\,ds
    \next
 	+ c \norma{\ov u_1^1}^2_{L^\infty(Q)}\norma{\psi_1}^2_{\L\infty V} \iot \norma{\bph}^2_{\pcol 6} \,ds
	\next
 	+ c \norma{\ov u_1^1}^2_{L^\infty(Q)} \I2 \psi
 	+ \cd  \norma{\psi_2}^2_{\L\infty V}\iot \norma{\ov u_1^1-\ov u_1^2}^2\,ds
 	\next
 	%+ \d \iot \norma{\nu}^2_4\,ds
 	+ \cd \norma{ \xi^\bk_1 }^2_{\L\infty V} \iot \norma{\bph}_{\pcol 6}^2\norma{h_1}^2\,ds
 	+ \cd  \iot \norma{\xi^\bk_1- \xi^\bk_2}_4^2\norma{h_1}^2\,ds
	\\
	& \leq
	\d \I2 {\nabla \nu} 
	+ \cd \I2 \nu
	+ c \I2 \psi
	\next
	+ \cd \norma{\ov \bu^1-\ov \bu^2}_{\L2 H ^2}^2\norma{\bh}_{\L2 H ^2}^2 \norma{\bk}_{\L2 H ^2}^2.
\end{align*}
Employing the separation property \eqref{ssbound2}, which entails the \Lip\ continuity of $F$ and of its derivatives, and Young's inequality,
 we obtain that
\begin{align*}
	I_3 
	& \leq 
	\d \I2 {\dt  \psi}
	+ \cd \intQt (| \r|^2+| \psi|^2)
	+ \cd \pcol{{}\norma{\psi_1}^2_{\L\infty V} \iot \norma{\bph}^2_4} \,ds
	\next
	\pcol{{}+ \cd \norma{\xi_1^\bh}^2_{\L\infty V}\norma{\xi_1^\bk}^2_{\L\infty V} \iot \norma{\bph}^2_6 \,ds}
	\next
	+ \cd \norma{\xi_1^\bk}^2_{\L\infty V}\iot \norma{\xi_1^\bh-\xi_2^\bh}^{\pcol 2}_4\,ds
	+ \cd \norma{\xi_2^\bh}^2_{\L\infty V}\iot \norma{\xi_1^\bk-\xi_2^\bk}^{\pcol 2}_4\,ds
	\\ & 
	\leq
	\d \I2 {\dt  \psi}
	+ \cd \intQt (| \r|^2+| \psi|^2)
	+  \cd  \norma{\ov \bu^1-\ov \bu^2}_{\L2 H ^2}^2\norma{\bh}_{\L2 H ^2}^2 \norma{\bk}_{\L2 H ^2}^2.
\end{align*}
Finally, by Young's inequality \pcol{we deduce that}
\begin{align*}
	I_4 & \leq 
		\d \intQt | {\nabla  \r}|^2
		+ \cd \intQt  |\nabla  \psi|^2.
\end{align*}
Therefore, choosing $\d>0$ small enough, and invoking Gronwall's lemma, we conclude that
\begin{align*}
	& \norma{ \nu }_{\L\infty H \cap \L2 V}
	+ \norma{ \psi }_{\H1 H \cap \L\infty V}
	+ \norma{ \r}_{\L\infty H \cap \L2 V} \non
	\\[0.3mm] 
	&\quad 	\leq 
	 c \norma{\ov \bu^1-\ov \bu^2}_{\L2 H ^2}\norma{\bh}_{\L2 H ^2} \norma{\bk}_{\L2 H ^2}.
\end{align*}

\noindent
{\sc Second estimate:}
The above estimate entails that the \pcol{norm of
$\dt  \psi $ in $\L2 H$ is bounded} by the expression
on the \rhs. 
\newju{Thus, in equation \eqref{cd:fre:2:bil} the term $\dt  \psi $ can be absorbed on the \rhs. A straightforward 
computation, which may be left to the reader, shows that the entire \rhs\ is bounded in $\L2 H$ by the same expression. Hence, we can infer from the elliptic regularity theory that also} 
\begin{align*}
	\norma{ \psi}_{\L2 {W_0}}\leq 
	 c \norma{\ov \bu^1-\ov \bu^2}_{\L2 H ^2}\norma{\bh}_{\L2 H ^2} \norma{\bk}_{\L2 H ^2}.
\end{align*}
This concludes the proof of the assertion.
\end{proof}

%%%%%%%%%%%%%%%%%%%%%%%%%%%%
\section{First-order Necessary Optimality Conditions}
\label{SEC:FOC}
\setcounter{equation}{0}

%%%%%%%%%%%%%%%%%%%%%%%%%%%
We now derive first-order necessary optimality conditions.
\newju{By the well-known characterization for differentiable maps on convex sets, it holds
(see, e.g., \cite{Fredibuch}) that}
\begin{align}
	\label{foc:abst}
	\newju{D \Jred(\ov\bu)}(\bu-\ov\bu) \geq 0 \quad \forall \bu \in \uad,
\end{align}
where $D \Jred$ denotes the derivative of the {\it reduced cost} functional given by   
\begin{align}
	\Jred (\bu):= \J(\S(\bu),\bu).
	% = \J_1(\S(\bu),\bu) 
	%+ \kappa g(\bu).
\end{align}
Therefore, \newju{using Theorem \ref{THM:FRECHET} and the chain rule, we have the following result:}
\begin{theorem}
\label{THM:FOC}
Suppose that \ref{const:weak}, \ref{F:strong:1}--\juergen{\ref{Ph:strong}}, \ref{ass:control:const}--\ref{ass:control:Uad}, \juergen{as well as \eqref{strong:initialdata}--\eqref{strong:sep:initialdata},} 
are fulfilled.
Moreover, let $\ov\bu$ be an optimal control for ($\cal CP$) with corresponding state $(\bm,\bph,\bs)$.
Then \juergen{it holds the variational inequality}
\begin{align}
	& %\non 
	{\anfirst{b_1}} \intQ (\bph-\hat \ph_Q)\xi
	+ {\anfirst{b_2}}\iO (\bph(T)-\hat \ph_\Omega) \xi(T)
	+ {\anfirst{b_0}} \intQ  \ov\bu \cdot (\bu -\ov\bu)
%	\\ & \quad 
%	+ \kappa (g(\bu) -g(\ov\bu))
	\geq 
	0 \quad \forall \bu \in \uad,
	\label{foc:first}
\end{align}
where $(\eta,\xi,\th)$ denotes the unique solution to the linearized system obtained from Lemma \ref{LEM:FRE} with $\anfirst{\bf h}=\bu -\ov\bu$ and $\lambda_1=\lambda_2=1,\lambda_3=\lambda_4=0$.
\end{theorem}
\newju{As usual, we simplify \eqref{foc:first} by means of the adjoint state variables $(p,q,r)$, which are
defined as the solution triple to the adjoint system whose strong form is given by the backward-in-time parabolic system} 
\begin{align}
	&\non - \dt p - \b \dt q - \Delta q + \chi \Delta r  + F''(\bph) q + \h'(\bph)\ov u_1 p
	\\ & \hspace{0.5cm}
	\label{adj1}
	-P'(\bph)(\bs+\chi(1-\bph)-\bm)(p-r) + \chi P(\bph)(p-r)
	={\anfirst{b_1}} (\bph - \hat \ph_Q) \quad&&\mbox{in }\,Q\,,\\[1mm]
	\label{adj2}
	&-\a \dt p-\Delta p -q + P(\bph)(p-r)=0 \quad&&\mbox{in }\,Q\,,\\[1mm]
	\label{adj3}
	& - \dt r -\Delta r  - \chi q - P(\bph)(p-r) =0 \quad&&\mbox{in }\,Q\,,\\[1mm]
	\label{adj4}
	&\dn p=\dn q=\dn r=0 \quad&&\mbox{on }\,\Sigma\,,\\[1mm] 
	\label{adj5}
	&(p+\b q)(T)= {\anfirst{b_2}}(\bph(T)-\hat \ph_\Omega),\quad \a p(T)= 0,\quad r(T)=0\quad &&\mbox{in }\,\Omega\,.
\end{align}
\Accorpa\Adjsys {adj1} {adj5}
Let us point out that since the \anfirst{terminal condition $(p+\b q)(T)$} prescribe\anfirst{s} a final datum in $\Lx2$ (cf. \ref{ass:control:targets}),
it is clear that the first equation \eqref{adj1} has to be considered in a weak \anfirst{sense}.
Here is the \anfirst{corresponding well-posedness result.}
\begin{theorem}[Well-posedness of the adjoint system]
\label{THM:ADJ}
\anfirst{Suppose} that \ref{const:weak}, \ref{F:strong:1}--\juergen{\ref{Ph:strong}, 
\ref{ass:control:const}--\ref{ass:control:Uad}, and 
\eqref{strong:initialdata}--\eqref{strong:sep:initialdata}} hold. Then the adjoint system \Adjsys\ admits a unique solution $(p,q,r)$ in the sense that
\begin{align*}
	p+ \b q  & \in \H1 {V^*},
	\\	
	p & \in \H1 H \cap \L\infty V \cap \L2 {W_0} \an{\cap L^\infty(Q)},
	\\
	q & \in \L\infty H \cap \L2 V,
	\\
	r & \in \H1 H \cap \L\infty V \cap \L2 {W_0}\an{\cap L^\infty(Q)},
\end{align*}
\juergen{where $(p,q,r)$ satisfies}
\begin{align*}
	& -\< \dt (p +\b q),v>
	+ \iO \nabla q \cdot \nabla v
	- \chi \iO \nabla r  \cdot \nabla v
	+ \iO F''(\bph) q v
	+ \iO \h'(\bph)\ov u_1 p v
	\\ 
	& \quad
	-\iO P'(\bph)(\bs+\chi(1-\bph)-\bm)(p-r) v 
	+ \chi \iO P(\bph)(p-r)v
	= {\anfirst{b_1}} \iO  (\bph - \hat \ph_Q) v,
	\\
	& -\< \an{\a} \dt p,v>
	+ \iO \nabla p \cdot \nabla v 
	-\iO qv
	 +\iO  P(\bph)(p-r)v=0,
	\\
	& 
	- \<\dt r,v>
	+\iO \nabla r\cdot \nabla v
	- \chi \iO q v
	 -\iO  P(\bph)(p-r) v=0 ,
\end{align*}
for every $v\in V$\juergen{ and almost every $t \in (0,T)$, as well as }the terminal conditions
\begin{align*}
		(p+\b q)(T)= {\anfirst{b_2}}(\bph(T)-\hat \ph_\Omega),\quad \a p(T)= 0,\quad r(T)=0,\quad \mbox{\anfirst{$a.e.$} in }\,\Omega\,.
\end{align*}
\end{theorem}
\Brem
\label{REM:ADJ}
\pcol{Before entering the proof of the above theorem, let us point out that the regularity conditions on the solution imply that both $p+\b q$ and $q$ are in $\H1 {V^*}\cap \C0 H$,
whence all terminal conditions make sense. About that, we observe that the first condition may be rewritten just in terms of $q$ as $\b q(T)= {\anfirst{b_2}}(\bph(T)-\hat \ph_\Omega)$.}
\Erem
\begin{proof}
{For brevity, we again argue formally, thus avoiding the introduction of approximation schemes like in the proof of Lemma 4.1 and just providing the relevant} a priori estimates.
Moreover, let us notice that the adjoint system \Adjsys\ is linear, so that the uniqueness part also follows from standard arguments as a consequence of the following estimates.

\noindent
{\sc First estimate:}
First, we add to both sides of \eqref{adj2} the term $p$.
We then multiply \eqref{adj1} by $q$, the new \eqref{adj2} by $\anfirst{-}\dt p$, and 
\eqref{adj3} by $ \chi^2 r$, add the resulting equalities, and integrate over $Q_t^T\pcol{{}= \Omega \times (t,T){}}$ and by parts to obtain that
\begin{align*}
	& \frac {\b}2 \IO2 q
	+ \IT2 {\nabla q}
	+ \a  \IT2 {\dt p}
	+   \frac {1}2 \norma{p(t)}_V^2
	+ \frac {\chi^2}2 \norma{r(t)}^2
	+ \chi^2 \IT2 {\nabla r}
	\\ & \quad
	=  
	\frac {\b}2 \norma{q(T)}^2
	+ {\anfirst{b_1}}\Qtt (\bph-\hat \ph_Q)q
	\anfirst{+}\chi \Qtt \nabla r \cdot \nabla q
	- \Qtt F''(\bph)|q|^2
	\\ & \qquad
	- \Qtt \h'(\bph)\ov u_1 pq
	+ \Qtt P'(\bph)(\bs+\chi(1-\bph)-\bm)(p-r) q
	- \chi \Qtt P(\bph)(p-r)q
	\\ & \qquad
	\anfirst{+}\Qtt P(\bph)(p-r)\dt p
	- \Qtt p \dt p
	+ \chi^3 \Qtt q r
	+ \chi^2 \Qtt P(\bph)(p-r)r
	=: \sum_{i=1}^{11}I_i.
\end{align*}

We estimate the terms on the \rhs\ individually.
The first summand is bounded by a constant, due to the terminal conditions \eqref{adj5}
and the assumptions \ref{ass:control:const}--\ref{ass:control:targets}.
For the other terms on the same line, we have  that
\begin{align*}
	I_2+I_3+I_4
	\leq 
	c \Qtt(|q|^2+1)
	+ \frac {\chi^2}2 \IT2 {\nabla r}
	+ \frac {1}2 \IT2 {\nabla q},
\end{align*}
by virtue of Young's inequality, the assumptions \ref{ass:control:const}--\ref{ass:control:targets},
and the separation property \eqref{ssbound2}.
Using the Young and \Holder\ inequalities, we bound the integrals in the second line of the \rhs\ by
\begin{align*}
	I_5+I_6+I_7
	& \leq 
	c \Qtt (|p|^2+|q|^2)
	+ c \int_t^T (\norma{\bs}_\infty + \norma{\bph}_\infty + \norma{\bm}_\infty + 1 )
	(\norma{p} + \norma{q})\norma{q}\juergen{\,ds}
	\\ & \quad 
	 + 	c \Qtt (|p|^2+|q|^2+|r|^2)	 \leq c \Qtt (|p|^2+|q|^2+|r|^2),
\end{align*}
where we also owe to the boundedness of \juergen{$P$, $P'$ and $\ov u_1$, and to the fact that}
$(\bmu,\bph,\bsigma)$\pcol{,} as a solution to \Statesys\ in the sense of Theorem \ref{THM:STRONG}\pcol{,}
\anfirst{satisfies \eqref{ssbound1}}.
Finally, the terms on the last line of the \rhs\ can be easily bounded by means of Young's inequality, namely, 
\begin{align*}
	I_8+I_9+I_{10}+I_{11}
	& \leq 
	\pcol{\frac \a 2\Qtt |\dt p|^2 
%	+ \cd \Qtt (|p|^2 + |r|^2)
	+c \Qtt (|p|^2+|q|^2+|r|^2).}
\end{align*}
%with a positive $\d$ yet to be determined.
\pcol{Now, we combine the above estimates %, adjust $\d>0$ small enough, 
and invoke Gronwall's lemma to infer that}
\begin{align*}
	\norma{p}_{\H1 H \cap \L\infty V}
	+ \norma{q}_{\L\infty H \cap \L2 V}
	+ \norma{r}_{\L\infty H \cap \L2 V}
	\leq c.
\end{align*}

\noindent
{\sc Second estimate:}
We can now rewrite \anfirst{equation} \eqref{adj3} as a parabolic equation with the source term 
$f_r\anfirst{:} = \chi q + P(\bph)(p-r)$, which is bounded in \pcol{$\L\infty H$} due to the above estimate. It is then a standard matter to infer that
\begin{align*}
	\norma{r}_{\H1 H \cap \L\infty V \cap \L2 {W_0}}
	\leq c.
\end{align*}
\pcol{In addition, since \elvis{$r(T)=0\in L^\infty(\Omega)$},
we can apply the regularity result~\cite[Thm.~7.1, p.~181]{LSU} to infer that $\norma{r}_{L^\infty (Q)}	\leq c$ as well.} 

\noindent
{\sc Third estimate:}
From equation \eqref{adj2} and the \pcol{parabolic regularity theory}, similarly we recover that
\begin{align*}
	\norma{p}_{\L2 {W_0}\pcol{{}\cap L^\infty (Q)}}
	\leq c.
\end{align*}

\noindent
{\sc Fourth estimate:}
Finally, comparison in equation \eqref{adj1}, along with the above estimates, produces that
\begin{align*}
	\norma{\dt (p+\b q)}_{\L2 {V^*}}\leq c,
\end{align*}
\pcol{and this further estimate} concludes the proof.
\end{proof}

It is then a standard matter to use the adjoint variables to simplify the first-order necessary conditions
obtained in Theorem \ref{THM:FOC}.
\begin{theorem}
Assume that \ref{const:weak}, \ref{F:strong:1}--\juergen{\ref{Ph:strong}, 
\ref{ass:control:const}--\ref{ass:control:Uad}, and 
\eqref{strong:initialdata}--\eqref{strong:sep:initialdata} are fulfilled,} 
and let $\ov\bu\in\uad$ be an optimal control for ($\cal CP$) with corresponding state $(\bm,\bph,\bs)$.
Then, setting
\begin{align*}
	{\bf d}(x,t) := \pcol{\big(\!-\h(\bph(x,t)) p(x,t) , r(x,t)\big)}, \quad \hbox{for $a.e. (x,t) \in Q,$ }
\end{align*}
we have that                                           
\begin{align}
	& \intQ ({\bf d} + {\anfirst{b_0}} \ov\bu )\cdot (\bu -\ov\bu)
	\geq 
	0 \quad \forall \bu \in \uad.
	\label{foc:final}
\end{align}
\anfirst{Moreover, $\ov \bu$ is the $\L2 H^2$-orthogonal projection of $ - b_0^{-1} {\bf d}$ onto 
\juergen{$\uad$}.}
\end{theorem}

\begin{proof}
By comparing the variational inequalities \eqref{foc:first} and \eqref{foc:final}\anfirst{,} it follows that
\pcol{it suffices} to show \anfirst{the identity}
\begin{align}
	\label{simplification}
	- \intQ \h(\bph)\juergen{h_1} p + \intQ \juergen{h_2} r
	= {\anfirst{b_1}}\intQ (\bph-\hat \ph_Q)\xi
	+ {\anfirst{b_2}}\iO (\bph(T)-\hat \ph_\Omega) \xi(T),
\end{align}
with the choices $\juergen{h_1} = u_1-\ov u_1$ and $\juergen{h_2}= u_2-\ov u_2$,
where $(\eta,\xi,\th)$ and $(p,q,r)$ denote the unique solutions to the linearized system and to the adjoint system obtained from Lemma \ref{LEM:FRE} \an{with $\lambda_1 = \lambda_2 = 1$ and $\lambda_3 = \lambda_4 = 0$,} and Theorem
\ref{THM:ADJ}, respectively.

In this direction, we multiply the linearized system, written with $(\m,\ph,\s)=(\eta,\xi,\th)$, by $p,q$ and $r$, in this order. Then we add the obtained equalities and integrate over $Q$ \juergen{and by parts} to obtain that
\begin{align*}
	0 & =\juergen{\int_0^T\langle -\dt(p+\beta q),\xi\rangle\,dt \,+\int_Q \nabla\xi\cdot\nabla q}\\ 
	&\quad \juergen{+
	\intQ \xi \big[\chi \Delta r  + F''(\bph) q + \h'(\bph)\ov u_1 p} 
	\\ & \qquad\qquad
	-P'(\bph)(\bs+\chi(1-\bph)-\bm)(p-r) + \chi P(\bph)(p-r)\big]
	\\ & \quad
	+ \intQ \eta \big[-\a \dt p-\Delta p -q + P(\bph)(p-r)\big]
	\\ & \quad
	+\intQ \th \big[ - \dt r -\Delta r  - \chi q - P(\bph)(p-r) \big]
	\\ & \quad
	+ \iO [\a p(T)\pcol \eta(T) + (p+\b q)(T)\pcol\xi(T) + r(T) \pcol\th(T)]	
	\\ & \quad
	+ \intQ (\h(\bph)\juergen{h_1} p - \juergen{h_2} r).
\end{align*}
\juergen{Now, we use the weak form of \eqref{adj1}, as well as \eqref{adj2}--\eqref{adj5}, to deduce}
that the above \pcol{equality} reduces to \anfirst{\eqref{simplification},} which concludes the proof.
\end{proof}

\smallskip
\section{Second-order Sufficient Optimality Conditions}
\label{SEC:SOC}
\setcounter{equation}{0}
\newju{We now establish second-order sufficient optimality conditions. Since the control-to-state mapping $\S$ is 
only known to be Fr\'echet differentiable on ${\cal U}$, we are faced with the so-called ``two-norm discrepancy'' 
(see also \cite[Sec.~4.10.2]{Fredibuch}). In order to overcome this difficulty, we follow the approach taken in
\cite[Chap.~5]{Fredibuch}. Since many of the arguments developed here are rather similar to those employed in 
\cite{Fredibuch}, we can afford to be sketchy here. For full details, we refer the reader to \cite{Fredibuch}, noting that 
there the case of one control variable is treated while in our case we have to deal with a pair of controls.
In order to simplify the analysis somewhat, we now make an additional assumption.}

\an{
\begin{enumerate}[label={\bf (C\arabic{*})}, ref={\bf (C\arabic{*})}, start=4]
\item \qquad It holds $b_2=0$. \label{C4}
\end{enumerate}
}
\vspace{2mm}\noindent
\newju{Notice that under the  assumption \an{\ref{C4}}  we have a zero terminal condition for $p+\b q$ in \eqref{adj5}. 
This easily leads to the conclusion that we have the additional regularity $q\in Z$, which, in turn, means that the
adjoint system is satisfied in the strong form \eqref{adj1}--\eqref{adj5}.}

By virtue of \an{\ref{C4}} , we readily infer that
for every $((\m,\ph,\s), \bu)\in (\C0 H)^3\times \Uh$ 
and ${\bf v}=(v_1,v_2,v_3),{\bf w}=(w_1,w_2,w_3)$ such that
$ ({\bf v}, \bh),({\bf w}, \bk)\in (\C0 H)^3\times \Uh$ we have 
\begin{align}
	\label{D2:0}
	D^2\anold{\J}((\m,\ph,\s), \bu)\anfirst{(}({\bf v},\bh),({\bf w},\bk)\anfirst{)}
	= {\anfirst{b_1}}\intQ  v_2 w_2 
	%\sgn{+ b_2 \iO v_2(T) w_2(T)}
	+ {\anfirst{b_0}} \intQ \bh\cdot \bk.
\end{align}
Using Theorem \ref{THM:FRECHET:II} along with the above expression, we can now derive the second-order derivative of the reduced cost functional $\Jred$. 
Namely, for a fixed control $\ov \bu$ we find that
\begin{align}
	\non D^2\Jred (\ov \bu)\anfirst{(\bh,\bk)} & = D_{(\m,\ph,\s)}\J((\bm,\bph,\bs),\ov \bu)(\nu,\psi,\rho)
	\next \label{D2:1}
	+ D^2 \J ((\bm,\bph,\bs), \ov \bu)\anfirst{(}((\eta^\bh,\xih,\th^\bh),\bh),((\eta^\bk,\xik,\th^\bk),\bk)\anfirst{)},
\end{align}
where $(\eta^\bh,\xih,\th^\bh)$, $(\eta^\bk,\xik,\th^\bk)$, and $(\nu,\psi,\rho)$ stand for the unique corresponding solutions to the linearized \anfirst{system associated with \anold{$\bh$ and $\bk$},} and \anfirst{to the} bilinearized system, respectively.
From the definition of \anold{the cost functional \eqref{cost}} (recall that now ${\anfirst{b_2}}=0$), we \anold{readily infer}~that
\begin{align}
	\label{D2:2} D_{(\m,\ph,\s)}\J((\bm,\bph,\bs),\ov \bu)(\nu,\psi,\rho) 
	= {\anfirst{b_1}}   \intQ (\bph - \hat \ph_Q) \psi
	%\sgn{+ b_2\iO (\bph(T) - \hat \ph_{\Omega})\psi(T)}
	.
\end{align}
\newju{We now claim that} 
\anold{\begin{align}
 	\non
 	{\anfirst{b_1}} \int_Q (\bph-\hat \ph_Q)\psi 
 	%\sgn{+ b_2\iO (\bph(T) - \hat \ph_{\Omega})\psi(T)}
 	& = 
 	\int_Q [ P'(\bph)\xik(\th^\bh-\chi\xih-\eta^\bh)(p-r)
 	\next \quad \non
 	+ P''(\bph)\xik\xih(\bs +\chi (1-\bph)-\bm)(p-r)
 	\next \quad \non
 	+ P'(\bph)\xih(\th^\bk-\chi\xik-\eta^\bk)(p-r)
 	- \h''(\bph)\xik\xih \ov u_1 p
 	\next\quad 
 	 \label{second:simplification}
 	- \h'(\bph)\xih \newju{k_1 p-\h'(\bph)\xik h_1 p}
 	- F^{(3)}(\bph)\xih\xik q].
\end{align}}%
To prove this claim, we multiply \eqref{bilin:1} by $p$, \eqref{bilin:2} by $q$, \eqref{bilin:3} by $r$, add the resulting equalities, and integrate over $Q$, to obtain that
\begin{align*}
	0& = \intQ p [\alpha\dt\nu+\dt\psi-\Delta\nu- g_1 -g_2]
	\\ & \quad
	+ \intQ q [\beta\dt\psi-\Delta\psi-\nu-\chi \rho + F''(\bvp)\psi + F^{(3)}(\bph) \xih \xik ]
	\\ & \quad
	+ \intQ r [\dt\rho-\Delta\rho+\chi\Delta\psi+g_1],
\end{align*} 
\newju{with the functions $g_1$ and $g_2$ defined in} \eqref{g1}--\eqref{g2}. Then, we integrate by parts and make use of the initial and terminal conditions \eqref{bilin:5} and \eqref{adj5} to find that
\begin{align*}
0& = \intQ \nu [-\a \dt p-\Delta p -q + P(\bph)(p-r)]
	\\ & \quad
	+ \intQ \psi [- \dt p - \b \dt q - \Delta q + \chi \Delta r  + F''(\bph) q + \h'(\bph)\ov u_1 p
	\next \quad
	-P'(\bph)(\bs+\chi(1-\bph)-\bm)(p-r) + \chi P(\bph)(p-r)]
	\\ & \quad
	+ \intQ \rho [- \dt r -\Delta r  - \chi q - P(\bph)(p-r)]
	\\ & \quad
	+ \intQ \big[- P'(\bph)\xik(\th^\bh-\chi\xih-\eta^\bh)(p-r)
 	\next \non
 	\qquad\quad - P''(\bph)\xik\xih(\bs +\chi (1-\bph)-\bm)(p-r)
 	\next \non
 	\qquad\quad - P'(\bph)\xih(\th^\bk-\chi\xik-\eta^\bk)(p-r)
 	+ \h''(\bph)\xik\xih \ov u_1 p
 	\next
	\qquad\quad + \h'(\bph)\newju{\xih k_1 p + \h'(\bph)\xik h_1 p}
 	+ F^{(3)}(\bph)\xih\xik q \big],
\end{align*}
\newju{whence the claim follows by using the adjoint system.}
From this characterization, along with \eqref{D2:1} and \eqref{D2:2}, we conclude that
\begin{align}	
	\non &
	D^2 \Jred (\ov \bu) \anfirst{(\bh,\bh)} 
	= {\anfirst{b_0}} \norma{\bh}^2_{\L2 H ^2}
	+ \intQ 2 P'(\bph)\xih(\th^\bh-\eta^\bh)(p-r)
	- \intQ \newju{2}\,\h'(\bph)\xih h_1 p
	\next \non
	+ \intQ \anold{\Big(}{\anfirst{b_1}} - P''(\bph)(\bs +\chi (1-\bph)-\bm)(p-r)
	\anold{- 2 \chi P'(\bph)(p-r)}
	\next \quad
	+ \h''(\bph)\ov u_1 p+ F^{(3)}(\bph)q \anold{\Big)}|\xih|^2
%	\sgn{+ b_2\norma{\xih(T)}^2}
	.
	\label{D2RED}
\end{align}
This explicit expression for the second-order derivative of $\Jred$
allows us to establish  sufficient conditions for optimality of $\ov \bu$.
We aim at showing that, under suitable assumptions, $D^2 \Jred (\ov \bu) $ \anold{is a positive definite operator
on a suitable subset of $\L2 H ^2$, meaning that for any admissible increment $\bh$ it holds that
\begin{align}
	D^2 \Jred (\ov \bu)\anfirst{(\bh,\bh)} > 0.
	\label{second:abs}
\end{align}}%
However, \eqref{second:abs} is rather restrictive as we need such a condition just along some suitable directions. 
To this end, for every $\tau>0$, we introduce the \newju{sets} of {\it strongly active constraints},
\begin{align*}
	&\newju{A^1_\tau (\ov\bu):= \big\{ (x,t) \in Q: \,\, \big| -\h(\bph(x,t)) p(x,t) + {b_0} \ov u_1(x,t) \big|>\tau\big\},
	}\\[0.5mm]
	&\newju{A^2_\tau( \ov\bu):=\big\{(x,t)\in Q: \,\,\big| r(x,t) + {{b_0}} \ov u_2(x,t) \big| >\tau  \big\}.}
\end{align*}
Moreover, for any {increment $\bh=(h_1,h_2)\in\Uh$} we introduce the {componentwise conditions}
%\anold{\begin{align}
%	\label{condition}
%	\bh(x,t)
%	\begin{cases}
%	= (0,0)  &\hbox{if $(x,t)\in A_\tau(\ov\bu),$}
%	\\
%	s.t. \,\, \min\{h_1(x,t), h_2(x,t)\}\geq 0  &\hbox{if $\ov\bu(x,t)=\underline\bu(x,t)$ and $(x,t) \not \in A_\tau(\ov\bu)$,}
%	\\
%	s.t. \,\, \max\{h_1(x,t), h_2(x,t)\}\leq 0   &\hbox{if $\ov\bu(x,t)=\hat\bu(x,t)$ and $(x,t) \not \in A_\tau(\ov\bu)$,}
%	\end{cases}
%\end{align}}
{\begin{align}
	\label{condition:1}
	h_1(x,t)
	\begin{cases}
	= 0  &\hbox{if $(x,t)\in A^1_\tau(\ov\bu),$}
	\\
	\geq 0  &\hbox{if $\ov u_1(x,t)=\underline u_1(x,t)$ and $(x,t) \not \in A^1_\tau(\ov\bu)$,}
	\\
	\leq 0   &\hbox{if $\ov u_1(x,t)=\hat u_1(x,t)$ and $(x,t) \not \in A^1_\tau(\ov\bu)$,}
	\end{cases}
	\\
	\label{condition:2}
	h_2(x,t)
	\begin{cases}
	= 0 &\hbox{if $(x,t)\in A^2_\tau(\ov\bu),$}
	\\
	\geq 0  &\hbox{if $\ov u_2(x,t)=\underline u_2(x,t)$ and $(x,t) \not \in A^2_\tau(\ov\bu)$,}
	\\
	\leq 0   &\hbox{if $\ov u_2(x,t)=\hat u_2(x,t)$ and $(x,t) \not \in A^2_\tau(\ov\bu)$,}
	\end{cases}
\end{align}}
and define the associated {\it $\tau$-critical cone} by
\begin{align}%\label{taucrit}
	{\CC}_\tau (\ov\bu):= \{\bh=(h_1,h_2) \in \Uh : h_1 \,\,\hbox{satisfies \eqref{condition:1} and} \,\, h_2 \,\,
	\hbox{satisfies \eqref{condition:2} a.e. in} \,\, Q \}.
\end{align}

The second-order sufficient condition for optimality then reads as follows:
\begin{align}\label{hinreichend}
	\exists \, \d, \tau >0\, : \, 
	D^2\Jred(\ov\bu)\anfirst{(\bh,\bh)} \geq \d \norma{\bh}^2_{\L2 H^2} \quad \forall \,\bh \in \CC_\tau (\ov \bu),
\end{align}
where $D^2\Jred(\ov\bu)\anfirst{(\bh,\bh)} $ is given by \eqref{D2RED} with the choices $(\bm,\bph,\bs)= \S(\ov\bu)$,
$(\eta^\bh,\xih,\th^\bh)= D\S(\ov \bu)\anfirst{(\bh)}$, and the corresponding adjoint variables $(p,q,r)$.
We have the following result.
\begin{theorem}
	\label{THM:SECOND:SUFF}
Assume that \ref{const:weak}, \ref{F:strong:1}--\ref{Ph:strong}, \ref{ass:control:const}--\newju{\an{\ref{C4}},
as well as \an{\eqref{strong:initialdata}--\eqref{strong:sep:initialdata}}} are fulfilled. 
Let $\ov\bu\in \uad$ be an admissible control which satisfies
 \eqref{foc:first} with the corresponding state $(\bm,\bph,\bs)$
and adjoint variables $(p,q,r)$, as obtained from Theorem \ref{THM:STRONG} and \ref{THM:ADJ}, respectively.
Then there exist positive constants $\eps_1,\eps_2$ such that
\begin{align*}
	\Jred (\bu) \geq \Jred (\ov \bu) + \eps_1 \norma{\bu-\ov\bu}^2_{\L2 H^2}
	\quad \hbox{for every}
	\quad \bu \in \Uh \,\mbox{ such that }\, \norma{\bu-\ov\bu}_\Uh \leq \eps_2.
\end{align*}
In particular, it follows that $\ov\bu$ is locally optimal for $({\cal CP})$ in the sense of $\Uh$.
\end{theorem}

For the proof we follow the line of argumentation employed in the proof of \cite[Thm.~5.17]{Fredibuch}\newju{,
where in our case we deal with a system of parabolic equations and a pair of controls, 
\pcol{with state and control nonlinearly coupled}. However, the techniques used in 
\cite{Fredibuch} can straightforwardly be adapted to our
more complicated situation. We therefore merely sketch the arguments}.
\begin{proof}
Given an arbitrary $\bu \in \uad$, we infer from Taylor's theorem with integral remainder that
\begin{align*}
	\Jred(\bu) - \Jred (\ov\bu) = D\Jred(\ov\bu)\anfirst{(}{\bf v}\anfirst{)} + \frac12 D^2\Jred(\ov\bu)\anfirst{(}{\bf v},{\bf v}\anfirst{)} + 
	\elvis{R^{\Jred}(\bu,\ov \bu)},
\end{align*}
where we have set ${\bf v} = \bu - \ov \bu $\anfirst{, and where the remainder $R^{\Jred}$ is \anold{given by}}
\begin{align*}
	\elvis{R^{\Jred}(\bu,\ov \bu)}= \int_0^1 (1-s) \big(  D^2\Jred(\ov\bu+s {\bf v}) -D^2\Jred(\ov\bu)\big) \anfirst{({\bf v},{\bf v})} \, ds.
\end{align*}
To estimate $\, \anold{\big(}D^2\Jred(\ov\bu+s {\bf v}) -D^2\Jred(\ov\bu)\big) \anfirst{({\bf v},{\bf v})} $, we set
\begin{align*}
	(\m^s, \ph^s,\s^s) &= \S (\ov \bu + s {\bf v}),
	\quad
	(\eta, \xi,\th) = D\S (\ov \bu )\anfirst{(\bf v)},	
	\quad
	(\eta^s, \xi^s,\th^s) = D\S (\ov \bu + s {\bf v})\anfirst{(\bf v)},
	\\
	(\nu, \psi,\rho) &= D^2\S (\ov \bu) \anfirst{({\bf v},{\bf v})},
	\quad
	(\nu^s, \psi^s,\rho^s) = D^2\S (\ov \bu + s {\bf v}) \anfirst{({\bf v},{\bf v})}.
\end{align*}
\newju{By arguing along the lines of the proof of \cite[Thm.~5.17]{Fredibuch}, we find that it suffices to show~that}
\begin{align}
	\frac{\Big| R^{\Jred}(\bu,\ov\bu)\Big|}{\norma{\bu-\ov\bu}^2_{\L2H^2}}\anold{ \,\longrightarrow\,\,} 0 \quad \hbox{as} \quad \norma{\bu-\ov\bu}_\Uh \to 0.
	\label{claim:second}
\end{align}
Using \newju{\eqref{D2:2}}, we see that
\begin{align}
	\non
 & D_{(\m,\ph,\s)}\J((\m^s, \ph^s,\s^s), \ov\bu+s {\bf v})(\nu^s, \psi^s,\rho^s)
	- D_{(\m,\ph,\s)}\J((\bm, \bph,\bs), \ov\bu)(\nu, \psi,\rho)
	\next 
	= {\anfirst{b_1}} \intQ (\ph^s - \bph)\psi
	+ {\anfirst{b_1}} \intQ (\ph^s - \hat \ph_Q)(\psi^s- \psi)
%	\next  \quad 
%	\sgn{+b_2 \iO (\ph^s(T) - \bph(T))\psi(T)
%	+ b_2 \iO (\ph^s(T) - \hat \ph_\Omega)(\psi^s(T)- \psi(T))}
	=: I_1,
	\label{proof:est:1}
\end{align}
and \anfirst{\pcol{from} \eqref{D2:0}--\eqref{D2:1}} we obtain that
\begin{align}
	\non & 
	D^2 \J((\m^s, \ph^s,\s^s), \ov\bu+s {\bf v})\anfirst{(}((\eta^s, \xi^s,\th^s), {\bf v}),((\eta^s, \xi^s,\th^s), {\bf v})\anfirst{)}
	\next 
	%\non
	- 	D^2 \J((\m, \ph,\s), \ov\bu)\anfirst{(}((\eta, \xi,\th), {\bf v}),((\eta, \xi,\th), {\bf v})\anfirst{)}
	%\next 
	= {\anfirst{b_1}} \intQ (\xi^s+\xi)(\xi^s-\xi)=: I_2.
	\label{proof:est:2}
\end{align}
It \anfirst{then} readily follows from the Cauchy--Schwarz inequality, \an{and the} stability estimates
 \eqref{cont:dep:weak}, \eqref{fre:cd:estimate} \anold{and \eqref{fre:cd:estimate:bil}}, that
\begin{align*}
	 \an{I_1} &\leq {\anfirst{b_1}} (\norma{\ph^s - \bph}_{\L2 H}\norma{\psi}_{\L2 H}
	+ \norma{\ph^s - \hat \ph_Q}_{\L2 H}\norma{\psi^s- \psi}_{\L2 H})
	\\ & \leq c s \norma{{\bf v}}^3_{\L2 H^2}.
\end{align*}
Moreover, owing to \ref{ass:control:targets}, Theorem \ref{THM:FRECHET}, and to 
the stability estimate \eqref{fre:cd:estimate},
\begin{align*}
	\an{I_2} & \leq {\anfirst{b_1}} (\norma{\xi^s+\xi}_{\L2 H}\norma{\xi^s-\xi}_{\L2 H})
	\leq  c s \norma{{\bf v}}^3_{\L2 H^2}.
\end{align*}
We thus can conclude  that
\begin{align*}
	\Big| R^{\Jred}(\bu,\ov\bu)\Big| \leq c \int_0^1 (1-s) s \norma{{\bf v}}^3_{\L2 H^2}\,ds
	\leq c \norma{{\bf v}}_{\Uh}\norma{{\bf v}}^2_{\L2 H^2},
\end{align*}
so that \eqref{claim:second} directly follows. With this, the proof can be completed by adapting 
the argumentation of \cite{Fredibuch} correspondingly to our situation.
\end{proof}

%%%%%%%%%%%%%%%%%%%%%%%%%%%%%%%%%%%%%%%%%%%%%%%%%%%%%%%%%%%%%%%%%%%%%%%%

\pcol{\section*{Acknowledgments}
This research was supported by the Italian Ministry of Education, 
University and Research~(MIUR): Dipartimenti di Eccellenza Program (2018--2022) 
-- Dept.~of Mathematics ``F.~Casorati'', University of Pavia. 
In addition, {PC gratefully mentions} some other support 
from the GNAMPA (Gruppo Nazionale per l'Analisi Matematica, 
la Probabilit\`a e le loro Applicazioni) of INdAM (Isti\-tuto 
Nazionale di Alta Matematica) and points out his affiliation 
as Research Associate to the IMATI -- C.N.R. Pavia, Italy.}

%%%%%%%%%%%%%%%%%%%%%%%%%%%%%%%%%
%% bibliography
%%%%%%%%%%%%%%%%%%%%%%%%%%%%%%%%%

\End{document}

%%%%%%%%%%%%%%%%%%%%%%%%%%%%%%%%%%%%%%%%%%%%%
\begin{thebibliography}{10}

{\footnotesize

\pier{\bibitem{Barbu}
V. Barbu,
``Nonlinear Differential Equations of Monotone Type in Banach Spaces'',
Springer,
London, New York, 2010.}

\pier{\bibitem{Brezis}
H. Brezis,
``Op\'erateurs maximaux monotones et semi-groupes de contractions
dans les espaces de Hilbert'',
North-Holland Math. Stud. Vol.
{\bf 5},
North-Holland,
Amsterdam,
1973.}

\bibitem{cartan1967}
H.~Cartan,
``Calcul diff\'erentiel. Formes diff\'erentielles''.
Hermann, Paris, 1967.

\bibitem{CRW} 
C. Cavaterra, E. Rocca and H. Wu,
Long-time Dynamics and Optimal Control of a Diffuse Interface Model for Tumor Growth. 
{\it Appl. Math. Optim.} (2019), https://doi.org/10.1007/s00245-019-09562-5.

\bibitem{CFGS}
P. Colli, M.H. Farshbaf-Shaker, G. Gilardi and J. Sprekels,
Second-order analysis of a boundary control problem for the viscous
Cahn--Hilliard equation with dynamic boundary condition,
{\it Ann. Acad. Rom. Sci. Ser. Math. Appl.} {\bf 7} (2015), 41--66.

\bibitem{CGH}
P. Colli, G. Gilardi and D. Hilhorst,
On a Cahn--Hilliard type phase field system related to tumor growth. {\em Discret. Cont. Dyn. Syst.} {\bf 35}
(2015), 2423--2442.                                          

\bibitem{CGRS1}
P. Colli, G. Gilardi, E. Rocca and J. Sprekels, 
Vanishing viscosities and error estimate for a Cahn--Hilliard
type phase field system related to tumor growth. {\em Nonlinear Anal. Real World Appl.} {\bf 26} (2015),
93--108.

\bibitem{CGRS2}
P. Colli, G. Gilardi, E. Rocca and J. Sprekels, 
Asymptotic analyses and   error estimates for a Cahn--Hilliard
type phase field system modelling tumor growth. {\em Discret. Contin. Dyn. Syst. Ser. S} {\bf 10} (2017),
37--54.

\bibitem{CGRS3}
P. Colli, G. Gilardi, E. Rocca and J. Sprekels,
 Optimal distributed control of a diffuse interface model of tumor growth. {\em
Nonlinearity} {\bf 30} (2017), 2518--2546.

\bibitem{CGS24}
P. Colli, G. Gilardi and J. Sprekels,
A distributed control problem for a fractional tumor growth model. {\em Mathematics} {\bf 7} (2019), 792.

\bibitem{CSS1}
P. Colli, A. Signori and J. Sprekels,
Optimal control of a phase field system
modelling tumor growth with chemotaxis and singular potentials. 
{\em Appl. Math. Optim.} (2019), https://doi.org/10.1007/s00245-019-09618-6.

\bibitem{CLLW}
V. Cristini, X. Li, J.S. Lowengrub and S.M. Wise,
Nonlinear simulations of solid tumor growth using a mixture model: invasion and branching.  {\em J. Math. Biol.} {\bf 58} (2009), 723--763.

\bibitem{CL}
V. Cristini and J. Lowengrub,
Multiscale Modeling of Cancer: An Integrated Experimental and Mathematical Modeling Approach. Cambridge University Press, 2010.

\bibitem{DFRGM}
M. Dai, E. Feireisl, E. Rocca, G. Schimperna and M.E. Schonbek,
Analysis of a diffuse interface model of multi-species tumor growth,
{\it Nonlinearity\/} {\bf  30} (2017), 1639--1658.

\bibitem{Dieu}
J. Dieudonn\'e, 
``Foundations of Modern Analysis''. Pure and Applied Mathematics, vol. 10,
Academic Press, New York, 1960.

\bibitem{EK_ADV}
M. Ebenbeck and P. Knopf,
Optimal control theory and advanced optimality conditions 
for a diffuse interface model of tumor growth,
preprint arXiv:1903.00333 [math.OC] (2019), 1--34.

\bibitem{EK}
M. Ebenbeck and P. Knopf,
Optimal medication for tumors modeled by a Cahn--Hilliard--Brinkman equation,
{\it Calc. Var. Partial Differential Equations} {\bf 58} (2019), 
https://doi.org/10.1007/s00526-019-1579-z.

\bibitem{EGAR}
M. Ebenbeck and H. Garcke,
Analysis of a Cahn--Hilliard--Brinkman model for tumour growth with chemotaxis,
{\it J. Differential Equations} \textbf{266} (2019), 5998--6036.

\bibitem{FGR}
S. Frigeri, M. Grasselli and E. Rocca,
On a diffuse interface model of tumor growth,
{\it  European J. Appl. Math.\/} {\bf 26} (2015), 215--243. 

\bibitem{FLR}
S. Frigeri, K.F. Lam and E. Rocca,
On a diffuse interface model for tumour growth with non-local interactions and degenerate mobilities. 
In {\it  Solvability, Regularity, and Optimal Control of Boundary Value Problems for PDEs},
P. Colli, A. Favini, E. Rocca, G. Schimperna, J. Sprekels (eds.),
{\it Springer INdAM Series,} {\bf 22}, Springer, Cham, 2017, pp. 217--254.

\bibitem{FLRS}
S. Frigeri, K.F. Lam, E. Rocca and G. Schimperna,
On a multi-species Cahn--Hilliard--Darcy tumor growth model with singular potentials,
{\it Commun. Math Sci.} {\bf  16} (2018), 821--856. 	

\bibitem{FLS}
S.~Frigeri, K.F.~Lam and A.~Signori,
Strong well-posedness and inverse identification problem of a non-local phase field tumor model with degenerate mobilities,
preprint arXiv:2004.04537 [math.AP] (2020), \pcol{1--41.}

\bibitem{GARL_1}
H. Garcke and K.F. Lam,
Well-posedness of a Cahn--Hilliard system modelling tumour
growth with chemotaxis and active transport,
{\it European. J. Appl. Math.} {\bf 28} (2017), 284--316.

\bibitem{GARL_3}
H. Garcke and K.F. Lam,
Global weak solutions and asymptotic limits 
of a Cahn--Hilliard--Darcy system modelling tumour growth,
{\it AIMS Mathematics} {\bf 1} (2016), 318--360.

\bibitem{GARL_2}
H. Garcke and K.F. Lam,
Analysis of a Cahn--Hilliard system with non--zero Dirichlet 
conditions modeling tumor growth with chemotaxis,
{\it Discrete Contin. Dyn. Syst.} {\bf 37} (2017), 4277--4308.

\bibitem{GARL_4}
H. Garcke and K.F. Lam,
On a Cahn--Hilliard--Darcy system for tumour growth 
with solution dependent source terms, 
in ``Trends on Applications of Mathematics to Mechanics'', 
E.~Rocca, U.~Stefanelli, L.~Truskinovski, A.~Visintin~(eds.), 
{\it Springer INdAM Series} {\bf 27}, Springer, Cham, 2018, pp. 243--264.

\bibitem{GAR}
H. Garcke, K.F. Lam, R. N\"urnberg and E. Sitka,
A multiphase Cahn--Hilliard--Darcy model for tumour growth with necrosis,
{\it Math. Models Methods Appl. Sci.} {\bf 28} (2018), 525--577.

\bibitem{GARLR}
H. Garcke, K.F. Lam and E. Rocca,
Optimal control of treatment time in a diffuse interface model of tumor growth,
{\it Appl. Math. Optim.} {\bf 78} (2018), {495--544}.

\bibitem{GLSS}
H. Garcke, K.\,F. Lam, E. Sitka and V. Styles,
 A Cahn--Hilliard--Darcy model for tumour
growth with chemotaxis and active transport, {\em Math. \pcol{Models} Methods Appl. Sci.} {\bf 26} (2016),
1095--1148.

\bibitem{GLS}
\anold{H.~Garcke, K.F.~Lam and A.~Signori,
On a phase field model of Cahn--Hilliard type for tumour growth with mechanical effects, 
{\it Nonlinear Anal. Real World Appl.}  {\bf 57} (2021), 103192, https://doi.org/10.1016/j.nonrwa.2020.103192.}
%Preprint arXiv:1912.01945 [math.AP], (2019).

\bibitem{HZO}
A. Hawkins-Daarud, K.\,G. van der Zee and J.\,T. Oden,
Numerical simulation of a thermodynamically
consistent four-species tumor growth model, {\em Int. J. Numer. Math. Biomed. Eng.} {\bf 28} (2011),
3--24.

\bibitem{HKNZ}
D. Hilhorst, J. Kampmann, T.N. Nguyen and K.G. van der Zee, 
Formal asymptotic limit of a diffuse-interface tumor-growth model, 
{\it Math. Models Methods Appl. Sci.} {\bf 25} (2015), 1011--1043.

\bibitem{KL}
C. Kahle and K. F. Lam,
Parameter identification via optimal control for a Cahn--Hilliard-chemotaxis system with a variable mobility,
{\it Appl. Math. Optim.} {\bf 82} (2020), 63--104.

\bibitem{LSU}
O.\,A. Lady\v{z}enskaja, V.\,A. Solonnikov and N.\,N. Uralceva,
``Linear and Quasilinear Equations of Parabolic
Type'', Mathematical Monographs, vol. 23, American Mathematical Society, Providence, Rhode Island, 1968.

\pier{\bibitem{L1}
J.-L. Lions,
``\'{E}quations diff\'{e}rentielles 
op\'{e}rationnelles et probl\`emes aux limites'',
Die Grundlehren der mathematischen Wissenschaften, Bd. 111,
Springer-Verlag, Berlin-G\"{o}ttingen-Heidelberg, 1961.}

\bibitem{SS}
L.~Scarpa and A.~Signori,
On a class of non-local phase-field models for tumor growth with possibly singular potentials, chemotaxis, and active transport,
preprint arXiv:2002.12702 [math.AP]  (2020), \an{1--40. }

\bibitem{S}
A. Signori,
Optimal distributed control of an extended model of tumor
growth with logarithmic potential,
{\it Appl. Math. Optim.} \elvis{{\bf 82} (2020), 517--549}.

\bibitem{S_DQ}
A. Signori,
Optimality conditions for an extended tumor growth model with 
double obstacle potential via deep quench approach, 
{\it Evol. Equ. Control Theory} \pcol{{\bf 9} (2020) 193--217}.

\bibitem{S_b}
A. Signori,
Optimal treatment for a phase field system of Cahn--Hilliard 
type modeling tumor growth by asymptotic scheme, 
{\it Math. Control Relat. Fields} \pcol{{\bf 10}  (2020) 305--331.}

\bibitem{S_a}
A. Signori,
Vanishing parameter for an optimal control problem modeling tumor growth,
{\it Asymptot. Anal.} 	\pcol{{\bf 117} (2020) 43--66.}

\bibitem{SigTime} A.~Signori, Penalisation of long treatment time and optimal control of a tumour growth model of Cahn--Hilliard type with singular potential, Preprint arXiv:1906.03460 [math.AP], 1--25.

\pier{\bibitem{Simon}
J. Simon,
{Compact sets in the space $L^p(0,T; B)$},
{\it Ann. Mat. Pura Appl.~(4)\/} 
{\bf 146} (1987) 65--96.}

\bibitem{ST}
J. Sprekels and F. Tr\"oltzsch, 
Sparse optimal control of a phase field system with singular potentials arising
in the modeling of tumor growth, preprint arXiv:2005.02784v3 [math.OC] (2020), \pcol{1--28.}

\bibitem{SW}
J. Sprekels and H. Wu,
Optimal distributed control of a Cahn--Hilliard--Darcy system with mass sources,
{\it Appl. Math. Optim.} (2019), https://doi.org/10.1007/s00245-019-09555-4.

\bibitem{Fredibuch}
F. Tr\"oltzsch, 
``Optimal Control of Partial Differential Equations: Theory, Methods and Applications'',
Graduate Studies in Mathematics vol. 112, American Mathematical Society, Providence, Rhode Island, 2010. 

\bibitem{WLFC}
S.M. Wise, J.S. Lowengrub, H.B. Frieboes and V. Cristini,
Three-dimensional multispecies nonlinear tumor growth I: Model and numerical method. 
{\it J. Theor. Biol.} {\bf 253(3)} (2008) 524--543.

}%

\end{thebibliography}
